\pdfoutput=1
\documentclass[11pt]{amsart}
\usepackage[utf8]{inputenc}
\usepackage[T1]{fontenc}
\usepackage{lmodern}
\ifdefined\pdfgentounicode\input{glyphtounicode}\pdfgentounicode=1\fi
\usepackage{amsmath,amssymb,amsthm}
\usepackage[margin=2.8cm]{geometry}
\usepackage{enumitem}
\usepackage{booktabs}
\usepackage{needspace}
\usepackage{xcolor}
\usepackage[colorlinks=true,linkcolor=blue!50!black,citecolor=blue!50!black,urlcolor=blue!50!black]{hyperref}

\newtheorem{theorem}{Theorem}[section]
\newtheorem{proposition}[theorem]{Proposition}
\newtheorem{lemma}[theorem]{Lemma}
\newtheorem{corollary}[theorem]{Corollary}
\newtheorem*{theoremA}{Theorem A}
\newtheorem*{theoremB}{Theorem B}
\newtheorem*{theoremC}{Theorem C}
\newtheorem*{theoremD}{Theorem D}
\newtheorem*{theoremE}{Theorem E}
\newtheorem*{theoremF}{Theorem F}
\newtheorem*{theoremG}{Theorem G}
\theoremstyle{definition}
\newtheorem{definition}[theorem]{Definition}
\newtheorem{example}[theorem]{Example}
\theoremstyle{remark}
\newtheorem{remark}[theorem]{Remark}

\newcommand{\PP}{\mathbb{P}}
\newcommand{\QQ}{\mathbb{Q}}
\newcommand{\RR}{\mathbb{R}}
\newcommand{\ZZ}{\mathbb{Z}}
\newcommand{\OO}{\mathcal{O}}
\newcommand{\cE}{\mathcal{E}}
\newcommand{\cF}{\mathcal{F}}
\newcommand{\cG}{\mathcal{G}}
\newcommand{\cI}{\mathcal{I}}
\newcommand{\cQ}{\mathcal{Q}}
\newcommand{\cU}{\mathcal{U}}
\newcommand{\cT}{\mathcal{T}}
\newcommand{\Tt}{\mathsf{T}}
\newcommand{\CC}{\mathcal{C}}
\newcommand{\BS}{\mathrm{BS}}
\DeclareMathOperator{\rk}{rk}
\DeclareMathOperator{\gr}{gr}
\DeclareMathOperator{\Ext}{Ext}
\DeclareMathOperator{\coker}{coker}

\title[Boij--Söderberg and Harder--Narasimhan polygons]{Boij--Söderberg and Harder--Narasimhan polygons of vector bundles on projective space}
\author{Cristian Anghel}
\address{``Simion Stoilow'' Institute of Mathematics of the Romanian Academy, Bucharest, Romania}
\email{cristian.anghel@imar.ro}
\author{Filip Chindea}
\address{``Simion Stoilow'' Institute of Mathematics of the Romanian Academy, Bucharest, Romania}
\email{fchindea@imar.ro}

\begin{document}

\begin{abstract}
Eisenbud and Schreyer proved that the cohomology table of a vector bundle $F$ on $\PP^n$ is a unique positive combination of supernatural tables along a chain of root sequences. The slopes of the pieces, weighted by their ranks, form a concave polygon $\Psi_F$, the Boij--Söderberg polygon, with the same endpoints as the Harder--Narasimhan polygon $P_F$. We prove that $P_F\le\Psi_F$, in any characteristic, by testing subsheaves against sheaves with pure resolutions on thickened lines; the bound is optimal among bounds depending only on the ray of the table. The chain majorizes every other supernatural decomposition, and $\Psi$ is monotone under cancellations, hence upper semicontinuous in flat families. The comparison extends to coherent sheaves and, by linear projection, to sheaves of positive dimension. On $\PP^2$ the defect is a sum of four cohomological terms; we characterize equality and show that the inequality survives the large volume limit of Bridgeland stability. For a smooth projective variety $X$ of dimension $n\ge2$ in characteristic zero, with $H$ ample and base-point-free, the mean normalized squared spread $G(\cE)$ of the roots in the cohomology table of every vector bundle $\cE$ satisfies $G(\cE)\ge(K_X^2-2c_2(X))\cdot H^{n-2}/((n+1)H^n)$, a bound independent of $\cE$. Since $G(\cE)=1$ exactly for arithmetically Cohen--Macaulay bundles, when this bound exceeds $1$ it is a uniform form of their absence.
\end{abstract}

\hypersetup{pdftitle={Boij--Söderberg and Harder--Narasimhan polygons of vector bundles on projective space},pdfauthor={Cristian Anghel and Filip Chindea}}
\maketitle

\enlargethispage{2pt}
\setcounter{tocdepth}{1}
\tableofcontents

\section{Introduction}\label{sec:intro}

Throughout, $K$ is an algebraically closed field of characteristic zero; the results that hold over any algebraically closed field are listed at the end of this introduction. For a coherent sheaf $\cF$ on $\PP^n$ let $\gamma(\cF)=(h^j(\cF(t)))_{0\le j\le n,\,t\in\ZZ}$ be its cohomology table. Eisenbud and Schreyer \cite{ES09} proved that the cohomology tables of vector bundles span, up to positive scalars, the cone generated by the tables $\gamma^z$ of supernatural bundles, and that every such table has a unique decomposition
\[
\gamma(F)=\sum_kq_k\,\gamma^{z^k},\qquad q_k>0,
\]
along a chain $z^1>z^2>\dots$ of root sequences; see \cite{ES10,ESi,Fl12} for the theory and its extension to coherent sheaves. A root sequence is a strictly decreasing sequence of integers $z=(z_1>\dots>z_n)$, and we normalize $\gamma^z$ to have Hilbert polynomial $\prod_j(t-z_j)$, i.e.\ rank $n!$. We attach to $z$ its \emph{center} $c(z)=\frac1n\sum_jz_j$ and its \emph{slope} $\mu(z)=-c(z)-\frac{n+1}2$; a supernatural bundle with root sequence $z$ has slope $\mu(z)$. Decreasing slopes correspond to increasing centers. Erman and Sam \cite{ErSa16} single out the question whether these numerical decompositions, whose coefficients are rational in general, can be categorified at the level of vector bundles, and they answer it in many cases. Here we take a different route and link the Eisenbud--Schreyer chain decomposition to slope stability, through the following polygon; on $\PP^2$ we determine exactly when the Harder--Narasimhan filtration realizes the chain decomposition (Corollary~\ref{cor:equalityP2}).

\begin{definition}
The \emph{Boij--Söderberg polygon} $\Psi_F\colon[0,\rk F]\to\RR$ of a vector bundle $F$ on $\PP^n$ is the concave polygon whose edges have slopes $\mu(z^k)$ and horizontal lengths $n!\,q_k$, arranged in decreasing order of slope.
\end{definition}

Its endpoints are $(0,0)$ and $(\rk F,\deg F)$, the same as those of the Harder--Narasimhan polygon $P_F$, whose edges have the slopes of the Harder--Narasimhan factors and lengths their ranks. The two polygons record different kinds of information: $P_F$ is defined by subsheaves, $\Psi_F$ by cohomology. Our main result compares them.

\begin{theoremA}[Theorem~\ref{thm:main}]
Let $F$ be a vector bundle on $\PP^n$, over an algebraically closed field of any characteristic. Then $\deg F'\le\Psi_F(\rk F')$ for every subsheaf $F'\subset F$; equivalently, $P_F\le\Psi_F$. In particular, a vector bundle whose cohomology table is a multiple of a single supernatural table is semistable.
\end{theoremA}

Theorem~A gives an upper bound for the slopes of all subsheaves in terms of the table alone. It is optimal among bounds that depend only on the ray of the table: every table has a multiple that is the table of a direct sum of supernatural bundles with $P=\Psi$. For a fixed table, integrality gives a sharper bound (Remark~\ref{rem:main}(1)). Moreover a saturated subsheaf of a supernatural bundle with the same slope, and its quotient, are again supernatural with the same roots (Proposition~\ref{prop:sameslope}), so supernatural bundles are Gieseker semistable. Theorem~A yields $V_{HN}(F)\le V_c(F)$ for the variances of the Harder--Narasimhan slopes and of the Boij--Söderberg slopes (Corollary~\ref{cor:supernatural}); equivalently, for $n\ge2$, the \emph{gap moment} $G(F)$, the mean over the chain of a normalized spread of the roots (Section~\ref{sec:gap}), satisfies $G(F)\ge1+\frac{12}{n+1}\sum_iw_i\,\delta(\cQ_i)$, where $\cQ_i$ are the Harder--Narasimhan factors, with weights $w_i=\rk\cQ_i/\rk F$ and normalized discriminants $\delta(\cQ_i)\ge0$ (Corollary~\ref{cor:Q}). Always $G(F)\ge1$, with equality exactly for split bundles.

The proof rests on a formula for the concave conjugate
\[
\Lambda_F(\lambda)=\sum_kn!\,q_k\big(\mu(z^k)-\lambda\big)_+ ,\qquad \Psi_F(m)=\min_\lambda\big(m\lambda+\Lambda_F(\lambda)\big).
\]
For a root sequence $z$ comparable with every element of the chain, $\Lambda_F(\mu(z))$ is the divided difference of truncated Euler characteristics
\[
\phi_z(\gamma(F))=(n-1)!\sum_{r=1}^n\frac{\sum_{j<r}(-1)^jh^j(F(z_r))}{\prod_{i\ne r}(z_r-z_i)} .
\]
By the theorem of Eisenbud, Fløystad, Weyman and Eisenbud, Schreyer on pure resolutions \cite{EFW11,ES09}, there is a sheaf $\Tt$, set-theoretically supported on a general line $L$ (in general on a thickening of $L$), whose minimal free resolution is pure, with twists $z_1,\dots,z_n$. For a saturated subsheaf $F'\subset F$ we may choose $L$ disjoint from the loci where $F'$ and $F/F'$ are not locally free; then $F'\otimes \Tt\subset F\otimes \Tt$, and
\[
\chi(F'\otimes \Tt)\ \le\ h^0(F'\otimes \Tt)\ \le\ h^0(F\otimes \Tt).
\]
Riemann--Roch turns the left-hand side into a positive multiple of $\deg F'-\mu(z)\rk F'$, and the hypercohomology spectral sequence of $F$ tensored with the resolution bounds the right-hand side by the same multiple of $\phi_z(\gamma(F))$. This spectral sequence argument is the one by which Eisenbud and Schreyer prove the nonnegativity of their pairing between Betti tables and cohomology tables \cite[Thm.~3.1]{ES09}, and the bound by $h^0(F\otimes\Tt)$ is the sheaf form of an inequality in their proof; here the sheaf $\Tt$ has one-dimensional support, and $h^0(F\otimes \Tt)$ replaces the vanishing of $F\otimes \Tt$. In the notation of \cite{ES09}, $\phi_z$ is a positive multiple of the pairing with the Betti table of $\Tt$ (see the paragraph after Lemma~\ref{lem:truncation}); for a categorified form of the Eisenbud--Schreyer functionals, as a single pairing between derived categories, see \cite{EE17}. On $\PP^2$ the test sheaf can be taken to be $\OO_{(a-b)L}(a)$, with resolution $\OO(b)\to\OO(a)$, for a line $L$ and $a>b$, and the argument is elementary (Section~\ref{sec:plane}).

Theorem~A has consequences in three directions. On $\PP^n$ it gives structural properties of the cone of cohomology tables: the chain decomposition is extremal among all decompositions into supernatural tables, and $\Psi$ is monotone under the cancellations that occur in extensions and in flat families (Theorems~B and~C). On $\PP^2$ the defect in Theorem~A is an explicit sum of four cohomological terms, the equality $P_F=\Psi_F$ is characterized, and the comparison persists in the large volume limit of Bridgeland stability (Theorems~D and~E). Finally, a projection to $\PP^n$ carries the comparison to sheaves of lower dimension on a larger projective space (Theorem~F) and to vector bundles on a polarized variety (Theorem~G). The last case was our starting point: after \cite{An1}, Eisenbud and Schreyer asked about the Boij--Söderberg cone of cohomology tables of vector bundles on the surfaces constructed there (Section~\ref{sec:twists}), and Theorem~G gives a lower bound, the same for all vector bundles on $(X,H)$, for the gap moment of their tables.

\medskip\noindent\textbf{Structural results on \texorpdfstring{$\PP^n$}{Pn}.}
Let $\CC$ be the set of nonnegative rational combinations of supernatural tables on $\PP^n$. Every $\gamma\in\CC$ has a unique chain decomposition, and $\Psi_\gamma$, $\Lambda_\gamma$ are defined as above. Applied to supernatural bundles, the inequality behind Theorem~A says that $\phi_z(\gamma^y)\ge n!\,(c(z)-c(y))_+$ for all root sequences $y,z$, with equality if $y$ and $z$ are comparable. So, for $\lambda\in\frac1n\ZZ$, $\Lambda_\gamma(\lambda)$ is the minimum of the linear functionals $\phi_z$ with $\mu(z)=\lambda$, attained at the root sequences comparable with the chain (Proposition~\ref{prop:minformula}). Two consequences follow.

\begin{theoremB}[Theorem~\ref{thm:concave}]
Let $\gamma\in\CC$, and let $\gamma=\sum_jq'_j\gamma^{y_j}$ be any decomposition with $q'_j\ge0$ and arbitrary root sequences $y_j$. Then the distribution of centers of the chain decomposition, with masses $q_k$, majorizes the distribution of the centers $c(y_j)$, with masses $q'_j$. Equivalently, $\Psi_\gamma$ is the largest slope polygon of any decomposition of $\gamma$; and for every convex function $f$ the map $\gamma\mapsto\sum_kq_kf(c(z^k))$ is concave on $\CC$.
\end{theoremB}

Theorem~B makes precise the sense in which the chain decomposition is distinguished among all decompositions of $\gamma$ into supernatural tables. That it is the natural one was suggested in \cite{BEKS12}, where the termwise order on root sequences is characterized by the existence of nonzero homomorphisms between supernatural sheaves. Equivalently, for $m\in\QQ$, $\Psi_\gamma(m)$ is the largest degree of a summand of $\gamma$ in $\CC$ of rank $m$ (Proposition~\ref{prop:splitting}).

Let $e_{j,t}$ be the table with a single entry $1$ in row $j$ and column $t$. We call the tables $e_{j,t}+e_{j+1,t}$, $0\le j<n$, \emph{cancellation vectors}. If $0\to F'\to F\to F''\to0$ is an exact sequence, the long exact cohomology sequences give
\[
\gamma(F')+\gamma(F'')-\gamma(F)=\sum_{t}\sum_{j=0}^{n-1}\rk\big(\partial^j_t\big)\,\big(e_{j,t}+e_{j+1,t}\big),
\]
where $\partial^j_t\colon H^j(F''(t))\to H^{j+1}(F'(t))$ are the connecting homomorphisms. The value of $\phi_z$ on $e_{r-1,t}+e_{r,t}$ is $(n-1)!/\prod_{j\ne r}|t-z_j|$ if $t=z_r$, and $0$ otherwise; in particular it is nonnegative.

\begin{theoremC}[Theorem~\ref{thm:monotone}, Corollary~\ref{cor:monoext}]
Let $\gamma,\gamma'\in\CC$, and suppose that $\gamma'-\gamma$ is a nonnegative combination of cancellation vectors. Then $\Psi_\gamma\le\Psi_{\gamma'}$. Consequently:
\begin{enumerate}[label=\textup{(\roman*)}]
\item $\Psi_F\le\Psi_{F'\oplus F''}$ for every exact sequence $0\to F'\to F\to F''\to0$ of vector bundles on $\PP^n$; in particular $\Psi_F\le\Psi_{\gr F}$ for every filtration with locally free subquotients;
\item if $(F_s)_{s\in S}$ is a flat family of vector bundles on $\PP^n$ over a variety $S$, then every $s_0\in S$ has a neighborhood $U$ with $\Psi_{F_s}\le\Psi_{F_{s_0}}$ for all $s\in U$; if $S$ is irreducible, $\Psi_{F_s}\le\Psi_{F_{s_0}}$ for general $s$ and every $s_0$.
\end{enumerate}
\end{theoremC}

On $\PP^2$ the hypothesis only says that $\gamma\le\gamma'$ entrywise and that both have the same Hilbert polynomial. For every $n\ge3$ this weaker hypothesis does not suffice; on $\PP^3$ the tables of $\OO(-4)^{6}\oplus\OO(-2)^{6}$ and $\OO(-5)\oplus\OO(-4)^{2}\oplus\OO(-3)^{6}\oplus\OO(-2)^{2}\oplus\OO(-1)$ are entrywise comparable, with the same Hilbert polynomial, but their polygons are not (Example~\ref{ex:weak}). Part (ii) is the analogue for $\Psi$ of Shatz's semicontinuity of the Harder--Narasimhan polygon \cite{Sh77}. Theorem~C gives a structural constraint on two problems of Erman and Sam \cite[Problem~4.3 and \S5]{ErS}, on the cohomology tables of the extensions of two given bundles and on the jumps of cohomology tables in flat families: from an extension to the direct sum, and from a general to a special fiber, the center distribution of the chain only moves by mean-preserving spreads. The same computation bounds how much cancellation can lower $\Lambda$ (Proposition~\ref{prop:defectineq}). With $\Lambda$ and $\Psi$ defined by the infimum formula, Theorem~A and Theorem~C(i) hold for all coherent sheaves on $\PP^n$ (Section~\ref{sec:coherent}), and Theorem~D for all coherent sheaves on $\PP^2$ (Proposition~\ref{prop:cohD}); through a general linear projection they extend to sheaves of lower dimension (Theorem~F below). On $\PP^2$ it follows that a coherent sheaf whose decomposition involves only tables of vector bundles is torsion-free, which answers in this case a question of Eisenbud and Schreyer \cite{ES10} (Corollary~\ref{cor:ESquestion}).

\medskip\noindent\textbf{The plane.}
On $\PP^2$ a root sequence is a pair $z=(z_1>z_2)$, and $\phi_{(a,b)}(\gamma)=\big(\gamma_{0,a}-\gamma_{0,b}+\gamma_{1,b}\big)/(a-b)$ involves three entries of the table.

\begin{theoremD}[Theorem~\ref{thm:three}, Proposition~\ref{prop:fourdefect}]
Let $F$ be a vector bundle on $\PP^2$, $\lambda\in\frac12\ZZ$, and $a>b$ integers with $a+b=-2\lambda-3$. Then
\[
(a-b)\,\Lambda_F(\lambda)\ \le\ h^0(F(a))-h^0(F(b))+h^1(F(b)),
\]
with equality whenever $(a,b)$ is comparable with every element of the chain of $F$; such pairs exist for every $\lambda$. For such a pair and a subsheaf $F'\subset F$ of rank $N'$ and slope $\mu'$ with torsion-free quotient $F''$,
\begin{multline*}
(a-b)\big(\Lambda_F(\lambda)-N'(\mu'-\lambda)\big)\\=\big(h^2_b(F')-h^2_a(F')\big)+\big(h^0_a(F'')-h^0_b(F'')\big)+\big(h^1_a(F')-\rk\partial^0_a\big)+\big(h^1_b(F'')-\rk\partial^1_b\big),
\end{multline*}
where $h^j_t(\cG)=h^j(\cG(t))$, and all four terms are nonnegative.
\end{theoremD}

The nonnegativity of the first two terms is the monotonicity of $h^0$ and $h^2$ under multiplication by a linear form, so Theorem~D gives an elementary proof of Theorem~A on $\PP^2$, together with a description of its equality cases, which we use in Section~\ref{sec:lvl}. On $\PP^2$ the majorization of Theorem~B also becomes quantitative: Corollary~\ref{cor:gain} computes how much $\Lambda$ gains on the chain over an arbitrary decomposition.

On $\PP^2$ the argument uses only the tables of $F'$ and $F''$ and the ranks of the connecting homomorphisms, so it gives necessary conditions on the tables of extensions of two given bundles (Corollary~\ref{cor:conicQ}). On $\PP^3$ no argument of this kind can prove Theorem~A: there are tables $\gamma_1,\gamma_2$ of direct sums of supernatural bundles and cancellations compatible with the source--target bounds for the ranks of the connecting homomorphisms, for which the analogous conclusion fails (Example~\ref{ex:conicP3}). The test sheaves of Theorem~A see the injection $F'\otimes \Tt\to F\otimes \Tt$, which the tables do not record, and they exclude these tables. On $\PP^2$ a stronger comparison holds: the Boij--Söderberg polygons of the steps of the Harder--Narasimhan filtration increase, $\Psi_{F_i}\le\Psi_{F_{i'}}$ on $[0,\rk F_i]$ for $i\le i'$ (Theorem~\ref{thm:HNmono}); and $P_F=\Psi_F$ holds exactly when the Harder--Narasimhan factors are supernatural with pairwise comparable roots (Corollary~\ref{cor:equalityP2}). In that case the chain decomposition of $\gamma(F)$ is realized by the Harder--Narasimhan filtration; for Betti tables, the question when the Boij--Söderberg decomposition comes from a filtration of the minimal free resolution is studied in \cite{EES13}.

On $\PP^2$ the Harder--Narasimhan polygon is the large volume limit of the Harder--Narasimhan polygons of the Bridgeland stability conditions $\sigma_{\alpha,\beta}$ \cite{Br08,ABCH13}. The chain decomposition defines a polygon $\Psi^\sigma_F$ in the same plane: for $\alpha$ sufficiently large the central charges $q_kZ_{\alpha,\beta}(E_{z^k})$ of the pieces lie in the open right half-plane, and the edges of $\Psi^\sigma_F$ are these vectors in order of decreasing argument (Definition~\ref{def:PsiSigma}).

\begin{theoremE}[Theorem~\ref{thm:lvl}]
Let $F$ be a vector bundle on $\PP^2$ over $\mathbb C$ and $\beta<\mu_{\min}(F)$. For $\alpha\gg0$ the region bounded by the $\sigma_{\alpha,\beta}$-Harder--Narasimhan polygon of $F$ is contained in the region bounded by $\Psi^\sigma_F$.
\end{theoremE}

The leading order is Theorem~A. The next order is controlled by the equality cases of Theorem~D: at a tight vertex, the pieces above the vertex have the Chern character of the Harder--Narasimhan step (Lemma~\ref{lem:Ktight}); and on a contact edge the Harder--Narasimhan factor is a supernatural bundle (Proposition~\ref{prop:contactedge}), hence Gieseker semistable (Proposition~\ref{prop:sameslope}), so only the vertices of the Gieseker polygon need to be compared. For small $\alpha$ the supernatural pieces may cross walls of their own, and the containment has to be formulated with the Harder--Narasimhan polygon of a supernatural model $\bigoplus_kE_{z^k}^{\oplus mq_k}$; whether it then holds for all $\sigma_{\alpha,\beta}$ with $\beta$ below the Boij--Söderberg slopes of $F$ is open. It does for $F=\OO\oplus E$ with $E$ a stable bundle with $c_1=0$, $c_2=2$ (Remark~\ref{rem:lvl}).

\medskip\noindent\textbf{Projections to \texorpdfstring{$\PP^n$}{Pn}.}
For a coherent sheaf $\cF$ of dimension $n$ on a projective space $\PP^N$, $N\ge n$, write $\chi(\cF(t))=\sum_i(-1)^ie_i(\cF)\binom{t+n-i}{n-i}$. The multiplicity $e_0$ plays the role of the rank and $-e_1$ that of the degree; for pure sheaves we use the normalized slope $\mu_n=-e_1/e_0=\hat\mu-\frac{n+1}2$, where $\hat\mu$ is the slope of \cite[\S1.6]{HL10}. The coefficients of every subsheaf are computed in the same binomial basis of dimension $n$. A general linear projection to $\PP^n$ is finite on the support of $\cF$, and it preserves the cohomology table, the numbers $e_0$ and $e_1$, and inclusions of subsheaves (Lemma~\ref{lem:projection}).

\begin{theoremF}[Theorem~\ref{thm:lowerdim}]
Let $\cF$ be a coherent sheaf of dimension $n\ge1$ on $\PP^N$, and let $\Psi_\cF$ be defined from its table by the infimum formula, with root sequences of length $n$. Then $-e_1(\cF')\le\Psi_\cF(e_0(\cF'))$ for every subsheaf $\cF'\subset\cF$; for $\cF$ pure, the Harder--Narasimhan polygon of $\cF$ with respect to $\mu_n=-e_1/e_0$ lies below $\Psi_\cF$. If $\cF$ is pure of dimension $n$ and locally Cohen--Macaulay, for instance a vector bundle on a Cohen--Macaulay subvariety, then $\Psi_\cF$ is the Boij--Söderberg polygon of a vector bundle on $\PP^n$.
\end{theoremF}

For a graded module $M$ of dimension $n+1$ over $K[x_0,\dots,x_N]$, Theorem~F bounds $e_0$ and $e_1$ of every submodule of the same dimension in terms of the Hilbert functions of $M$ and of its local cohomology modules (Remark~\ref{rem:lowerdim}(2)). The multiplicity bounds of Herzog, Huneke and Srinivasan \cite{HS98}, proved in \cite{ES09,BS12}, concern $e_0(M)$ alone; their naive analogue for submodules, on the side of Betti tables, fails (Remark~\ref{rem:lowerdim}(3)).

A polarized variety is treated in the same way, with a finite projection in place of the linear one. Let $X$ be a smooth projective variety of dimension $n\ge2$, $H$ an ample and base-point-free divisor, $d=H^n$, and $\tau=(K_X^2-2c_2(X))\cdot H^{n-2}/(n+1)$ (denoted $\sigma$ in \cite{An1,An2}); for a surface, $\tau=K_X^2-8\chi(\OO_X)$ is the topological signature over $\mathbb C$. A general projection $\pi\colon X\to\PP^n$ is finite and flat, and $H^j(X,\cE(tH))=H^j(\PP^n,\pi_*\cE(t))$. In \cite{An2} it is shown that if $d<\tau$ then $X$ carries no arithmetically Cohen--Macaulay bundle, i.e.\ no $\cE$ with $\pi_*\cE$ split. The gap moment $G(\cE)=G(\pi_*\cE)$ satisfies $G\ge1$, with equality iff $\cE$ is aCM. Like the variance $V_c(\cE)$ of the centers, it depends only on the table $(h^j(\cE(tH)))_{j,t}=\gamma(\pi_*\cE)$, not on $\pi$. For a torsion-free sheaf $\cQ$ of rank $m$ on $X$ put $v=c_1(\cQ)/m$, $\eta(\cQ)=\big(v-\frac12K_X\big)\cdot H^{n-1}/d$, $\xi(\cQ)=v-\frac12K_X-\eta(\cQ)H$ and $\delta(\cQ)=\big(v^2-2\operatorname{ch}_2(\cQ)/m\big)\cdot H^{n-2}$, and write $\xi(\cQ)^2$ for the number $\xi(\cQ)^2\cdot H^{n-2}$. For the $\mu_H$-Harder--Narasimhan factors $\cQ_i$ of a vector bundle $\cE$ of rank $r$ let $w_i=\rk\cQ_i/r$, and let $V_\eta(\cE)$ be the variance of the numbers $\eta(\cQ_i)$ with weights $w_i$.

\begin{theoremG}[Proposition~\ref{prop:defect}, Corollary~\ref{cor:cases}]
For every vector bundle $\cE$ of rank $r$ on $X$, with $\mu_H$-Harder--Narasimhan factors $\cQ_i$,
\[
\frac{n+1}{12}\big(d\,G(\cE)-\tau\big)=d\big(V_c(\cE)-V_\eta(\cE)\big)+\sum_iw_i\big(\delta(\cQ_i)-\xi(\cQ_i)^2\big),
\]
where $V_c$ is the variance of the centers of the chain; here $\delta(\cQ_i)\ge0$ by Bogomolov's inequality and $\xi(\cQ_i)^2\le0$ by the Hodge index theorem. The terms $d(V_c-V_\eta)$, $\sum_iw_i\delta(\cQ_i)$ and $-\sum_iw_i\xi(\cQ_i)^2$ are nonnegative, the first by Theorem~A. Consequently
\[
G(\cE)\ \ge\ \frac\tau d\qquad\text{and}\qquad\sum_{0<j<n}\sum_th^j(\cE(tH))\ \ge\ C_n\,r\,(\tau-d)_+ .
\]
\end{theoremG}

So the absence of aCM bundles in \cite{An2}, which is the case $G=1$, comes with a uniform lower bound for the gap moment of every cohomology table of a vector bundle on $(X,H)$. Here $C_n=n\,(\lfloor n/2\rfloor-1)!\,(\lceil n/2\rceil-1)!/(12\,(n-2)!)$, which is optimal in the intermediate-cohomology estimate of Lemma~\ref{lem:gapPn}(iv), and $C_2=\frac16$. For the Hesse surfaces $Y_\nu$, $\nu\ge3$, of \cite{An2}, $\tau/d\ge16/7$. For surfaces the integrality of the gaps gives the stronger bound $\sum_th^1(\cE(tH))\ge rd\,L(\kappa)$, $\kappa=\max(1,\tau/d)$, where $L$ interpolates $g^2\mapsto\frac{g^3-g}{12}$ linearly.

The structure of $\Psi$ also bears on a question of Eisenbud and Schreyer (Corollary~\ref{cor:ES}): if all cohomology of $\cE$ vanishes at $n$ twists $a_1<\dots<a_n$, then $\pi_*\cE$ has a supernatural table with roots $a_i$, $\cE$ is $\mu_H$-semistable, and $\frac{n+1}{12}(d\,\rho-\tau)=\delta(\cE)-\xi(\cE)^2\ge0$, where $\rho\ge1$ is the normalized spread of the $a_i$, equal to $e^2$ for an arithmetic progression of difference $e$. On the Hesse surfaces $Y_\nu$ of \cite{An2} with $\nu\ge4$ no vector bundle $\cE$ has $H^j(\cE(-H))=H^j(\cE(H))=0$ for all $j$.

\medskip\noindent\textbf{Relation to known results.}
The last assertion of Theorem~A, the semistability of bundles with supernatural tables, also follows from known results: in any characteristic from the semistability of Ulrich bundles \cite[Thm.~2.9]{CH12}, applied to an Ulrich bundle on $\PP^n\times(\PP^1)^W$ built from $F$, and in characteristic zero from the work of Erman and Sam \cite{ErSa16} (Remark~\ref{rem:main}(3)). The content of Theorem~A is the comparison for arbitrary chains. For Theorem~G, the underlying Riemann--Roch identity on $X$ is \cite[Prop.~10.6, (10.3)]{An3} (for aCM bundles on surfaces, \cite[Prop.~2.6]{An2}); what is new here is its expression through the table, by $G$ and $V_c$, and the comparison $V_\eta\le V_c$ for every bundle, which is Theorem~A (for split push-forwards it is the majorization of \cite[Lemma~10.8]{An3}). In Corollary~\ref{cor:ES}, if the twists $a_i$ form an arithmetic progression of difference $e$, a twist of $\cE$ satisfies the cohomological condition for Ulrich bundles on $(X,eH)$, so when $eH$ is very ample the semistability and the inequality also follow from \cite{CH12} and from the inequality of Lopez for Ulrich bundles \cite[Lem.~3.1]{Lo26}; what is new is the case of root sequences not in arithmetic progression, which occurs only for $n\ge3$ and does not reduce by a twist to the Ulrich condition for a multiple of $H$, and the interpretation as a necessary numerical condition for the realization of a supernatural ray.

\medskip\noindent\textbf{Organization.}
Section~\ref{sec:tables} collects notation on $\PP^n$. Section~\ref{sec:main} proves Theorem~A and introduces the gap moment. Section~\ref{sec:monotone} proves Theorems~B and~C and the decomposition criterion, and extends Theorems~A and~C(i) to coherent sheaves (Section~\ref{sec:coherent}). Section~\ref{sec:plane} specializes to $\PP^2$: gap moment, extreme slopes, Theorem~D, also for coherent sheaves, the tables of extensions, monotonicity along the Harder--Narasimhan filtration, the example on $\PP^3$, and the characterization of the equality $P_F=\Psi_F$. Section~\ref{sec:lvl} proves Theorem~E, Section~\ref{sec:lowerdim} proves Theorem~F, and Section~\ref{sec:polarized} contains the application to polarized varieties (Theorem~G). Appendix~\ref{app:P3} describes, by Beilinson's theorem, the vector bundles on $\PP^3$ whose table is a multiple of $\gamma^z$ with $z_1-z_3\le4$; by Theorem~A they are all semistable, and for $z_1-z_3=4$ they are, up to a twist, instanton bundles, Steiner bundles or their duals (Corollary~\ref{cor:appP3}).

\medskip\noindent\textbf{Ground field.}
The characteristic enters through Bogomolov's inequality in Corollary~\ref{cor:Q} and in Section~\ref{sec:polarized}, through Theorem~\ref{thm:lvl} and Remark~\ref{rem:lvl}, where we work over $\mathbb C$ as in the references on Bridgeland stability, and through results we quote: the comparison with \cite{ErSa16} in Remark~\ref{rem:main}(3) and the stability statements in Appendix~\ref{app:P3}. Theorem~A, Sections~\ref{sec:tables}--\ref{sec:monotone} (except for the inequality $\delta(\cQ_i)\ge0$ in Corollary~\ref{cor:Q}), Theorems~\ref{thm:three} and~\ref{thm:HNmono}, Proposition~\ref{prop:fourdefect}, the results of Section~\ref{sec:equality} and Theorem~F hold over any algebraically closed field (see the table), since the pure resolutions and the supernatural bundles of \cite[Thms.~5.1, 6.1]{ES09} exist over any field.

\begin{table}[ht]\centering\small
\textsc{Ground fields of the main results}\par\smallskip
\begin{tabular}{@{}lll@{}}
\toprule
Result & Section & Ground field\\
\midrule
Theorem~A, Propositions~\ref{prop:sameslope}, \ref{prop:marginPn}, Lemma~\ref{lem:gapPn} & \ref{sec:main} & any\\
Corollary~\ref{cor:Q} (Bogomolov) & \ref{sec:main} & characteristic $0$\\
Theorems~B, C & \ref{sec:monotone} & any\\
Propositions~\ref{prop:cohA}, \ref{prop:cohseries}, \ref{prop:cohC}, Corollary~\ref{cor:ESquestion} & \ref{sec:coherent} & any\\
Theorem~D, Proposition~\ref{prop:cohD}, Theorem~\ref{thm:HNmono} & \ref{sec:plane} & any\\
Lemmas~\ref{lem:contact}--\ref{lem:gapdisc}, Proposition~\ref{prop:contactedge}, Corollary~\ref{cor:equalityP2} & \ref{sec:equality} & any\\
Theorem~E, Lemmas~\ref{lem:LV}, \ref{lem:support}, Remark~\ref{rem:lvl} & \ref{sec:lvl} & $\mathbb C$\\
Theorem~F & \ref{sec:lowerdim} & any\\
Lemma~\ref{lem:push}, Proposition~\ref{prop:defect} & \ref{sec:polarized} & any\\
Theorem~G, Corollary~\ref{cor:ES} (Bogomolov) & \ref{sec:polarized} & characteristic $0$\\
Stability of general bundles in Appendix~\ref{app:P3} & \ref{app:P3} & characteristic $0$\\
\bottomrule
\end{tabular}
\end{table}

\section{Cohomology tables on \texorpdfstring{$\PP^n$}{Pn}}\label{sec:tables}

\subsection{Supernatural tables}
A \emph{root sequence} is a strictly decreasing sequence of integers $z=(z_1>\dots>z_n)$, with center $c(z)=\frac1n\sum_jz_j\in\frac1n\ZZ$ and slope $\mu(z)=-c(z)-\frac{n+1}2\in\frac1n\ZZ$; note that $\frac{n+1}2=\frac1n\cdot\frac{n(n+1)}2\in\frac1n\ZZ$. Root sequences are ordered termwise, and a \emph{chain} is a totally ordered set. Two comparable root sequences with the same center are equal. The normalized supernatural table $\gamma^z$ is the cohomology table of a supernatural bundle with root sequence $z$ \cite[Thm.~0.4, \S6]{ES09}, scaled so that its Hilbert polynomial is $\prod_j(t-z_j)$:
\[
\gamma^z_{j,t}=\begin{cases}\big|\prod_{i=1}^n(t-z_i)\big| & \text{if } j=\#\{i:z_i>t\},\\ 0 & \text{otherwise.}\end{cases}
\]
So each column has at most one nonzero entry, the column $t$ vanishes iff $t$ is a root, and $\sum_j(-1)^j\gamma^z_{j,t}=\prod_i(t-z_i)$.

\begin{theorem}[Eisenbud--Schreyer {\cite[Thm.~0.5]{ES09}}]\label{thm:ES}
For every vector bundle $F$ on $\PP^n$ there are a unique finite chain $z^1>z^2>\dots$ and unique $q_k\in\QQ_{>0}$ with $\gamma(F)=\sum_kq_k\gamma^{z^k}$.
\end{theorem}

Since supernatural bundles exist for every root sequence \cite[Thm.~0.4]{ES09}, a positive integer multiple of every nonnegative rational combination of supernatural tables is the table of a direct sum of supernatural bundles. So Theorem~\ref{thm:ES} applies to all such combinations. We write $\CC$ for the set of nonnegative rational combinations of the $\gamma^z$; every $\gamma\in\CC$ has a unique chain decomposition with positive rational coefficients. We compare tables entrywise. For $\gamma=0$ the chain is empty, $\Psi_0$ is the point $(0,0)$ and $\Lambda_0=0$; statements about the first and last elements of a chain, which we denote $z^{\max}$ and $z^{\min}$, assume $\gamma\ne0$. Throughout, (semi)stable means slope (Mumford--Takemoto) (semi)stable with respect to $\OO(1)$; $P_F$, $\mu_{\max}(F)$ and $\mu_{\min}(F)$ refer to the corresponding Harder--Narasimhan filtration.

\subsection{Rank, degree and the polygon}
For $\gamma=\sum_kq_k\gamma^{z^k}\in\CC$ put
\[
\rk\gamma=n!\sum_kq_k,\qquad \deg\gamma=\sum_kn!\,q_k\,\mu(z^k).
\]
Both are linear in $\gamma$, because they are read off from the two leading coefficients of the Hilbert polynomial $\sum_kq_k\prod_j(t-z^k_j)$:
\[
\sum_kq_k\prod_j(t-z^k_j)=\frac{\rk\gamma}{n!}\,t^n+\Big(\deg\gamma+\frac{n+1}2\rk\gamma\Big)\frac{t^{n-1}}{(n-1)!}+\dots
\]
By Riemann--Roch, $\chi(\cF(t))=\frac{\rk\cF}{n!}t^n+\big(\deg\cF+\frac{n+1}2\rk\cF\big)\frac{t^{n-1}}{(n-1)!}+\dots$ for every coherent sheaf $\cF$, where $\deg\cF=c_1(\cF)\cdot h^{n-1}$ and $\rk\cF$ may be zero, so $\rk\gamma(F)=\rk F$ and $\deg\gamma(F)=\deg F$ for a vector bundle $F$. The \emph{Boij--Söderberg slopes} of $\gamma$ are the $\mu(z^k)$, with ranks $n!\,q_k$, and
\[
\Psi_\gamma(m)=\max\Big\{\sum_kx_k\,\mu(z^k)\ :\ 0\le x_k\le n!\,q_k,\ \sum_kx_k=m\Big\},\qquad 0\le m\le\rk\gamma .
\]
We also use the concave conjugate in the form
\[
\Lambda_\gamma(\lambda)=\sum_kn!\,q_k\,\big(\mu(z^k)-\lambda\big)_+=\sum_kn!\,q_k\,\big(c_\lambda-c(z^k)\big)_+,\qquad c_\lambda=-\lambda-\tfrac{n+1}2 ,
\]
so that $c(z)=c_\lambda$ iff $\mu(z)=\lambda$. It is convex and piecewise linear in $\lambda$, with breakpoints at the Boij--Söderberg slopes, which lie in $\frac1n\ZZ$. For a polygon $P$ with slopes $\mu_i$ and lengths $n_i$ put $\Lambda_P(\lambda)=\sum_in_i(\mu_i-\lambda)_+$.

\begin{lemma}\label{lem:legendre}
\begin{enumerate}[label=\textup{(\roman*)}]
\item Let $\gamma\ne0$. Then $\Psi_\gamma(m)=\min_\lambda\big(m\lambda+\Lambda_\gamma(\lambda)\big)$ for $0\le m\le\rk\gamma$, and the minimum is attained at a Boij--Söderberg slope of $\gamma$; in particular at a point of $\frac1n\ZZ$.
\item For two concave polygons $P,\Psi$ on $[0,N]$ with the same endpoints, $P\le\Psi$ iff $\Lambda_P\le\Lambda_\Psi$ on $\RR$, iff the distribution of slopes of $P$ is majorized by that of $\Psi$. Since $\Psi$ is concave and $P$ is linear between its vertices, $P\le\Psi$ iff $P(m)\le\Psi(m)$ at the vertices $m$ of $P$.
\item Let $\gamma=\gamma(F)$ and let $F_1\subset F$ have rank $n_1>0$ and slope $\mu_1$. Then $\Psi_F(n_1)\ge n_1\mu_1$ iff $\Lambda_F(\lambda)\ge n_1(\mu_1-\lambda)$ for all $\lambda\in\frac1n\ZZ$. If $F/F_1$ has rank $n_2>0$ and slope $\mu_2<\mu_1$, it suffices to check $\lambda\in(\mu_2,\mu_1)\cap\frac1n\ZZ$.
\end{enumerate}
\end{lemma}

\begin{proof}
(i) is the duality of the fractional knapsack problem: $m\lambda+\Lambda_\gamma(\lambda)$ is convex and piecewise linear in $\lambda$, with slope $m-\sum_{\mu(z^k)>\lambda}n!\,q_k$, which is negative for $\lambda$ below the smallest and positive above the largest Boij--Söderberg slope when $0<m<\rk\gamma$; so its minimum is at a breakpoint, and the endpoints $m=0,\rk\gamma$ are checked directly. Since $\Lambda_P(\lambda)=\max_m(P(m)-\lambda m)$ for a concave polygon with $P(0)=0$, the minimum equals $\Psi_\gamma(m)$. (ii) The first equivalence is the comparison of concave conjugates; the second is the theorem of Hardy, Littlewood and P\'olya \cite[Thm.~108]{HLP}, since both distributions have the same total mass and mean (for rational masses, repeat the points; the general case follows by continuity). The last assertion holds because $\Psi$ is concave. (iii) The first assertion follows from (i). For $\lambda\ge\mu_1$ the right-hand side is $\le0$. For $\lambda\le\mu_2$, the inequality $(x-\lambda)_+\ge x-\lambda$ gives $\Lambda_F(\lambda)\ge\deg F-\rk F\cdot\lambda=n_1(\mu_1-\lambda)+n_2(\mu_2-\lambda)\ge n_1(\mu_1-\lambda)$.
\end{proof}

\subsection{Duality}
By Serre duality $h^j(F^\vee(t))=h^{n-j}(F(-t-n-1))$. On tables this is the involution $\gamma\mapsto\gamma^\vee$, $\gamma^\vee_{j,t}=\gamma_{n-j,-t-n-1}$, and $(\gamma^z)^\vee=\gamma^{z^\vee}$ with $z^\vee=(-n-1-z_n,\dots,-n-1-z_1)$. The map $z\mapsto z^\vee$ reverses the order and $\mu(z^\vee)=-\mu(z)$. So the chain of $F^\vee$ is $\{(z^k)^\vee\}$ with the same coefficients, its Boij--Söderberg slopes are the negatives of those of $F$, and
\begin{equation}\label{eq:dualpsi}
\Psi_{F^\vee}(N-m)=\Psi_F(m)-\deg F,\qquad N=\rk F .
\end{equation}
The Harder--Narasimhan polygon satisfies the same relation.

\subsection{The ends of the chain}
The first and last elements of the chain are read off from the supports of the cohomology rows.

\begin{lemma}[ends of the chain]\label{lem:ends}
Let $\gamma=\sum_lq'_l\gamma^{y_l}$ with all $q'_l>0$, for arbitrary root sequences $y_l$. Then the termwise maximum $y^{\max}$ and the termwise minimum $y^{\min}$ of the $y_l$ are given by
\[
y^{\max}_i=1-i+\max\{t+j:\ j\ge i,\ \gamma_{j,t}\ne0\},\qquad y^{\min}_i=\min\{t+j:\ j<i,\ \gamma_{j,t}\ne0\}-i .
\]
In particular they depend only on $\gamma$, and they are the first and the last element, $z^{\max}$ and $z^{\min}$, of the chain of $\gamma$.
\end{lemma}

\begin{proof}
All entries are nonnegative, so $\gamma_{j,t}\ne0$ iff $\gamma^{y_l}_{j,t}\ne0$ for some $l$; it suffices to prove the formulas for a single root sequence $y$. If $\gamma^y_{j,t}\ne0$ with $j\ge i$, then exactly $j$ roots exceed $t$, so $y_j>t$ and $y_i\ge y_j+(j-i)\ge t+1+j-i$. Conversely, let $j\ge i$ be maximal with $y_j=y_i-(j-i)$, and $t=y_j-1$. Then $t$ is not a root, because $y_{j+1}<y_j-1$ by maximality (or $j=n$), and exactly $y_1,\dots,y_j$ exceed $t$; so $\gamma^y_{j,t}\ne0$, and $1-i+t+j=y_i$. Symmetrically, if $\gamma^y_{j,t}\ne0$ with $j<i$, then $y_{j+1}<t$ and $y_i\le y_{j+1}-(i-j-1)\le t+j-i$; and if $j+1\le i$ is minimal with $y_{j+1}=y_i+(i-j-1)$, then $t=y_{j+1}+1$ is not a root, exactly $j$ roots exceed it, and $t+j-i=y_i$. Applied to the chain, whose elements are totally ordered, this identifies $y^{\max}$ and $y^{\min}$ with its first and last elements.
\end{proof}

\section{The comparison theorem}\label{sec:main}

\subsection{Divided differences}
For a root sequence $z$ put $D_r(z)=\prod_{i\ne r}(z_r-z_i)$ and $P_r(z)=|D_r(z)|$. Exactly $r-1$ of the $z_i$ exceed $z_r$, so $D_r(z)=(-1)^{r-1}P_r(z)$. For a table $\gamma$ put
\begin{equation}\label{eq:phi}
\phi_z(\gamma)=(n-1)!\sum_{r=1}^n\frac{1}{D_r(z)}\sum_{j<r}(-1)^j\gamma_{j,z_r} .
\end{equation}
This is a linear functional, which involves the entries of the rows $0,\dots,r-1$ in the column $z_r$. For a polynomial $f$, $\sum_rf(z_r)/D_r(z)$ is the divided difference $f[z_1,\dots,z_n]$; if $f$ has degree at most $n$, $f(t)=a_nt^n+a_{n-1}t^{n-1}+\dots$, then
\begin{equation}\label{eq:divdiff}
\sum_{r=1}^n\frac{f(z_r)}{D_r(z)}=a_n\sum_iz_i+a_{n-1}.
\end{equation}
For $n=2$ and $z=(a,b)$, $\phi_z(\gamma)=\big(\gamma_{0,a}-\gamma_{0,b}+\gamma_{1,b}\big)/(a-b)$. For $n=1$, $\phi_z(\gamma)=\gamma_{0,z_1}$.

\begin{lemma}\label{lem:comparable}
If $y$ and $z$ are comparable root sequences, then $\phi_z(\gamma^y)=n!\,\big(c(z)-c(y)\big)_+$.
\end{lemma}

\begin{proof}
If $y\le z$, then for $i\ge r$ we have $y_i\le z_i\le z_r$, so at most $r-1$ roots of $y$ exceed $z_r$, and the nonzero entry of the column $z_r$ of $\gamma^y$, if any, lies in a row $j<r$. Hence $\sum_{j<r}(-1)^j\gamma^y_{j,z_r}=\prod_i(z_r-y_i)$, also when $z_r$ is a root of $y$. By \eqref{eq:divdiff} applied to the monic polynomial $\prod_i(t-y_i)$, $\phi_z(\gamma^y)=(n-1)!\big(\sum_iz_i-\sum_iy_i\big)=n!\,(c(z)-c(y))\ge0$. If $y\ge z$, then $y_i\ge z_i\ge z_r$ for $i\le r$, so either $z_r$ is a root of $y$ or at least $r$ roots of $y$ exceed it; in both cases $\gamma^y_{j,z_r}=0$ for $j<r$, and $\phi_z(\gamma^y)=0=n!\,(c(z)-c(y))_+$.
\end{proof}

For root sequences $x\le y$ let $[x,y]=\{z:x\le z\le y\}$, and write $\epsilon_1,\dots,\epsilon_n$ for the unit vectors of $\ZZ^n$.

\begin{lemma}\label{lem:chains}
Let $x\le y$ be root sequences and $\ell=\sum_i(y_i-x_i)$. Lowering first the coordinate $n$ from $y_n$ to $x_n$, then the coordinate $n-1$, and so on, one unit at a time, passes only through root sequences. So there is a chain $y=z^0>z^1>\dots>z^\ell=x$ in $[x,y]$ with $z^k-z^{k+1}$ a unit vector and $c(z^k)=c(y)-\frac kn$.
\end{lemma}

\begin{proof}
Suppose the coordinates $n,\dots,l+1$ have been lowered to $x$ and the coordinate $l$ is being lowered from a value $z_l>x_l$. Then $z_l-1\ge x_l>x_{l+1}=z_{l+1}$ and $z_{l-1}=y_{l-1}>y_l\ge z_l$, so $z-\epsilon_l$ is a root sequence.
\end{proof}

\begin{lemma}[adapted root sequences]\label{lem:adapted}
Let $\gamma\in\CC$, $\gamma\ne0$, with chain $z^1>z^2>\dots$, and let $\lambda\in\frac1n\ZZ$. There is a root sequence $z$ with $\mu(z)=\lambda$ that is comparable with every $z^k$; we call it \emph{adapted} to $\gamma$ at $\lambda$. For every such $z$, $\phi_z(\gamma)=\Lambda_\gamma(\lambda)$. In particular $\phi_{z^k}(\gamma)=\Lambda_\gamma(\mu(z^k))$ for every $k$.
\end{lemma}

\begin{proof}
Put $c=c_\lambda\in\frac1n\ZZ$. If $c\ge c(z^{\max})$, take $z=z^{\max}+n(c-c(z^{\max}))\,\epsilon_1$, which lies above $z^{\max}$, hence above every $z^k$. If $c\le c(z^{\min})$, take $z=z^{\min}-n(c(z^{\min})-c)\,\epsilon_n$. Otherwise $c(z^{k+1})\le c\le c(z^k)$ for some $k$, and the chain of Lemma~\ref{lem:chains} for $[z^{k+1},z^k]$ passes through a root sequence $z$ with $c(z)=c$; it satisfies $z^{k+1}\le z\le z^k$, so it is comparable with every element of the chain. In all cases Lemma~\ref{lem:comparable} and linearity give $\phi_z(\gamma)=\sum_kq_k\,n!\,(c-c(z^k))_+=\Lambda_\gamma(\lambda)$.
\end{proof}

\subsection{Test sheaves}
A sheaf set-theoretically supported on a line, with a pure resolution of prescribed type, exists by the theorem on pure resolutions.

\begin{lemma}[test sheaves]\label{lem:testsheaf}
Let $z$ be a root sequence and $L\subset\PP^n$ a line. There are a nonzero coherent sheaf $\Tt$ with support $L$, an integer $B>0$ and an exact sequence
\[
0\to\OO(z_n)^{b_{n-1}}\to\dots\to\OO(z_2)^{b_1}\to\OO(z_1)^{b_0}\to \Tt\to0,\qquad b_i=\frac B{P_{i+1}(z)} .
\]
\end{lemma}

\begin{proof}
For $n=1$ take $L=\PP^1$ and $\Tt=\OO(z_1)$. Let $n\ge2$, $R=K[x_1,\dots,x_{n-1}]$ and $d_i=-z_{i+1}$, so that $d_0<\dots<d_{n-1}$. By the existence of pure resolutions, due to Eisenbud, Fløystad and Weyman in characteristic zero \cite{EFW11} and to Eisenbud and Schreyer in general \cite[Thm.~5.1]{ES09}, there is a graded $R$-module $M\ne0$ of finite length with a pure resolution $0\to R(-d_{n-1})^{b_{n-1}}\to\dots\to R(-d_0)^{b_0}\to M\to0$. By Herzog and Kühl \cite{HK84}, the Betti numbers of a pure resolution are proportional to $1/\prod_{j\ne i}|d_j-d_i|=1/P_{i+1}(z)$. Tensoring with the flat extension $R\to S=R[x_0,x_n]$ and passing to sheaves on $\PP^n=\operatorname{Proj}S$ gives a resolution of the sheaf associated with $M\otimes_RS$, which is supported on the line $x_1=\dots=x_{n-1}=0$. A linear change of coordinates moves this line to $L$.
\end{proof}

\begin{lemma}[truncation]\label{lem:truncation}
Let $F$ be a vector bundle on $\PP^n$, $z$ a root sequence and $\Tt$, $B$ as in Lemma~\textup{\ref{lem:testsheaf}}. Then
\[
h^0(F\otimes \Tt)\ \le\ \frac B{(n-1)!}\,\phi_z\big(\gamma(F)\big).
\]
\end{lemma}

\begin{proof}
Since $F$ is locally free, $F$ tensored with the resolution of $\Tt$ is a resolution $K^\bullet$ of $F\otimes \Tt$, with $K^{-i}=F(z_{i+1})^{b_i}$, $0\le i\le n-1$. The hypercohomology spectral sequence of $K^\bullet$ has $E_1^{-i,j}=H^j(F(z_{i+1}))^{b_i}$ and converges to $H^{j-i}(F\otimes \Tt)$. On each page put $S_r=\sum_{k\le0}(-1)^k\sum_{j-i=k}\dim E_r^{-i,j}$. The differentials raise the total degree $k$ by one; those between total degrees $k-1$ and $k$ with $k\le0$ do not change $S_r$, and those from total degree $0$ to total degree $1$ lower it by their rank. So $S_{r+1}\le S_r$, and at the limit $S_\infty=\sum_{k\le0}(-1)^kh^k(F\otimes \Tt)=h^0(F\otimes \Tt)$, because the hypercohomology of a sheaf vanishes in negative degrees. Hence
\[
h^0(F\otimes \Tt)\le S_1=\sum_{i=0}^{n-1}b_i\sum_{j\le i}(-1)^{j-i}h^j\big(F(z_{i+1})\big)=B\sum_{r=1}^n\frac{(-1)^{r-1}}{P_r(z)}\sum_{j<r}(-1)^jh^j\big(F(z_r)\big),
\]
with $r=i+1$, and $(-1)^{r-1}/P_r(z)=1/D_r(z)$.
\end{proof}

In the notation of \cite{ES09}, where $\langle\beta,\gamma\rangle=\sum_{j\le i}\sum_k(-1)^{i-j}\beta_{i,k}\gamma_{j,-k}$ for a Betti table $\beta$ and a cohomology table $\gamma$, the sum $S_1$ is $\langle\beta,\gamma(F)\rangle$ for the Betti table $\beta$ of the resolution of Lemma~\ref{lem:testsheaf}, with entries $\beta_{i,-z_{i+1}}=b_i$. So $\phi_z=\frac{(n-1)!}B\,\langle\beta,-\rangle$ is, up to a positive factor, the Eisenbud--Schreyer pairing with a pure Betti table of codimension $n-1$, and Lemma~\ref{lem:truncation} says that $\langle\beta,\gamma(F)\rangle\ge h^0(F\otimes\Tt)$. This is the sheaf form of an inequality contained in the proof of \cite[Thm.~3.1]{ES09}, where the pairing of a free resolution of a module $M$ with a free complex $E$ is written as the dimension of the degree zero part of $H^0_{\mathrm{tot}}$, a submodule of $M\otimes E^0$, plus a nonnegative term. A pairing between derived categories from which all the functionals of \cite{ES09} are obtained is constructed in \cite{EE17}.

\subsection{The theorem}
\begin{theorem}\label{thm:main}
Let $F$ be a vector bundle on $\PP^n$ and $F'\subset F$ a subsheaf of rank $N'$.
\begin{enumerate}[label=\textup{(\alph*)}]
\item For every root sequence $z$, $\deg F'-\mu(z)\,N'\le\phi_z(\gamma(F))$.
\item $\deg F'\le\Psi_F(N')$. Equivalently, $P_F\le\Psi_F$.
\end{enumerate}
\end{theorem}

\begin{proof}
(a) If $F'=0$, the assertion is $\phi_z(\gamma(F))\ge0$, which holds by Lemma~\ref{lem:truncation}. Otherwise, the saturation of $F'$ has the same rank and at least the same degree, so we may assume that $F/F'$ is torsion-free. Then $F'$ and $F/F'$ are locally free outside a closed subset of codimension at least two, and a general line $L$ is disjoint from it. Having chosen $L$, let $\Tt$ be the sheaf of Lemma~\ref{lem:testsheaf} supported on $L$. Near $L$ the sequence $0\to F'\to F\to F/F'\to0$ is a locally split exact sequence of vector bundles, and away from $L$ the resolution of $\Tt$ is a locally split exact complex; so $F'\otimes \Tt\to F\otimes \Tt$ is injective, and $F'$ tensored with the resolution of $\Tt$ is a resolution of $F'\otimes \Tt$. Since $F'\otimes \Tt$ has one-dimensional support,
\[
\chi(F'\otimes \Tt)\ \le\ h^0(F'\otimes \Tt)\ \le\ h^0(F\otimes \Tt)\ \le\ \frac B{(n-1)!}\,\phi_z(\gamma(F))
\]
by Lemma~\ref{lem:truncation}. On the other hand $\chi(F'\otimes \Tt)=\sum_i(-1)^ib_i\chi(F'(z_{i+1}))=B\sum_r\chi(F'(z_r))/D_r(z)$, and by Riemann--Roch and \eqref{eq:divdiff},
\[
(n-1)!\sum_r\frac{\chi(F'(z_r))}{D_r(z)}=(n-1)!\Big(\frac{N'}{n!}\sum_iz_i+\frac{\deg F'+\frac{n+1}2N'}{(n-1)!}\Big)=\deg F'-\mu(z)\,N'.
\]
(b) If $F'=0$ there is nothing to prove. Otherwise apply (a) to the elements $z^k$ of the chain of $F$: by Lemma~\ref{lem:adapted}, $\deg F'\le N'\mu(z^k)+\Lambda_F(\mu(z^k))$ for every $k$, and by Lemma~\ref{lem:legendre}(i) the minimum over $k$ of the right-hand side is $\Psi_F(N')$. Since $P_F$ is the upper concave envelope of the points $(\rk F',\deg F')$ for subsheaves $F'\subset F$ (its vertices are the steps of the Harder--Narasimhan filtration, and every subsheaf lies below it, see \cite[\S1.3]{HL10} and \cite{Sh77}), this is $P_F\le\Psi_F$.
\end{proof}

\begin{corollary}\label{cor:supernatural}
Let $F$ be a vector bundle of rank $N$ on $\PP^n$.
\begin{enumerate}[label=\textup{(\alph*)}]
\item If $\gamma(F)$ is a multiple of a single supernatural table, then $F$ is semistable.
\item $\mu_{\max}(F)\le\mu^\BS_{\max}(F)=\mu(z^{\min})$ and $\mu_{\min}(F)\ge\mu^\BS_{\min}(F)=\mu(z^{\max})$, where $z^{\max}$ and $z^{\min}$ are the first and the last element of the chain, given by Lemma~\textup{\ref{lem:ends}} in terms of the supports of the rows of $\gamma(F)$.
\item $V_{HN}(F)\le V_c(F)$, where $V_{HN}(F)=\sum_iw_i(\mu_i-\mu(F))^2$ is the variance of the slopes $\mu_i$ of the Harder--Narasimhan factors $\cQ_i$ of $F$, with weights $w_i=\rk\cQ_i/N$, and $V_c(F)$ is the variance of the Boij--Söderberg slopes, with weights $p_k=n!\,q_k/N$. Equality holds if and only if $P_F=\Psi_F$.
\end{enumerate}
\end{corollary}

\begin{proof}
(a) Then $\Psi_F$ is linear of slope $\mu(F)$, and Theorem~\ref{thm:main} says $\mu(F')\le\mu(F)$ for every subsheaf. (b) For small $x>0$ we have $P_F(x)=x\,\mu_{\max}(F)$ and $\Psi_F(x)=x\,\mu^\BS_{\max}(F)$, so $P_F\le\Psi_F$ gives the first inequality; the second follows at the other end, or by \eqref{eq:dualpsi}. (c) Both polygons are concave with the same endpoints, so $P_F\le\Psi_F$ says that the distribution of the Harder--Narasimhan slopes is majorized by that of the Boij--Söderberg slopes, both with mean $\mu(F)$; apply Karamata's inequality \cite{Ka32} to $x\mapsto x^2$. After clearing denominators both distributions have the same finite total mass; for the strictly convex function $x^2$, equality in the inequality of Hardy, Littlewood and P\'olya \cite{HLP} forces the two distributions to coincide, i.e.\ $P_F=\Psi_F$. The converse is clear.
\end{proof}

\begin{proposition}[subsheaves of the same slope]\label{prop:sameslope}
Let $F$ be a vector bundle on $\PP^n$ whose table is a multiple of $\gamma^z$, and let $\cG\subset F$ be a saturated subsheaf with $0\ne \cG\ne F$ and $\mu(\cG)=\mu(F)$. Then $\cG$ and $F/\cG$ are vector bundles whose tables are multiples of $\gamma^z$. Consequently $F$ is Gieseker semistable, and $F$ is Gieseker stable if and only if it is $\mu$-stable.
\end{proposition}

\begin{proof}
Let $D=z_1-z_n+1$, let $W$ be the set of integers in $[z_n,z_1]$ that are not roots, so $|W|=D-n$, and let $X=\PP^n\times(\PP^1)^W$ with the very ample $A=\OO_X(1,\dots,1)$; $\dim X=D$. For a coherent sheaf $\cT$ on $\PP^n$ put
\[
\cU_\cT=\Big(\cT\boxtimes\mathop{\boxtimes}_{w\in W}\OO_{\PP^1}(-w-1)\Big)\big((z_1+1)A\big).
\]
For $1\le j\le D$ put $t=z_1+1-j\in[z_n,z_1]$; by the Künneth formula $H^*(\cU_\cT(-jA))=H^*(\cT(t))\otimes\bigotimes_{w}H^*(\OO_{\PP^1}(t-w-1))$. If $t$ is a root, $H^*(F(t))=0$; if $t\in W$, the factor $\OO_{\PP^1}(-1)$ is acyclic. So $\cU_F$ has no cohomology at the twists $-A,\dots,-DA$, i.e.\ it is an Ulrich bundle on $(X,A)$ \cite[Prop.~2.1]{ESW03}. (This goes in the opposite direction to the construction of supernatural bundles in \cite[\S6]{ES09}, as push-forwards of line bundles from products of projective spaces under finite linear projections.) The functor $\cT\mapsto\cU_\cT$ is exact and preserves torsion-freeness, so $0\to\cU_\cG\to\cU_F\to\cU_{F/\cG}\to0$ is exact with torsion-free quotient. For $\cT$ of rank $s>0$, $c_1(\cU_\cT)\cdot A^{D-1}=\frac{(D-1)!}{(n-1)!}\deg\cT+s\kappa$ with $\kappa$ independent of $\cT$, so $\mu_A(\cU_\cG)=\mu_A(\cU_F)$. By \cite[Thm.~2.9(b)]{CH12} (Casanellas and Hartshorne work over an algebraically closed field of arbitrary characteristic, and the proof uses Bertini's theorem for general hyperplane sections, Serre vanishing and Serre duality), $\cU_\cG$ and $\cU_{F/\cG}$ are Ulrich bundles; restricting to $\PP^n\times\{p\}$ shows that $\cG$ and $F/\cG$ are locally free. For $t=z_i$, every factor $\OO_{\PP^1}(z_i-w-1)$, $w\in W$, has a nonzero cohomology group, while $\cU_\cG(-(z_1+1-z_i)A)$ is acyclic; so $H^*(\cG(z_i))=0$, and likewise $H^*((F/\cG)(z_i))=0$. Hence the columns $z_1,\dots,z_n$ of the tables of $\cG$ and $F/\cG$ vanish, and since these tables are positive combinations of supernatural tables, every element of their chains is $z$.

The reduced Hilbert polynomial of a bundle with table $q\gamma^z$ is $\prod_i(t-z_i)/n!$. If $\cG'\subset F$ is a subsheaf with $\mu(\cG')<\mu(F)$, its reduced Hilbert polynomial is smaller. If $\mu(\cG')=\mu(F)$, its saturation $\cG''$ has $\mu(\cG'')\ge\mu(\cG')=\mu(F)$, hence $\mu(\cG'')=\mu(F)$ by Corollary~\ref{cor:supernatural}(a); if $\cG''\ne F$ it has the reduced Hilbert polynomial of $F$ by the first part, and if $\cG''=F$ this is trivial. The reduced Hilbert polynomial of $\cG'$ is at most that of $\cG''$, since $\cG''/\cG'$ is a torsion sheaf. So $F$ is Gieseker semistable; and if it is not $\mu$-stable, a saturated $\cG\subset F$ with $\mu(\cG)=\mu(F)$ has the same reduced Hilbert polynomial as $F$, so $F$ is not Gieseker stable. A $\mu$-stable sheaf is Gieseker stable.
\end{proof}

\begin{remark}\label{rem:main}
\begin{enumerate}[label=(\arabic*)]
\item Equality $P_F=\Psi_F$ holds whenever the Harder--Narasimhan factors $\cQ_i=F_i/F_{i-1}$, $1\le i\le s$, of $F$ are vector bundles with tables $q_i\gamma^{z(i)}$ and pairwise comparable $z(i)$, for instance for a direct sum of supernatural bundles with comparable root sequences. Indeed the slopes decrease, so $z(1)<\dots<z(s)$; for $z'<z$ the row of $\gamma^{z'}$ in each column is at most the row of $\gamma^z$, so, by induction on $i$, at each twist the cohomology of $F_{i-1}$ lies in rows at most the row of $\cQ_i$, every connecting homomorphism $H^j(\cQ_i(t))\to H^{j+1}(F_{i-1}(t))$ vanishes, and $\gamma(F_i)=\gamma(F_{i-1})+\gamma(\cQ_i)$. Then $\gamma(F)=\sum_i\gamma(\cQ_i)$ is a positive combination along a chain, hence the chain decomposition, and its slopes and lengths are those of $P_F$. Consequently Theorem~\ref{thm:main} is optimal among bounds that depend only on the ray of the table: by \cite[Thm.~6.1]{ES09} every ray $\RR_{>0}\gamma^z$ contains the table of a supernatural bundle, so every $\gamma\in\CC$ has a positive integer multiple $M\gamma$ that is the table of a direct sum $F$ of supernatural bundles along the chain of $\gamma$, and then $P_F(Mm)=M\Psi_\gamma(m)$. So $\Psi_\gamma(m)$ is the maximum of $P_F(Mm)/M$ over all $M\ge1$ and all vector bundles $F$ with $\gamma(F)=M\gamma$. For a fixed table integrality gives more: the vertices of $P_F$ are lattice points, so $P_F\le\Psi^{\ZZ}_F$, the upper concave hull of the points $(m,\lfloor\Psi_F(m)\rfloor)$, $m=0,\dots,\rk F$, which also depends only on $\gamma(F)$. In particular $F$ is semistable if $\Psi^{\ZZ}_F$ is the segment from $(0,0)$ to $(\rk F,\deg F)$. This happens for the table of Example~\ref{ex:sum} below: there $\Psi_F$ takes the values $0,-\frac52,-\frac{16}3,-\frac{17}2,-12$ at $m=0,\dots,4$, so every bundle with this table is semistable, although $\Psi_F$ is not linear.
\item The root sequences of the chain lie within the cohomology range of $F$ (proof of \cite[Cor.~0.6]{ES09}), so $\mu^\BS_{\max}(F)$ and $\mu^\BS_{\min}(F)$ lie between the two bounds that \cite[Cor.~0.6]{ES09} gives for $\mu(F)$ in terms of the cohomology range. By Corollary~\ref{cor:supernatural}(b), the upper bound holds for $\mu_{\max}(F)$, i.e.\ for the slope of every subsheaf of $F$, and the lower bound for $\mu_{\min}(F)$, i.e.\ for the slope of every torsion-free quotient. The reverse bounds fail: on $\PP^2$, $0\to\OO(-d)\to\OO^2\to\cI_Z(d)\to0$, with $Z$ a complete intersection of two curves of degree $d$, has a quotient of slope $d$ and a subsheaf of slope $-d$.
\item Two comparable root sequences with the same slope are equal, so $\Psi_F$ is linear iff $\gamma(F)$ is a multiple of a single $\gamma^z$; Corollary~\ref{cor:supernatural}(a) is the case of a linear polygon. If the roots form an arithmetic progression, $z_i=z_1-(i-1)e$, then the bundle $F(z_1+e)$ has no cohomology at the twists $-e,-2e,\dots,-ne$, so it is an Ulrich bundle on $(\PP^n,\OO(e))$, and its semistability also follows from \cite[Thm.~2.9]{CH12}. For arbitrary $z$ the Segre product of Proposition~\ref{prop:sameslope} gives the same reduction: a destabilizing subsheaf $\cG\subset F$ would give a destabilizing subsheaf $\cU_\cG\subset\cU_F$ of the Ulrich bundle $\cU_F$, contradicting \cite[Thm.~2.9(a)]{CH12}; this is a second proof of Corollary~\ref{cor:supernatural}(a), in any characteristic. In characteristic zero, Corollary~\ref{cor:supernatural}(a) can also be deduced from \cite{ErSa16}: for $n+1$ integers containing the roots $z$, the transform $\Phi^W_2$ of Erman and Sam is exact on coherent sheaves \cite[Prop.~4.4]{ErSa16}, multiplies classes in $K_0$ by a positive integer $N_W$ \cite[Cor.~1.4]{ErSa16}, and takes a bundle $F$ with table $q\gamma^z$ to a direct sum of copies of the homogeneous bundle $\mathbf S_\lambda(\mathcal Q^*)(w)$ of \cite[(3.2), Cor.~1.9]{ErSa16}, which is stable, being irreducible homogeneous \cite{Ra66,Um78}; a destabilizing subsheaf of $F$ would give a destabilizing subsheaf of $\Phi^W_2(F)$. Theorem~\ref{thm:main} needs neither the characteristic nor a single ray. Appendix~\ref{app:P3} describes the bundles with supernatural tables on $\PP^3$ with $z_1-z_3\le4$.
\item The proof gives the equality cases of Theorem~\ref{thm:main}(a): for $F/F'$ torsion-free and $L$ general, $\deg F'-\mu(z)N'=\phi_z(\gamma(F))$ iff $h^1(F'\otimes \Tt)=0$, $H^0(F'\otimes \Tt)=H^0(F\otimes \Tt)$, and every differential of the spectral sequence of Lemma~\ref{lem:truncation} from total degree $0$ to total degree $1$ vanishes. Proposition~\ref{prop:marginPn} below expresses the defect in terms of the cohomology of $F'$ and $F/F'$.
\end{enumerate}
\end{remark}

The defect of (a) has an exact expression in terms of the cohomology of $F'$ and $F/F'$ at the $n$ twists $z_1,\dots,z_n$.

\begin{proposition}[the defect]\label{prop:marginPn}
Let $0\to F'\to F\to F''\to0$ be an exact sequence of coherent sheaves on $\PP^n$ with $F$ a vector bundle and $F'\ne0$, $N'=\rk F'$, $\mu'=\mu(F')$, let $\lambda\in\frac1n\ZZ$, and let $z$ be adapted to $\gamma(F)$ at $\lambda$. Write $h^j_t(\cG)=h^j(\cG(t))$, and let $\partial^j_t\colon H^j(F''(t))\to H^{j+1}(F'(t))$ be the connecting homomorphisms. Then
\begin{equation}\label{eq:marginPn}
\Lambda_F(\lambda)-N'(\mu'-\lambda)\ =\ (n-1)!\sum_{r=1}^n\frac{B_r}{P_r(z)},
\end{equation}
where
\[
B_r=\sum_{j=r}^{n}(-1)^{j-r}h^j_{z_r}(F')+\sum_{j=0}^{r-1}(-1)^{r-1-j}h^j_{z_r}(F'')-\rk\partial^{r-1}_{z_r}.
\]
\end{proposition}

\begin{proof}
For a sheaf $\cG$ put $A_r(\cG)=\sum_{j<r}(-1)^jh^j_{z_r}(\cG)$. By Lemma~\ref{lem:adapted},
\[
\Lambda_F(\lambda)=(n-1)!\sum_r\frac{(-1)^{r-1}A_r(F)}{P_r(z)} .
\]
Fix $r$ and $t=z_r$. In the long exact cohomology sequence, $h^j_t(F)=h^j_t(F')+h^j_t(F'')-\rk\partial^j_t-\rk\partial^{j-1}_t$, and the alternating sum over $j<r$ telescopes:
\[
A_r(F)=A_r(F')+A_r(F'')-(-1)^{r-1}\rk\partial^{r-1}_t .
\]
Moreover $A_r(F')=\chi(F'(t))-\sum_{j\ge r}(-1)^jh^j_t(F')$. Hence $(-1)^{r-1}A_r(F)=(-1)^{r-1}\chi(F'(z_r))+B_r$, and $(n-1)!\sum_r(-1)^{r-1}\chi(F'(z_r))/P_r(z)=N'(\mu'-\lambda)$, as in the proof of Theorem~\ref{thm:main}.
\end{proof}

For $n\ge3$ the intermediate cohomology enters $B_r$ with both signs, and the nonnegativity of the right-hand side of \eqref{eq:marginPn}, which Theorem~\ref{thm:main} implies, is not visible term by term. For $n=2$ the $B_r$ can still be negative, but after regrouping the summands the defect becomes a sum of four nonnegative terms: see Proposition~\ref{prop:fourdefect} and the paragraph after it.

\subsection{The gap moment}\label{sec:gap}
In this subsection $n\ge2$. For a root sequence $z=(z_1>\dots>z_n)$ put
\[
\rho(z)=\frac{12}{n(n^2-1)}\sum_{i=1}^n\big(z_i-c(z)\big)^2 ,
\]
the \emph{normalized spread} of $z$; for $n=2$ it is $(z_1-z_2)^2$. For a vector bundle $F$ of rank $N$ on $\PP^n$ with chain decomposition $\gamma(F)=\sum_kq_k\gamma^{z^k}$ let $p_k=n!\,q_k/N$ be the weights of Corollary~\ref{cor:supernatural}(c), and put $c_k=c(z^k)$, $\bar c=\sum_kp_kc_k$ and
\[
G(F)=\sum_kp_k\,\rho(z^k),
\]
the \emph{gap moment} of $F$. Since $\mu(z^k)=-c_k-\frac{n+1}2$, the variance of the Boij--Söderberg slopes in Corollary~\ref{cor:supernatural}(c) is the variance of the centers, $V_c(F)=\sum_kp_k(c_k-\bar c)^2$. For a torsion-free sheaf $\cF$ of rank $N$ on $\PP^n$ let $c_1$ and $\operatorname{ch}_2$ be the degrees of $c_1(\cF)\cdot h^{n-1}$ and $\operatorname{ch}_2(\cF)\cdot h^{n-2}$, $\mu(\cF)=c_1/N$, and $\delta(\cF)=\mu^2-2\operatorname{ch}_2/N$, the degree of $\Delta(\cF)\cdot h^{n-2}/N^2$ with $\Delta=2Nc_2-(N-1)c_1^2$.

\begin{lemma}\label{lem:gapPn}
Let $F$ be a vector bundle of rank $N$ on $\PP^n$, $n\ge2$.
\begin{enumerate}[label=\textup{(\roman*)}]
\item $\displaystyle\frac{\chi(F(t))}N=\frac{t^n}{n!}-\bar c\,\frac{t^{n-1}}{(n-1)!}+\Big(\frac{V_c(F)+\bar c^{\,2}}2-\frac{(n+1)\,G(F)}{24}\Big)\frac{t^{n-2}}{(n-2)!}+O(t^{n-3})$.
\item $\bar c=-\mu(F)-\frac{n+1}2$, and $G(F)-1=\dfrac{12}{n+1}\big(V_c(F)+\delta(F)\big)$.
\item $G(F)\ge1$, with equality iff $F$ has no intermediate cohomology, iff $F$ splits.
\item $\displaystyle\sum_{0<j<n}\sum_th^j(F(t))\ \ge\ C_n\,N\big(G(F)-1\big)$, where $C_n=\dfrac{n\,b_n}{12\,(n-2)!}$ and $b_n=(\lfloor n/2\rfloor-1)!\,(\lceil n/2\rceil-1)!$; so $C_2=\frac16$, $C_3=\frac14$, $C_4=\frac16$. The constant is optimal: equality holds for $F=\Omega^u_{\PP^n}(u)$, $u=\lfloor n/2\rfloor$, and in general iff every element of the chain of $F$ has all gaps $z_u-z_{u+1}$ equal to $1$, or exactly one equal to $2$, in a position $u\in\{\lfloor n/2\rfloor,\lceil n/2\rceil\}$.
\item With the Harder--Narasimhan factors $\cQ_i$, their weights $w_i$ and slopes $\mu_i$, and $V_{HN}(F)$ as in Corollary~\textup{\ref{cor:supernatural}(c)},
\[
G(F)-1=\frac{12}{n+1}\Big(V_c(F)-V_{HN}(F)+\sum_iw_i\,\delta(\cQ_i)\Big).
\]
\end{enumerate}
\end{lemma}

\begin{proof}
(i) The Hilbert polynomial is linear in the table and $N=n!\sum_kq_k$, so $\chi(F(t))/N=\frac1{n!}\sum_kp_k\prod_i(t-z^k_i)$. For a root sequence $z$ with center $c$,
\[
\prod_i(t-z_i)=t^n-nc\,t^{n-1}+e_2(z)\,t^{n-2}-\dots,\qquad e_2(z)=\sum_{i<j}z_iz_j=\frac{n(n-1)}2c^2-\frac12\sum_i(z_i-c)^2 ,
\]
where $\frac12\sum_i(z_i-c)^2=\frac{n(n^2-1)}{24}\rho(z)$; and $\sum_kp_k c_k^2=V_c(F)+\bar c^{\,2}$.

(ii) By Riemann--Roch, $\chi(F(t))=\int\operatorname{ch}(F)e^{th}\operatorname{td}(\PP^n)$ with $\operatorname{td}(\PP^n)=1+\frac{n+1}2h+\frac{(n+1)(3n+2)}{24}h^2+\dots$, so the coefficients of $t^{n-1}/(n-1)!$ and of $t^{n-2}/(n-2)!$ in $\chi(F(t))/N$ are $\mu+\frac{n+1}2$ and $\operatorname{ch}_2/N+\frac{n+1}2\mu+\frac{(n+1)(3n+2)}{24}$. Comparing with (i) gives $\bar c=-\mu-\frac{n+1}2$, and then, since $\bar c^{\,2}=\mu^2+(n+1)\mu+\frac{(n+1)^2}4$,
\[
\frac{V_c+\mu^2}2+\frac{n+1}{24}-\frac{(n+1)G}{24}=\frac{\operatorname{ch}_2}N ,
\]
which is the identity.

(iii) For a root sequence $z$, $n\sum_i(z_i-c)^2=\sum_{i<j}(z_i-z_j)^2\ge\sum_{i<j}(j-i)^2=\frac{n^2(n^2-1)}{12}$, since $z_i-z_j\ge j-i$ for $i<j$; so $\rho(z)\ge1$, with equality iff the roots are consecutive, iff $\gamma^z$ has no intermediate cohomology. The entries of $\gamma(F)=\sum_kq_k\gamma^{z^k}$ are sums of nonnegative terms, so $F$ has no intermediate cohomology iff every $z^k$ is consecutive, iff $G(F)=1$; and a bundle on $\PP^n$ without intermediate cohomology splits by Horrocks' criterion \cite{Ho64}.

(iv) Let $z$ be a root sequence, $g_u=z_u-z_{u+1}=1+a_u$ with integers $a_u\ge0$, $1\le u<n$, and let $I(z)=\sum_{0<j<n}\sum_t\gamma^z_{j,t}$ be its total intermediate cohomology. The left side of (iv) is $\sum_kq_k I(z^k)$, and $\sum_kq_k(\rho(z^k)-1)=N(G(F)-1)/n!$, so it suffices to prove $I(z)\ge n!\,C_n(\rho(z)-1)=\frac{b_nn^2(n-1)}{12}(\rho(z)-1)$. Put $D(z)=\sum_{i<j}\big((z_i-z_j)^2-(j-i)^2\big)$; by the proof of (iii), $D(z)=\frac{n^2(n^2-1)}{12}(\rho(z)-1)$. Writing $z_i-z_j=(j-i)+\sum_{u=i}^{j-1}a_u$ and expanding,
\begin{equation}\label{eq:Dexp}
D(z)=\sum_uu(n-u)\big(na_u+a_u^2\big)+2\sum_{u<v}u(n-v)\,a_ua_v .
\end{equation}
An integer $t$ with $z_{u+1}<t<z_u$ contributes $\prod_j|t-z_j|=(z_u-t)(t-z_{u+1})\prod_{j<u}(z_j-t)\prod_{j>u+1}(t-z_j)$ to $I(z)$, and
\[
z_j-t\ge u-j+1+\sum_{v=j}^{u-1}a_v\quad(j<u),\qquad t-z_j\ge j-u+\sum_{v=u+1}^{j-1}a_v\quad(j>u+1).
\]
Multiplying these lower bounds and keeping the constant and the linear terms, $\prod_{j\ne u,u+1}|t-z_j|\ge B_u\big(1+\sum_{v\ne u}c_{uv}a_v\big)$, where $B_u=u!\,(n-u)!$, $c_{uv}=\sum_{j=v+1}^n\frac1{j-u}$ for $v>u$ and $c_{uv}=\sum_{j=1}^v\frac1{u-j+1}$ for $v<u$. Summing over the gap, $\sum_t(z_u-t)(t-z_{u+1})=f(a_u)$ with $f(a)=\frac{a(a+1)(a+2)}6\ge a$, so
\[
I(z)\ \ge\ \sum_uB_uf(a_u)+\sum_{u<v}\big(B_uc_{uv}+B_vc_{vu}\big)a_ua_v .
\]
Now $B_u=u(n-u)\,(u-1)!\,(n-u-1)!\ge u(n-u)\,b_n$, and for $u<v$ we have $c_{uv}\ge\frac{n-v}{n-u}$ and $c_{vu}\ge\frac uv$, so $B_uc_{uv}+B_vc_{vu}\ge2b_nu(n-v)$. Moreover, for integers $a\ge0$,
\[
f(a)-\frac{na+a^2}{n+1}=\frac{a(a-1)\big((n+1)a+4n-2\big)}{6(n+1)}\ \ge\ 0 .
\]
Comparing with \eqref{eq:Dexp} gives $I(z)\ge\frac{b_n}{n+1}D(z)=\frac{b_nn^2(n-1)}{12}(\rho(z)-1)$. For $F=\Omega^u_{\PP^n}(u)$ the roots are $\{0,-1,\dots,-n\}\setminus\{-u\}$ by Bott's formula, so $a_u=1$ and all other $a_v$ vanish; $\rk F=\binom nu$, the only intermediate cohomology is $H^u(F(-u))=K$, and $G(F)-1=\rho-1=\frac{12u(n-u)}{n^2(n-1)}$. Since $\binom nu u(n-u)=n!/b_n$ for $u=\lfloor n/2\rfloor$, this gives equality. For the equality cases, it suffices by additivity to treat one root sequence. Since $B_uc_{uv}+B_vc_{vu}\ge2b_nu(n-v)>2b_nu(n-v)/(n+1)$, equality forces at most one $a_u>0$; the displayed identity excludes $a_u\ge2$; and for $a_u=1$ the remaining condition $B_u=u(n-u)b_n$ means $(u-1)!(n-u-1)!=b_n$, i.e.\ $u\in\{\lfloor n/2\rfloor,\lceil n/2\rceil\}$. Conversely, if all $a_v$ vanish then $I(z)=D(z)=0$, and if only $a_u=1$ then $I(z)=B_u$ and $D(z)=(n+1)u(n-u)$, so equality holds for central $u$.

(v) The numbers $c_1$ and $\operatorname{ch}_2$ are additive, so $\delta(F)=\sum_iw_i\delta(\cQ_i)+\mu(F)^2-\sum_iw_i\mu_i^2=\sum_iw_i\delta(\cQ_i)-V_{HN}(F)$; insert this in (ii).
\end{proof}

\begin{corollary}\label{cor:Q}
For every vector bundle $F$ on $\PP^n$, $n\ge2$, in characteristic zero, with Harder--Narasimhan factors $\cQ_i$ of weights $w_i$,
\[
G(F)\ \ge\ 1+\frac{12}{n+1}\sum_iw_i\,\delta(\cQ_i),
\]
where $\delta(\cQ_i)\ge0$ by Bogomolov's inequality \cite{Bo78,La04}, \cite[Thm.~7.3.1]{HL10}.
\end{corollary}

\begin{proof}
By Lemma~\ref{lem:gapPn}(v) the inequality is equivalent to $V_{HN}(F)\le V_c(F)$, which is Corollary~\ref{cor:supernatural}(c).
\end{proof}

\section{Majorization and monotonicity}\label{sec:monotone}

\subsection{The functionals as upper bounds}
Applied to supernatural bundles, Theorem~\ref{thm:main}(a) and Lemma~\ref{lem:truncation} give lower bounds for $\phi_z$ on all supernatural tables.

\begin{lemma}\label{lem:positivity}
For all root sequences $y$ and $z$,
\[
\phi_z(\gamma^y)\ \ge\ n!\,\big(c(z)-c(y)\big)_+ ,
\]
with equality if $y$ and $z$ are comparable.
\end{lemma}

\begin{proof}
Let $E$ be a supernatural bundle with root sequence $y$ \cite[Thm.~0.4]{ES09}, so that $\gamma(E)=q\gamma^y$ with $q>0$, $\rk E=n!\,q$ and $\mu(E)=\mu(y)$. Lemma~\ref{lem:truncation} gives $\phi_z(\gamma(E))\ge0$, and Theorem~\ref{thm:main}(a) with $F'=F=E$ gives $\phi_z(\gamma(E))\ge\rk E\,(\mu(y)-\mu(z))=n!\,q\,(c(z)-c(y))$. Divide by $q$. Equality for comparable $y,z$ is Lemma~\ref{lem:comparable}.
\end{proof}

\begin{proposition}\label{prop:minformula}
Let $\gamma\in\CC$ and $\lambda\in\frac1n\ZZ$. For every decomposition $\gamma=\sum_yq'_y\gamma^y$ with $q'_y\ge0$ and every root sequence $z$ with $\mu(z)=\lambda$,
\[
\sum_yq'_y\,n!\,\big(\mu(y)-\lambda\big)_+\ \le\ \phi_z(\gamma).
\]
For the chain decomposition the left-hand side is $\Lambda_\gamma(\lambda)$, and equality holds when $z$ is adapted to $\gamma$ at $\lambda$. Hence
\[
\Lambda_\gamma(\lambda)=\min\{\phi_z(\gamma):\ \mu(z)=\lambda\}\ \ge\ \sum_yq'_y\,n!\,\big(\mu(y)-\lambda\big)_+ .
\]
\end{proposition}

\begin{proof}
The first inequality is Lemma~\ref{lem:positivity} and linearity, since $c(z)-c(y)=\mu(y)-\lambda$. The equality is Lemma~\ref{lem:adapted}.
\end{proof}

So each $\phi_z$ is an upper bound for $\Lambda_\lambda=\Lambda_{(\cdot)}(\lambda)$ on all of $\CC$, and $\Lambda_\lambda$ is their minimum, attained at the root sequences comparable with the chain.

\begin{theorem}\label{thm:concave}
Let $\gamma\in\CC$, and let $\gamma=\sum_jq'_j\gamma^{y_j}$ be any decomposition with $q'_j\ge0$.
\begin{enumerate}[label=\textup{(\alph*)}]
\item The distribution of centers of the chain, with masses $q_k$, majorizes the distribution of the $c(y_j)$, with masses $q'_j$. Equivalently, the slope polygon of the decomposition lies below $\Psi_\gamma$, and $\Lambda_\gamma(\lambda)\ge\sum_jn!\,q'_j(\mu(y_j)-\lambda)_+$ for all $\lambda$.
\item For each $\lambda$ the function $\gamma\mapsto\Lambda_\gamma(\lambda)$ is concave, positively homogeneous and superadditive on $\CC$. For every convex function $f$ on $\frac1n\ZZ$ the function $\gamma\mapsto\sum_kq_kf(c(z^k))$, computed on the chain, is concave; in particular so is $\gamma\mapsto\rk\gamma\cdot V_c(\gamma)$. The function $(\gamma,m)\mapsto\Psi_\gamma(m)$, on $\{(\gamma,m):0\le m\le\rk\gamma\}$, is jointly concave and positively homogeneous; in particular $\Psi_{\gamma+\gamma'}(m+m')\ge\Psi_\gamma(m)+\Psi_{\gamma'}(m')$.
\end{enumerate}
\end{theorem}

\begin{proof}
(a) Both sides of the last inequality are convex and piecewise linear in $\lambda$, with breakpoints in $\frac1n\ZZ$, and they agree for $\lambda$ large; by Proposition~\ref{prop:minformula} the inequality holds on $\frac1n\ZZ$, hence everywhere. Both distributions have the same total mass and mean, because $\rk$ and $\deg$ are linear, so the equivalences are Lemma~\ref{lem:legendre}(ii).

(b) By Proposition~\ref{prop:minformula}, $\Lambda_\lambda$ is a minimum of linear functionals for $\lambda\in\frac1n\ZZ$, hence concave and positively homogeneous, so superadditive; for other $\lambda$ it is a convex combination of its values at the neighboring points of $\frac1n\ZZ$. Let $\gamma,\gamma'\in\CC$ be nonzero (if one of them is zero there is nothing to prove). The chains of all $s\gamma+t\gamma'$ with $s,t\in\QQ_{\ge0}$, not both zero, lie between the termwise minimum of the last elements and the termwise maximum of the first elements of the chains of $\gamma$ and $\gamma'$ (Lemma~\ref{lem:ends}), so only finitely many centers occur on the segment. On a finite set of centers, a convex function $f$ on $\frac1n\ZZ$ is an affine function plus a nonnegative combination of the functions $x\mapsto(a-x)_+$, $a\in\frac1n\ZZ$, and $\sum_kq_k(a-c(z^k))_+=\Lambda_\gamma(\lambda_a)/n!$ with $c_{\lambda_a}=a$, while $\sum_kq_k$ and $\sum_kq_kc(z^k)$ are linear in $\gamma$. Further, with $\rk\gamma\cdot V_c(\gamma)=0$ for $\gamma=0$, $\rk\gamma\cdot V_c(\gamma)=n!\sum_kq_kc(z^k)^2-\big(n!\sum_kq_kc(z^k)\big)^2/\rk\gamma$, and the second term is convex in $\gamma$, being a square of a linear function divided by a linear function. Finally, by Lemma~\ref{lem:legendre}(i), $\Psi_\gamma(m)$ is the minimum over $\lambda$ of $m\lambda+\Lambda_\gamma(\lambda)$, which is jointly concave and positively homogeneous in $(\gamma,m)$.
\end{proof}

For fixed $m>0$, $\gamma\mapsto\Psi_\gamma(m)$ is concave but neither homogeneous nor superadditive: $\Psi_\gamma(1)=\Psi_{2\gamma}(1)=1$ for $\gamma=\gamma(\OO(1)\oplus\OO(-1))$.

\begin{example}\label{ex:sum}
On $\PP^2$, let $F=\pi_*\OO_{\PP^1\times\PP^1}(-4,-1)\oplus\OO(-3)^{\oplus2}$, where $\pi$ is a finite linear projection of the quadric \cite[Thm.~6.1]{ES09}. Then $\gamma(F)=\gamma^{(3,0)}+\gamma^{(2,1)}$, with incomparable root sequences, and the chain decomposition is $\gamma(F)=\frac23\big(\gamma^{(3,1)}+\gamma^{(3,0)}+\gamma^{(2,0)}\big)$. The given decomposition has centers $\frac32,\frac32$; the chain has centers $2,\frac32,1$, each with mass $\frac23$, a mean-preserving spread.
\end{example}

\begin{proposition}[decomposition criterion]\label{prop:splitting}
Let $\gamma\in\CC$, $m\in\QQ$ with $0\le m\le\rk\gamma$, and $d\in\QQ$. Then $\Psi_\gamma(m)\ge d$ if and only if
\[
\gamma=x+y\qquad\text{with }x,y\in\CC,\quad \rk x=m,\quad \deg x\ge d .
\]
\end{proposition}

\begin{proof}
Let $\gamma=\sum_kq_k\gamma^{z^k}$ be the chain decomposition. If $\Psi_\gamma(m)\ge d$, choose rational $0\le x_k\le q_k$ with $\sum_kn!\,x_k=m$ and $\sum_kn!\,x_k\mu(z^k)=\Psi_\gamma(m)$, and put $x=\sum_kx_k\gamma^{z^k}$, $y=\gamma-x$. Conversely, decompositions of $x$ and of $y$ together form a decomposition of $\gamma$; by Theorem~\ref{thm:concave}(a) its slope polygon lies below $\Psi_\gamma$, and taking all the pieces of $x$ shows that this polygon is at least $\deg x\ge d$ at $m$.
\end{proof}

So $\Psi_\gamma(m)=\max\{\deg x:\ \gamma=x+y,\ x,y\in\CC,\ \rk x=m\}$ for $m\in\QQ$: the polygon $\Psi_\gamma$ is determined by the cone $\CC$ alone, in the same way as the Harder--Narasimhan polygon $P_F$ is the upper concave envelope of the points $(\rk F',\deg F')$ of the subsheaves $F'\subset F$. In these terms Theorem~\ref{thm:main} says that for every subsheaf $F'\subset F$ there is a decomposition $\gamma(F)=x+y$ in $\CC$ with $\rk x=\rk F'$ and $\deg x\ge\deg F'$.

\subsection{Monotonicity}
A table $\Delta$ is \emph{$\chi$-neutral} if $\sum_j(-1)^j\Delta_{j,t}=0$ for all $t$. Two tables in $\CC$ have the same Hilbert polynomial iff their difference is $\chi$-neutral. The cancellation vectors $e_{j,t}+e_{j+1,t}$ are $\chi$-neutral, and by the long exact sequences, for every exact sequence $0\to F'\to F\to F''\to0$,
\begin{equation}\label{eq:extdiff}
\gamma(F')+\gamma(F'')-\gamma(F)=\sum_t\sum_{j=0}^{n-1}\rk\big(\partial^j_t\big)\big(e_{j,t}+e_{j+1,t}\big),\qquad \partial^j_t\colon H^j(F''(t))\to H^{j+1}(F'(t)).
\end{equation}
On $\PP^2$ every nonnegative $\chi$-neutral table is a nonnegative combination of cancellation vectors, since a column $(a,a+b,b)$ is $a(1,1,0)+b(0,1,1)$. For $n\ge3$ this fails: $e_{0,t}+e_{3,t}$ is nonnegative and $\chi$-neutral, but not a nonnegative combination of cancellation vectors.

\begin{lemma}[local values]\label{lem:localvalues}
Let $z$ be a root sequence, $1\le r\le n$ and $t\in\ZZ$. Then
\[
\phi_z\big(e_{r-1,t}+e_{r,t}\big)=v_{r,t}(z)=\begin{cases}\dfrac{(n-1)!}{P_r(z)} & \text{if } z_r=t,\\[2ex] 0 & \text{otherwise.}\end{cases}
\]
\end{lemma}

\begin{proof}
By \eqref{eq:phi}, $\phi_z(e_{j,t})=(n-1)!(-1)^j/D_s(z)$ if $t=z_s$ and $j<s$, and $\phi_z(e_{j,t})=0$ otherwise. If $t=z_s$ with $s>r$, the contributions of $j=r-1$ and $j=r$ cancel; if $s=r$, only $j=r-1$ contributes, with $(-1)^{r-1}/D_r(z)=1/P_r(z)$; if $s<r$, neither contributes.
\end{proof}

\begin{theorem}[monotonicity]\label{thm:monotone}
Let $\gamma,\gamma'\in\CC$, and suppose that $\gamma'-\gamma$ is a nonnegative combination of cancellation vectors. Then $\Lambda_\gamma(\lambda)\le\Lambda_{\gamma'}(\lambda)$ for every $\lambda$; equivalently $\Psi_\gamma\le\Psi_{\gamma'}$, i.e.\ the center distribution of the chain of $\gamma'$ majorizes that of the chain of $\gamma$. On $\PP^2$ the hypothesis says that $\gamma\le\gamma'$ entrywise and that $\gamma,\gamma'$ have the same Hilbert polynomial.
\end{theorem}

\begin{proof}
Let $\lambda\in\frac1n\ZZ$ and let $z$ be adapted to $\gamma'$ at $\lambda$. By Lemma~\ref{lem:adapted}, Lemma~\ref{lem:localvalues} and Proposition~\ref{prop:minformula},
\[
\Lambda_{\gamma'}(\lambda)=\phi_z(\gamma')=\phi_z(\gamma)+\phi_z(\gamma'-\gamma)\ \ge\ \phi_z(\gamma)\ \ge\ \Lambda_\gamma(\lambda).
\]
Both sides are piecewise linear with breakpoints in $\frac1n\ZZ$, so the inequality holds for all $\lambda$. Since $\gamma'-\gamma$ is $\chi$-neutral, $\gamma$ and $\gamma'$ have the same rank and degree, and the equivalence with $\Psi_\gamma\le\Psi_{\gamma'}$ is Lemma~\ref{lem:legendre}(ii). On $\PP^2$ the hypothesis is equivalent to the stated one by the remark before the lemma.
\end{proof}

\begin{corollary}\label{cor:monoext}
\begin{enumerate}[label=\textup{(\alph*)}]
\item For every exact sequence $0\to F'\to F\to F''\to0$ of vector bundles on $\PP^n$, $\Psi_F\le\Psi_{F'\oplus F''}$. Hence $\Psi_F\le\Psi_{\gr F}$ for every filtration of $F$ with locally free subquotients, and $V_c(F)\le V_c(\gr F)$.
\item Let $(F_s)_{s\in S}$ be a flat family of vector bundles on $\PP^n$ over a variety $S$. Every $s_0\in S$ has an open neighborhood $U$ with $\Psi_{F_s}\le\Psi_{F_{s_0}}$ for all $s\in U$. If $S$ is irreducible, then $\Psi_{F_s}\le\Psi_{F_{s_0}}$ for all $s$ in a dense open subset and every $s_0\in S$.
\item Let $F',F''$ be vector bundles and let $F$ be an extension of $F''$ by $F'$ with general class in $\Ext^1(F'',F')$. Then $\Psi_F\le\Psi_{F_1}$ for every extension $F_1$ of $F''$ by $F'$.
\end{enumerate}
\end{corollary}

\begin{proof}
(a) By \eqref{eq:extdiff}, $\gamma(F'\oplus F'')-\gamma(F)$ is a nonnegative combination of cancellation vectors; apply Theorem~\ref{thm:monotone}. This relation between tables is transitive and is preserved by adding the same table to both sides, so we may induct along the filtration. The center distributions of $F$ and $\gr F$ have the same mean, and the second majorizes the first; applying the convex function $x\mapsto x^2$ gives $V_c(F)\le V_c(\gr F)$.

(b) We may shrink $S$ to an affine neighborhood $\operatorname{Spec}A$ of $s_0$. Being $m$-regular is an open condition on $S$, since it is the vanishing of finitely many cohomology groups, which are upper semicontinuous in $s$; so after shrinking $S$ the fibers $F_s$ and $F_s^\vee$ are $m$-regular for a fixed $m$. This gives a bounded range of $t$ outside which the column $t$ of $\gamma(F_s)$ has a single nonzero row, equal to $\pm\chi(F_s(t))$, which does not depend on $s$. Fix $t$ in this bounded range. There is a complex $K^0\to K^1\to\dots\to K^n$ of finitely generated projective, hence after shrinking $S$ free, $A$-modules computing $H^j(F_s(t))$ for every $s$, after base change to $s$ \cite[Ch.~II, \S5]{Mu70}. With $\rho_j(s)$ the rank of $K^j\otimes k(s)\to K^{j+1}\otimes k(s)$, we get $h^j(F_s(t))=\rk K^j-\rho_j(s)-\rho_{j-1}(s)$. Each $\rho_j$ is lower semicontinuous, so $\rho_j(s)\ge\rho_j(s_0)$ for all $j$ and all $s$ in an open neighborhood $U_t$ of $s_0$, and for these $s$
\[
\gamma(F_{s_0})-\gamma(F_s)=\sum_{j=0}^{n-1}\big(\rho_j(s)-\rho_j(s_0)\big)\big(e_{j,t}+e_{j+1,t}\big)\quad\text{in the column }t .
\]
On the intersection $U$ of the $U_t$, $t$ in the bounded range, $\gamma(F_{s_0})-\gamma(F_s)$ is a nonnegative combination of cancellation vectors, and Theorem~\ref{thm:monotone} applies. If $S$ is irreducible, then, $S$ being quasi-compact, finitely many open sets as above cover it, so there is one $m$ for which all $F_s$ and $F_s^\vee$, $s\in S$, are $m$-regular; the locus $V$ where the finitely many $h^j(F_s(t))$ with $t$ in the resulting bounded range take their minimal values is dense and open, $\gamma(F_s)$ is constant on $V$, and $V$ meets every neighborhood $U$ as above. (c) Apply (b) to the family of extensions over $\Ext^1(F'',F')$. The proof of (b) works for a vector bundle on $\PP^n_S$ over any noetherian scheme $S$, with the tables, and their chain decompositions, of the geometric fibers, over the algebraic closures of the residue fields; for instance, over the spectrum of a discrete valuation ring of mixed characteristic it gives $\Psi_{F_\eta}\le\Psi_{F_s}$ for the generic and the special fiber.
\end{proof}

The same computation bounds how much cancellation can lower $\Lambda$.

\begin{proposition}[defect inequality]\label{prop:defectineq}
Let $\gamma,\gamma'\in\CC$ with
\[
\gamma'-\gamma=\sum_{r=1}^n\sum_td^r_t\,\big(e_{r-1,t}+e_{r,t}\big),\qquad d^r_t\ge0 ,
\]
let $\lambda\in\frac1n\ZZ$, and let $z$ be adapted to $\gamma$ at $\lambda$. Then
\begin{equation}\label{eq:defect}
\Lambda_\gamma(\lambda)\ \ge\ \Lambda_{\gamma'}(\lambda)-\sum_{r=1}^nd^r_{z_r}\,\frac{(n-1)!}{P_r(z)}.
\end{equation}
If $z'$ is adapted to $\gamma'$ at $\lambda$, then also $\Lambda_{\gamma'}(\lambda)-\Lambda_\gamma(\lambda)\ge\sum_{r=1}^nd^r_{z'_r}\,(n-1)!/P_r(z')$.
\end{proposition}

\begin{proof}
By Lemma~\ref{lem:adapted}, Proposition~\ref{prop:minformula} and Lemma~\ref{lem:localvalues}, $\Lambda_\gamma(\lambda)=\phi_z(\gamma)=\phi_z(\gamma')-\phi_z(\gamma'-\gamma)\ge\Lambda_{\gamma'}(\lambda)-\sum_{r,t}d^r_t\,v_{r,t}(z)$. Similarly $\Lambda_{\gamma'}(\lambda)=\phi_{z'}(\gamma)+\phi_{z'}(\gamma'-\gamma)\ge\Lambda_\gamma(\lambda)+\sum_{r,t}d^r_t\,v_{r,t}(z')$.
\end{proof}

For an exact sequence $0\to F'\to F\to F''\to0$ of vector bundles, \eqref{eq:extdiff} gives $d^r_t=\rk\partial^{r-1}_t$. So Theorem~\ref{thm:monotone} and \eqref{eq:defect} bound $\Lambda_F$ from both sides: cancellation lowers $\Lambda_F(\lambda)$ below $\Lambda_{F'\oplus F''}(\lambda)$, but at most by the ranks of the connecting homomorphisms at the $n$ twists $z_1,\dots,z_n$, weighted by $(n-1)!/P_r(z)$.

For $n\ge3$ the hypothesis of Theorem~\ref{thm:monotone} cannot be weakened to entrywise comparison with the same Hilbert polynomial. The local value of $e_{0,t}+e_{3,t}$ at a root sequence $z$ with $z_2=t$ is $\phi_z(e_{0,t})=(n-1)!/D_2(z)<0$, and this already happens for split bundles.

\begin{example}[entrywise comparison is not enough]\label{ex:weak}
On $\PP^3$ let $F=\OO(-4)^{6}\oplus\OO(-2)^{6}$ and $\cG=\OO(-5)\oplus\OO(-4)^{2}\oplus\OO(-3)^{6}\oplus\OO(-2)^{2}\oplus\OO(-1)$. In the Grothendieck group $\cG-F=\OO(-5)-4\,\OO(-4)+6\,\OO(-3)-4\,\OO(-2)+\OO(-1)$, and the Hilbert polynomials agree because the fourth difference of the cubic $a\mapsto\chi(\OO(a+t))$ vanishes. Moreover $\gamma(\cG)-\gamma(F)=e_{0,1}+e_{3,1}$: in the column $1$ the table of $F$ vanishes, while $h^0(\OO(-1)(1))=h^3(\OO(-5)(1))=1$, and in all other columns the two tables agree. So $\gamma(F)\le\gamma(\cG)$ entrywise, with the same Hilbert polynomial. For split bundles $\Psi=P$, and $\Psi_F(6)=-12>-14=\Psi_\cG(6)$, while $\Psi_F(1)=-2<-1=\Psi_\cG(1)$. So the polygons are not comparable. In terms of chains, $\gamma(F)=\gamma^{(3,2,1)}+\gamma^{(1,0,-1)}$ and $\gamma(\cG)=\frac16\gamma^{(4,3,2)}+\frac13\gamma^{(3,2,1)}+\gamma^{(2,1,0)}+\frac13\gamma^{(1,0,-1)}+\frac16\gamma^{(0,-1,-2)}$.
\end{example}

The same happens on every $\PP^n$, $n\ge3$: appending the tail $(-3,-4,\dots,-(n-1))$ to each root sequence of Example~\ref{ex:weak}, with the same coefficients, gives tables $\gamma\le\gamma'$ in $\CC$ with the same Hilbert polynomial and $\gamma'-\gamma=\frac{n!}6(e_{0,1}+e_{3,1})$, whose polygons are not comparable: $\Psi_\gamma-\Psi_{\gamma'}$ is $-\frac{(n-1)!}2$ at $m=\frac{n!}6$ and $(n-1)!$ at $m=n!$.

Theorem~\ref{thm:monotone} has no analogue for the Harder--Narasimhan polygon as a function of the table, since $P_F$ is not determined by $\gamma(F)$. For instance, let $\pi\colon\PP^1\times\PP^1\to\PP^2$ be the projection of a smooth quadric surface from a point off it, so $\pi^*\OO(1)=\OO(1,1)$, $\cG=\pi_*\OO(-1,3)$ and $H=\cG^\vee$. Then $\gamma(\cG)=\gamma^{(0,-4)}$ by the Künneth formula, $\cG$ is stable (semistable by Corollary~\ref{cor:supernatural}(a), of rank $2$ and odd degree), $\gamma(H)=\gamma^{(1,-3)}$, and $\Ext^1(\cG,H)$ has dimension $-\chi(\cG,H)=15$. For every extension $0\to H\to E\to \cG\to0$ the connecting maps vanish, since $H^0(\cG(t))$ and $H^1(H(t))$, and likewise $H^1(\cG(t))$ and $H^2(H(t))$, are never both nonzero; so $\gamma(E)=\gamma(\cG\oplus H)$. But a nonsplit $E$ is semistable: $H^0(E)=H^0(E^\vee)=0$ excludes destabilizing subsheaves of rank $1$ and $3$, and a saturated subsheaf $A\subset E$ of rank $2$ with $c_1(A)>0$ would split the extension. Indeed, $A$ is reflexive, hence locally free. If $A\cap H$ had rank $2$, then $A=H$ by saturation, but $c_1(H)=-1$; if it had rank $1$, then $c_1(A\cap H)\le-1$ and the image of $A$ in $\cG$ has $c_1\le0$, by the stability of $H$ and $\cG$, so $c_1(A)\le-1$. Hence $A\to\cG$ is injective, and $c_1(A)\ge1=c_1(\cG)$ shows that its determinant is a nonzero constant; so $A\to\cG$ is an isomorphism. So $P_E$ is a segment, while $P_{\cG\oplus H}$ has slopes $\pm\frac12$. For extensions and specializations, however, $P_F\le P_{F'\oplus F''}$ and Shatz's semicontinuity give the analogous monotonicity \cite{Sh77}. Nor is the chain itself monotone: on $\PP^2$, $\gamma=\gamma^{(4,3)}+\gamma^{(4,2)}$ and $\gamma'=\gamma+e_{1,2}+e_{2,2}=\frac32\gamma^{(4,3)}+\frac12\gamma^{(4,1)}$, and neither chain contains the other.

\subsection{Coherent sheaves}\label{sec:coherent}
Let $\cF$ be a coherent sheaf on $\PP^n$, with $\rk\cF\ge0$, $\deg\cF=c_1(\cF)\cdot h^{n-1}$ and torsion subsheaf $\cF_{\rm tors}$; $\deg\cF_{\rm tors}\ge0$ is the degree of the $(n-1)$-dimensional part of $\cF_{\rm tors}$. Following Proposition~\ref{prop:minformula}, put
\[
\Lambda_\cF(\lambda)=\inf\big\{\phi_z(\gamma(\cF)):\ \mu(z)=\lambda\big\},\ \ \lambda\in\tfrac1n\ZZ,\qquad \Psi_\cF(m)=\inf_{\lambda\in\frac1n\ZZ}\big(m\lambda+\Lambda_\cF(\lambda)\big),\ \ 0\le m\le\rk\cF .
\]
For vector bundles these are the functions defined before, by Proposition~\ref{prop:minformula} and Lemma~\ref{lem:legendre}(i). Let $P_\cF$ be the upper concave envelope of the points $(\rk\cF',\deg\cF')$ for coherent subsheaves $\cF'\subset\cF$; it is the Harder--Narasimhan polygon of $\cF/\cF_{\rm tors}$ raised by $\deg\cF_{\rm tors}$.

\begin{proposition}[Theorem~A for coherent sheaves]\label{prop:cohA}
Let $\cF$ be a coherent sheaf on $\PP^n$ and $\cF'\subset\cF$ a coherent subsheaf.
\begin{enumerate}[label=\textup{(\alph*)}]
\item $\deg\cF'-\mu(z)\rk\cF'\le\phi_z(\gamma(\cF))$ for every root sequence $z$. Hence $\deg\cF'\le\Psi_\cF(\rk\cF')$, i.e.\ $P_\cF\le\Psi_\cF$.
\item $\Lambda_\cF(\lambda)\ge\max(\deg\cF_{\rm tors},\deg\cF-\lambda\rk\cF)$ for all $\lambda\in\frac1n\ZZ$, and $\Lambda_\cF(\lambda)=\deg\cF-\lambda\rk\cF$ for $\lambda\ll0$. So $\Psi_\cF$ is concave, $\Psi_\cF(0)\ge\deg\cF_{\rm tors}$ and $\Psi_\cF(\rk\cF)=\deg\cF$.
\end{enumerate}
\end{proposition}

\begin{proof}
(a) For $\cF'=0$ the inequality says that $\phi_z(\gamma(\cF))\ge0$, which is a case of the positivity theorem of Eisenbud and Schreyer for coherent sheaves \cite[Thm.~4.1]{ES10}; as in its proof, we pass to a general translate to make the sheaves homologically transverse. Let $\Tt$ be the sheaf of Lemma~\ref{lem:testsheaf}. By the Kleiman--Bertini theorem for tensor products of Miller and Speyer \cite{MiSp08}, which holds in any characteristic, $\operatorname{Tor}_i(\mathcal G,g_*\Tt)=0$ for $i>0$ and $\mathcal G=\cF',\cF,\cF/\cF'$, for a general $g\in\mathrm{PGL}_{n+1}(K)$; and $g_*\Tt$ has the same resolution as $\Tt$. So $\mathcal G$ tensored with this resolution is a resolution of $\mathcal G\otimes g_*\Tt$, and $\cF'\otimes g_*\Tt\to\cF\otimes g_*\Tt$ is injective. With $g_*\Tt$ in place of $\Tt$, and without saturating $\cF'$, the proofs of Lemma~\ref{lem:truncation} and Theorem~\ref{thm:main}(a) apply verbatim, with Riemann--Roch in the form of Section~\ref{sec:tables}, where the rank may be zero. This gives the first inequality, and the second follows by taking the infimum over $z$ and then over $\lambda$. For torsion-free $\cF$ the original proofs apply as they stand, with a general line $L$ that is also disjoint from the locus where $\cF$ is not locally free, and \cite{MiSp08} is not needed.

(b) The lower bound is (a) for $\cF'=\cF_{\rm tors}$ and for $\cF'=\cF$. If $z_n$ is at least the regularity of $\cF$, the higher cohomology of $\cF$ vanishes at the twists $z_r$, and the Riemann--Roch computation in the proof of Theorem~\ref{thm:main} gives $\phi_z(\gamma(\cF))=\deg\cF-\mu(z)\rk\cF$; this gives the equality for $\lambda\ll0$. The function $\Psi_\cF$ is concave as an infimum of affine functions, and the values at $0$ and $\rk\cF$ follow.
\end{proof}

\begin{example}\label{ex:Ip}
Let $\cI_p\subset\OO_{\PP^2}$ be the ideal sheaf of a point. In the rows $0$ and $1$ of its table, $h^0(\cI_p(t))=\binom{t+2}2-1$ for $t\ge0$ and $h^1(\cI_p(t))=1$ for $t<0$, and the other entries vanish. So for $a>b$ with $a+b=-2\lambda-3$, $\phi_{(a,b)}(\gamma(\cI_p))$ equals $-\lambda$ if $b\ge0$, $\frac{(a+1)(a+2)}{2(a-b)}$ if $b<0\le a$, and $\frac1{a-b}$ if $a<0$. Minimizing, with the lower bound $-\lambda$ of Proposition~\ref{prop:cohA}(b),
\[
\Lambda_{\cI_p}(\lambda)=\begin{cases}-\lambda, & \lambda\le-\frac12,\\[0.5ex] \dfrac1{2\lambda+3}, & \lambda\ge-\frac12,\end{cases}\qquad \Psi_{\cI_p}(m)=\inf_{x\in\ZZ,\,x\ge2}\Big(\frac{x-3}2\,m+\frac1x\Big).
\]
So $\Psi_{\cI_p}$ has infinitely many vertices $\big(\frac2{x(x+1)},\frac{2(x-1)}{x(x+1)}\big)$, $x\ge1$, accumulating at $(0,0)$, with an edge of slope $\frac{x-3}2$ between the vertices for $x$ and $x-1$, for $x\ge2$; $\Psi_{\cI_p}(0)=0$ and $\Psi_{\cI_p}(m)=\sqrt{2m}+O(m)$ as $m\to0$, while $P_{\cI_p}=0$.
\end{example}

By \cite[Thm.~0.2]{ES10}, $\gamma(\cF)=\sum_{y\in Y}q_y\gamma^y$, entrywise convergent, for a unique chain $Y$ and unique $q_y>0$. Here $y=(y_1>\dots>y_s)$, $0\le s\le n$, is completed by $y_i=-\infty$ for $i>s$, and $\gamma^y$ is the table of a supernatural sheaf of dimension $s$ with Hilbert polynomial $\prod_{i\le s}(t-y_i)$; such sheaves exist for every $y$ \cite{ES10}. We call $Y$ the chain of $\cF$, and we write $s=s(y)$.

\begin{proposition}[the series of Eisenbud and Schreyer]\label{prop:cohseries}
Let $\cF$ be a coherent sheaf on $\PP^n$ with chain $Y$, and put $\zeta_\cF=(n-1)!\sum_{y:\,s=n-1}q_y$.
\begin{enumerate}[label=\textup{(\alph*)}]
\item For every $\lambda\in\frac1n\ZZ$ there is a root sequence $z$ with $\mu(z)=\lambda$ comparable with every element of $Y$, and for every such $z$
\[
\Lambda_\cF(\lambda)=\phi_z(\gamma(\cF))=\sum_{y:\,s=n}n!\,q_y\big(\mu(y)-\lambda\big)_++\zeta_\cF .
\]
\item $\sum_{s=n}n!\,q_y=\rk\cF$ and $\zeta_\cF+\sum_{s=n}n!\,q_y\mu(y)=\deg\cF$. The graph of $\Psi_\cF$ starts at $(0,\zeta_\cF)$, where $\zeta_\cF=\Psi_\cF(0)\ge\deg\cF_{\rm tors}$, and has edges of slopes $\mu(y)$ and horizontal lengths $n!\,q_y$, $s=n$, in decreasing order of slope, possibly with infinitely many edges accumulating at the left endpoint.
\item If $\cF$ is torsion-free, only $s=n$ occurs and $\zeta_\cF=0$. The chain $Y$ is finite with $s=n$ throughout if and only if $\cF$ is locally free.
\end{enumerate}
\end{proposition}

\begin{proof}
The assertions are immediate for $\cF=0$; assume $\cF\ne0$. (a) Since $\phi_z$ involves finitely many entries, $\phi_z(\gamma(\cF))=\sum_yq_y\phi_z(\gamma^y)$. Let $s<n$. If $y\le z$, the computation in the proof of Lemma~\ref{lem:comparable}, with \eqref{eq:divdiff} applied to the polynomial $\prod_{i\le s}(t-y_i)$ of degree $s$, gives $\phi_z(\gamma^y)=(n-1)!$ for $s=n-1$ and $\phi_z(\gamma^y)=0$ for $s\le n-2$. For arbitrary $z$, Proposition~\ref{prop:cohA}(a) applied to a supernatural sheaf $\mathcal S$ of dimension $s$, with $\cF'=\mathcal S$ if $s=n-1$ and $\cF'=0$ if $s\le n-2$, gives $\phi_z(\gamma^y)\ge(n-1)!$, resp.\ $\phi_z(\gamma^y)\ge0$. The elements of $Y$ with $s=n$ lie above the others, and $Y$ is bounded above, since the rows $j\ge1$ of $\gamma(\cF)$ vanish in large degrees; so only finitely many elements of $Y$ with $s=n$ have center $\ge c_\lambda$. Hence, as in Lemma~\ref{lem:adapted}, there is a root sequence $z$ with $\mu(z)=\lambda$ comparable with every element of $Y$: let $x$ be the smallest element of $Y$ with $s=n$ and $c(x)\ge c_\lambda$. If $x$ is followed in $Y$ by an element $x'$ with $s=n$, take $z$ in $[x',x]$ as in Lemma~\ref{lem:adapted}; if $x$ exists and is not followed by such an element, take $z=x-n(c(x)-c_\lambda)\epsilon_n$; and if there is no such $x$, take $z$ above the largest element $y$ of $Y$: if $y$ has $s=n$, take $z=y+n(c_\lambda-c(y))\epsilon_1$; if $y$ has $s<n$, choose integers $z_1>\dots>z_{n-1}$ with $z_i\ge y_i$ for $i\le s$, so large that $z_n=nc_\lambda-\sum_{i<n}z_i$ is smaller than $z_{n-1}$ (for $n=1$ take $z=(c_\lambda)$). For a comparable $z$ the displayed value of $\phi_z(\gamma(\cF))$ follows from Lemma~\ref{lem:comparable} and the values above, and by Lemma~\ref{lem:positivity} and the lower bounds above it is a lower bound for $\phi_{z'}(\gamma(\cF))$ for every $z'$ with $\mu(z')=\lambda$.

(b) The sum in (a) converges, since $\Lambda_\cF(\lambda)<\infty$. The slopes $\mu(y)$, $s=n$, are bounded below. Choose $\lambda_0\in\frac1n\ZZ$ below all of them and in the range where $\Lambda_\cF(\lambda)=\deg\cF-\lambda\rk\cF$ (Proposition~\ref{prop:cohA}(b)). Subtracting the formulas at $\lambda_0$ and $\lambda_0-1/n$ gives $\sum_{s=n}n!\,q_y=\rk\cF$, and substitution gives $\zeta_\cF+\sum_{s=n}n!\,q_y\mu(y)=\deg\cF$. These series have finite absolute first moment. For $\lambda\ge\lambda_0$, the summands $n!\,q_y(\mu(y)-\lambda)_+$ are dominated by the summable terms $n!\,q_y(\mu(y)-\lambda_0)$; hence their sum tends to $0$ as $\lambda\to\infty$, and $\Psi_\cF(0)=\inf_\lambda\Lambda_\cF(\lambda)=\zeta_\cF$, which is at least $\deg\cF_{\rm tors}$ by Proposition~\ref{prop:cohA}(b). The description of the graph follows from (a) by the duality of Lemma~\ref{lem:legendre}(i). One may adjoin the vertical segment from $(0,0)$ to $(0,\zeta_\cF)$ to obtain an augmented polygon.

(c) If $\cF$ is torsion-free, only $s=n$ occurs and the $q_y$ are rational \cite[Thm.~0.2]{ES10}. If $Y$ is finite with $s=n$ throughout, then $H^j(\cF(t))=0$ for $j<n$ and $t\ll0$, and then $\mathcal{E}xt^i(\cF,\OO)=0$ for $i>0$: Serre duality gives $\operatorname{Ext}^i(\cF,\OO(k))=0$ for $k\gg0$, and the local-to-global Ext spectral sequence, together with Serre vanishing, identifies these groups with $H^0(\mathcal{E}xt^i(\cF,\OO)(k))$. Vanishing for all large $k$ forces the Ext sheaves to vanish, so $\cF$ is locally free. The converse is Theorem~\ref{thm:ES}.
\end{proof}

By Example~\ref{ex:Ip} and Proposition~\ref{prop:cohseries}(b), the chain of $\cI_p$ is infinite. It is computed in \cite[Example~0.3]{ES10}: $\gamma(\cI_p)=\sum_{k\ge2}\frac2{(k-1)k(k+1)}\,\gamma^{(0,-k)}$. The piece with root sequence $(0,-k)$ has slope $\frac{k-3}2$ and rank $\frac4{(k-1)k(k+1)}$, in agreement with the edges of $\Psi_{\cI_p}$ found in Example~\ref{ex:Ip}.

\begin{corollary}\label{cor:ESquestion}
Let $\cF$ be a coherent sheaf on $\PP^n$, $n\ge1$, such that all the tables in the decomposition $\gamma(\cF)=\sum_{y\in Y}q_y\gamma^y$ of \cite[Thm.~0.2]{ES10} are tables of vector bundles on $\PP^n$, i.e.\ $s=n$ for every $y\in Y$ (the criterion of \cite{ES10} for $\gamma^y$ to be the table of a vector bundle). Then $\cF$ has no associated subvarieties of dimension $n-1$ or $0$. In particular, for $n\le2$ the sheaf $\cF$ is torsion-free.
\end{corollary}

\begin{proof}
By Proposition~\ref{prop:cohseries}(b), $\zeta_\cF=0$ and $\zeta_\cF\ge\deg\cF_{\rm tors}\ge0$; so $\deg\cF_{\rm tors}=0$, i.e.\ $\cF_{\rm tors}$ has no component of dimension $n-1$. Let $\cF_0\subset\cF$ be the largest subsheaf of finite length, and $\ell$ its length. Then $h^0(\cF_0(t))=\ell$ and $H^i(\cF_0(t))=0$ for $i>0$ and all $t$, so $\gamma(\cF)=\gamma(\cF/\cF_0)+\ell\sum_te_{0,t}$, and $\sum_te_{0,t}$ is the supernatural table with $s=0$, the table of a point. Its sequence lies below every other one, so adding it to the chain of $\cF/\cF_0$ gives a chain, which for $\ell>0$ is the chain of $\cF$ by the uniqueness in \cite[Thm.~0.2]{ES10}. By hypothesis the chain of $\cF$ has no element with $s=0$, so $\ell=0$. Finally, the associated subvarieties of $\cF$ of dimension less than $n$ are those of $\cF_{\rm tors}$, and for $n\le2$ they have dimension $n-1$ or $0$.
\end{proof}

Eisenbud and Schreyer observe that the decomposition of a torsion-free sheaf involves only tables of vector bundles, and ask whether the converse is true \cite[after Thm.~0.2]{ES10}. Corollary~\ref{cor:ESquestion} shows that it is for $n\le2$; for $n\ge3$ it leaves open the associated subvarieties of dimension $1,\dots,n-2$.

\begin{proposition}[Theorem~C for coherent sheaves]\label{prop:cohC}
\begin{enumerate}[label=\textup{(\alph*)}]
\item For every exact sequence $0\to\cF'\to\cF\to\cF''\to0$ of coherent sheaves on $\PP^n$, $\Lambda_\cF\le\Lambda_{\cF'\oplus\cF''}$ and $\Psi_\cF\le\Psi_{\cF'\oplus\cF''}$. Hence $\Psi_\cF\le\Psi_{\gr\cF}$ for every filtration of $\cF$ by coherent subsheaves, for instance for the Harder--Narasimhan filtration of a vector bundle, whose factors need not be locally free.
\item If $n\ge2$ and $\cF\subset\cE$ are coherent sheaves with $\cE/\cF$ of finite length, then $\Psi_\cE\le\Psi_\cF$.
\item Let $(\cF_s)_{s\in S}$ be a flat family of coherent sheaves on $\PP^n$ over a variety $S$, $s_0\in S$ and $m_0>0$. There is a neighborhood $U$ of $s_0$ with $\Psi_{\cF_s}\le\Psi_{\cF_{s_0}}$ on $[m_0,\rk\cF_{s_0}]$ for all $s\in U$.
\end{enumerate}
\end{proposition}

\begin{proof}
(a) The identity \eqref{eq:extdiff} holds for every exact sequence of coherent sheaves, and $\phi_z$ is nonnegative on cancellation vectors (Lemma~\ref{lem:localvalues}). So $\phi_z(\gamma(\cF'\oplus\cF''))\ge\phi_z(\gamma(\cF))$ for every $z$, and the inequalities follow by taking infima. The same holds along a filtration, as in the proof of Corollary~\ref{cor:monoext}(a).

(b) Let $\ell$ be the length of $\cE/\cF$. The table of $\cE/\cF$ is $\ell\sum_te_{0,t}$, on which every $\phi_z$ vanishes by \eqref{eq:divdiff}, so $\phi_z(\gamma(\cF))\ge\phi_z(\gamma(\cE))$ by \eqref{eq:extdiff}; and $\cE$, $\cF$ have the same rank.

(c) By Proposition~\ref{prop:cohseries}, for $m\ge m_0$ the infimum defining $\Psi_{\cF_{s_0}}(m)$ is attained at one of the finitely many slopes of the edges of $\Psi_{\cF_{s_0}}$ that meet $[m_0,\rk\cF_{s_0}]$. The argument of Corollary~\ref{cor:monoext}(b), applied to the columns $z_1,\dots,z_n$ of root sequences $z$ comparable with the chain of $\cF_{s_0}$ at these slopes, gives a neighborhood $U$ of $s_0$ on which $\phi_z(\gamma(\cF_s))\le\phi_z(\gamma(\cF_{s_0}))$ for these $z$; so $\Psi_{\cF_s}(m)\le m\lambda+\phi_z(\gamma(\cF_s))\le m\lambda+\Lambda_{\cF_{s_0}}(\lambda)=\Psi_{\cF_{s_0}}(m)$ for $m\ge m_0$ and the corresponding slope $\lambda$.
\end{proof}

On $\PP^2$, (b) applies to a torsion-free sheaf $\cF$ and $\cE=\cF^{\vee\vee}$, a vector bundle with the same Harder--Narasimhan polygon; so there Theorem~\ref{thm:main} for $\cF^{\vee\vee}$ bounds the subsheaves of $\cF$ at least as well as $\Psi_\cF$ does. In (c) we do not know whether $U$ can be chosen independently of $m_0$; the tables may change in infinitely many columns, as for the extensions $0\to\OO(-1)\to\cF_s\to\cI_p\to0$ on $\PP^2$, $s\in\Ext^1(\cI_p,\OO(-1))\cong K$, where $\cF_s\cong T_{\PP^2}(-2)$ and $h^1(\cF_s(t))=0$ for $s\ne0$ and $t\le-2$, while $h^1(\cF_0(t))=1$ for all $t<0$.

\section{The plane}\label{sec:plane}

In this section and the next, $n=2$, except in Example~\ref{ex:conicP3}. A root sequence is a pair $z=(z_1>z_2)$, with center $c(z)=\frac12(z_1+z_2)$, gap $g(z)=z_1-z_2\ge1$ and slope $\mu(z)=-c(z)-\frac32$, and
\[
\gamma^z_{j,t}=\begin{cases}(t-z_1)(t-z_2) & j=0,\ t>z_1,\\ (z_1-t)(t-z_2) & j=1,\ z_2<t<z_1,\\ (t-z_1)(t-z_2) & j=2,\ t<z_2,\\ 0 & \text{otherwise.}\end{cases}
\]
Its Hilbert polynomial is $(t-z_1)(t-z_2)=(t-c)^2-\frac{g^2}4$, and on the $h^1$ row $\gamma^z_{1,t}=\big(\frac{g^2}4-(t-c)^2\big)_+$ for all $t\in\ZZ$. Ranks of pieces are $2q_k$, and Boij--Söderberg slopes lie in $\frac12\ZZ$.

\subsection{Gap moment and centers}
For $n=2$ the normalized spread of Section~\ref{sec:gap} is $\rho(z)=g(z)^2$. So for a vector bundle $F$ of rank $N$ on $\PP^2$ with chain decomposition $\sum_kq_k\gamma^{z^k}$, with $p_k=2q_k/N$, $c_k=c(z^k)$ and $g_k=g(z^k)$,
\[
G(F)=\sum_kp_k\,g_k^2=\tfrac2N\Gamma(F),\qquad\text{where}\quad\Gamma(F)=\sum_kq_k\,g_k^2 .
\]
For a torsion-free sheaf $\cF$ of rank $N$ on $\PP^2$, with $\mu=\mu(\cF)$ and $\delta=\delta(\cF)=\Delta(\cF)/N^2$ as in Section~\ref{sec:gap}, Riemann--Roch reads
\begin{equation}\label{eq:RRP2}
\frac{\chi(\cF(t))}{N}=\frac12\Big(t+\mu+\frac32\Big)^2-\frac18-\frac{\delta}{2}.
\end{equation}

\begin{lemma}\label{lem:basic}
Let $F$ be a vector bundle of rank $N$ on $\PP^2$.
\begin{enumerate}[label=\textup{(\roman*)}]
\item $\chi(F(t))/N=\frac12\sum_kp_k(t-c_k)^2-G(F)/8$.
\item $\bar c=-\mu(F)-\frac32$, i.e.\ $\deg\gamma(F)=\deg F$, and $G(F)-1=4V_c(F)+4\delta(F)$.
\item $G(F)\ge1$, with equality iff $F$ splits.
\item $\sum_th^1(F(t))=\sum_kq_k\frac{g_k^3-g_k}6=\frac N{12}\sum_kp_k(g_k^3-g_k)$, and this is at least $\frac N6\big(G(F)-1\big)$ and at least $\frac N{12}\sqrt{G(F)}\,\big(G(F)-1\big)$.
\item $G(F)-1=4\big(V_c(F)-V_{HN}(F)\big)+4\sum_iw_i\,\delta(\cQ_i)$, with the notation of Lemma~\textup{\ref{lem:gapPn}(v)}.
\end{enumerate}
\end{lemma}

\begin{proof}
(i)--(iii) and (v) are Lemma~\ref{lem:gapPn}(i)--(iii) and (v) for $n=2$, since $\sum_kp_k(t-c_k)^2=t^2-2\bar c\,t+V_c(F)+\bar c^{\,2}$.
(iv) $\sum_{u=1}^{g-1}u(g-u)=\frac{g^3-g}6$, and $\frac{g^3-g}6\ge\frac{g^2-1}3$ is $(g^2-1)(g-2)\ge0$ for integers $g\ge1$; this first bound is Lemma~\ref{lem:gapPn}(iv) for $n=2$. For the second bound, $g^3-g=h(g^2)$ with $h(x)=x^{3/2}-x^{1/2}$, which is convex on $[0,\infty)$; Jensen's inequality gives $\sum_kp_kh(g_k^2)\ge h(G(F))$.
\end{proof}

For $n=2$, Corollary~\ref{cor:Q} reads $G(F)\ge1+4\sum_iw_i\,\delta(\cQ_i)$. By (ii) and Corollary~\ref{cor:monoext}(a), $G(F)\le G(\gr F)$ for every filtration with locally free subquotients, since $\delta$ depends only on the Hilbert polynomial. The constant term of the Hilbert polynomial, $\sum_kq_k(c_k^2-g_k^2/4)$, is linear in the table, so $\Gamma=4\sum_kq_kc_k^2$ minus a linear function; by Theorem~\ref{thm:concave}, $\Gamma$ is concave on $\CC$, and $\sum_jq'_jg(y_j)^2\le\Gamma(\gamma)$ for every decomposition $\gamma=\sum_jq'_j\gamma^{y_j}$.

\begin{lemma}[two twists]\label{lem:twotwist}
For every vector bundle $F$ of rank $N$ on $\PP^2$ and every $t\in\ZZ$,
\[
G(F)-1\ \ge\ \frac4N\Big(h^1(F(t))+h^1(F(t+1))\Big).
\]
\end{lemma}

\begin{proof}
Since $G-1=\frac2N\sum_kq_k(g_k^2-1)$, it suffices to show $g^2-1\ge2(\gamma^z_{1,t}+\gamma^z_{1,t+1})$ for every $z$. If $t,t+1\in(z_2,z_1)$, put $a=z_1-t$, $b=t-z_2$; then $g^2-1-2[ab+(a-1)(b+1)]=(a-b-1)^2$. If exactly one of $t,t+1$ lies in $(z_2,z_1)$, it is $z_2+1$ or $z_1-1$, the right-hand side is $2(g-1)$, and $g^2-1-2(g-1)=(g-1)^2$. If neither does, the right-hand side is $0\le g^2-1$.
\end{proof}

\subsection{The extreme slopes}
For a nonzero table $\gamma\in\CC$ on $\PP^2$ put
\[
M_j=\max\{t:\gamma_{j,t}\ne0\}\ \ (j=1,2),\qquad m_j=\min\{t:\gamma_{j,t}\ne0\}\ \ (j=0,1),
\]
with $M_1=-\infty$ and $m_1=+\infty$ if the $h^1$ row vanishes.

\begin{corollary}\label{cor:endsP2}
For $\gamma\in\CC$ on $\PP^2$, $\gamma\ne0$, the first and the last element of the chain are
\[
z^{\max}=\big(\max(M_1+1,\,M_2+2),\ M_2+1\big),\qquad z^{\min}=\big(m_0-1,\ \min(m_1-1,\,m_0-2)\big),
\]
and the largest and smallest Boij--Söderberg slopes of $\gamma$ are
\[
\mu^\BS_{\max}(\gamma)=\max\Big(-m_0,\ -\frac{m_0+m_1+1}2\Big),\qquad
\mu^\BS_{\min}(\gamma)=\min\Big(-M_2-3,\ -\frac{M_1+M_2+5}2\Big).
\]
\end{corollary}

\begin{proof}
This is Lemma~\ref{lem:ends} for $n=2$; the slopes are $-c-\frac32$.
\end{proof}

By Corollary~\ref{cor:supernatural}(b), $\mu_{\max}(F)\le\mu^\BS_{\max}(F)$. At the first step of the Harder--Narasimhan filtration one can say more: the chain of $F$ reaches at least as far as that of the subsheaf.

\begin{proposition}\label{prop:firstedge}
Let $F$ be a vector bundle on $\PP^2$ and $F_1\subset F$ a saturated subsheaf with $\mu_{\max}(F/F_1)<\mu(F_1)$, for instance the maximal destabilizing subsheaf. Then
\[
\mu^\BS_{\max}(F)\ \ge\ \mu^\BS_{\max}(F_1)\ \ge\ \mu(F_1).
\]
\end{proposition}

\begin{proof}
$F_1$ is reflexive, hence locally free. Put $F_2=F/F_1$, $\mu_1=\mu(F_1)$, $\bar c_1=-\mu_1-\frac32$, and use the invariants $m_0,m_1$ of Corollary~\ref{cor:endsP2}. First, $m_0(F)\le m_0(F_1)$, because $H^0(F_1(t))\subset H^0(F(t))$. Second, for $t<-\mu_{\max}(F_2)$ we have $H^0(F_2(t))=0$, so $H^1(F_1(t))\to H^1(F(t))$ is injective.

If $m_1(F_1)-1\le m_0(F_1)-2$, the last chain element of $F_1$ is $(m_0(F_1)-1,m_1(F_1)-1)$, so $m_1(F_1)\le c(z^{\min}(F_1))+\frac12$. The last center is the smallest, hence at most the mean center $\bar c_1$ (Lemma~\ref{lem:basic}(ii)); so $m_1(F_1)\le \bar c_1+\frac12=-\mu_1-1<-\mu_{\max}(F_2)$. Hence $m_1(F)\le m_1(F_1)$, and by Corollary~\ref{cor:endsP2}
\[
c(z^{\min}(F))\le\tfrac12\big(m_0(F)+m_1(F)-2\big)\le\tfrac12\big(m_0(F_1)+m_1(F_1)-2\big)=c(z^{\min}(F_1)).
\]
Otherwise $c(z^{\min}(F_1))=m_0(F_1)-\frac32\ge m_0(F)-\frac32\ge c(z^{\min}(F))$. In both cases $\mu^\BS_{\max}(F)\ge\mu^\BS_{\max}(F_1)$. The second inequality holds because the rank-weighted mean of the Boij--Söderberg slopes of $F_1$ is $\mu_1$.
\end{proof}

The first inequality compares the ends of the chains of $F$ and $F_1$; Theorem~\ref{thm:HNmono} below compares the whole polygons when $F_1$ is a step of the Harder--Narasimhan filtration.

\begin{remark}[no monotonicity for arbitrary subbundles]\label{rem:nomono}
The slope hypothesis cannot be dropped. Let $\cG=\ker(\OO^3\to\OO(2))$ for a regular sequence of three conics $q_1,q_2,q_3$ (for instance three general conics); then $H^1(\cG(t))$ is the degree $t+2$ part of $K[x_0,x_1,x_2]/(q_1,q_2,q_3)$, of dimensions $1,3,3,1$ for $t=-2,\dots,1$, and $H^0(\cG(1))=0$, since there are no linear syzygies. So $\cG$ is stable of rank $2$ and slope $-1$, with chain $\frac13\gamma^{(2,-2)}+\frac13\gamma^{(1,-2)}+\frac13\gamma^{(1,-3)}$. The subbundle $F_1=\OO\oplus \cG\subset\OO^4$ has $m_0(F_1)=0$ and $m_1(F_1)=-2$, so by Corollary~\ref{cor:endsP2} its chain ends with $(-1,-3)$; indeed it is
\[
\tfrac13\gamma^{(2,-2)}+\tfrac7{12}\gamma^{(1,-2)}+\tfrac16\gamma^{(1,-3)}+\tfrac1{12}\gamma^{(0,-3)}+\tfrac13\gamma^{(-1,-3)},
\]
with largest slope $\frac12$, while all Boij--Söderberg slopes of $\OO^4$ are $0$. So $\Psi_{F_1}(1)=\frac14>0=\Psi_{\OO^4}(1)$ and $\Lambda_{F_1}(0)=\frac13>0=\Lambda_{\OO^4}(0)$. Here $\mu(F_1)=-\frac23<2=\mu(\OO^4/F_1)$. Moreover $\gamma(\OO)=\frac12\gamma^{(-1,-2)}$, and $(-1,-2)$ is incomparable with $(1,-3)$; so the union of the chains of $\OO$ and $\cG$ is not a chain, and the re-chaining creates the root sequences $(0,-3)$ and $(-1,-3)$, which occur in neither.
\end{remark}

\subsection{The three-term formula}
For integers $a>b$ and a table $\gamma$ on $\PP^2$ put
\[
\phi_{a,b}(\gamma)=\phi_{(a,b)}(\gamma)=\frac{\gamma_{0,a}-\gamma_{0,b}+\gamma_{1,b}}{a-b}.
\]
With the Hilbert polynomial $\chi_\gamma(t)=\sum_j(-1)^j\gamma_{j,t}$, also $\phi_{a,b}(\gamma)=\big(\gamma_{0,a}+\gamma_{2,b}-\chi_\gamma(b)\big)/(a-b)$.

\begin{theorem}[three-term formula]\label{thm:three}
Let $\gamma\in\CC$ be a table on $\PP^2$ with chain decomposition $\sum_kq_k\gamma^{z^k}$, let $\lambda\in\frac12\ZZ$, and put $c=c_\lambda$.
\begin{enumerate}[label=\textup{(\alph*)}]
\item If $a>b$ are integers with $a+b=2c$ and $y$ is a root sequence comparable with $(a,b)$, then $\phi_{a,b}(\gamma^y)=2\big(c-c(y)\big)_+$; for arbitrary $y$, $\phi_{a,b}(\gamma^y)\ge2\big(c-c(y)\big)_+$.
\item There are integers $a>b$ with $a+b=2c$ such that $(a,b)$ is comparable with every $z^k$, and for each of them $\Lambda_\gamma(\lambda)=\phi_{a,b}(\gamma)$.
\item $\Lambda_\gamma(\lambda)=\min\{\phi_{a,b}(\gamma):\ a>b,\ a+b=2c\}$.
\end{enumerate}
Consequently, for a vector bundle $F$ on $\PP^2$ and $0\le m\le\rk F$,
\[
\Psi_F(m)=\min_{a>b}\Big(\frac{h^0(F(a))-h^0(F(b))+h^1(F(b))}{a-b}-m\Big(\frac{a+b}2+\frac32\Big)\Big).
\]
\end{theorem}

\begin{proof}
(a) For comparable $y$ this is Lemma~\ref{lem:comparable}. Otherwise put $P(u)=(u-y_1)(u-y_2)$, so that $P(a)-P(b)=2(a-b)(c-c(y))$, and recall $\phi_{a,b}(\gamma)=(\gamma_{0,a}-\gamma_{0,b}+\gamma_{1,b})/(a-b)$. If $y_1>a$ and $y_2<b$, both columns $a,b$ lie in row $1$, where the entries are $-P$, so $\phi_{a,b}(\gamma^y)=\gamma^y_{1,b}/(a-b)$, while $\gamma^y_{1,a}-\gamma^y_{1,b}=2(a-b)(c(y)-c)$. If $y_1<a$ and $y_2>b$, then $\phi_{a,b}(\gamma^y)=\gamma^y_{0,a}/(a-b)$, while $\gamma^y_{0,a}-\gamma^y_{2,b}=2(a-b)(c-c(y))$. In both cases
\begin{equation}\label{eq:incomparable}
\phi_{a,b}(\gamma^y)=2\big(c-c(y)\big)_++\frac{\min\big(\gamma^y_{\bullet,a},\gamma^y_{\bullet,b}\big)}{a-b}\ \ge\ 2\big(c-c(y)\big)_+,
\end{equation}
where $\gamma^y_{\bullet,u}=|P(u)|$ is the entry of $\gamma^y$ in the column $u$. (b) is Lemma~\ref{lem:adapted} for $n=2$. (c) By (a), $\phi_{a,b}(\gamma)=\sum_kq_k\phi_{a,b}(\gamma^{z^k})\ge\Lambda_\gamma(\lambda)$ for every pair, with equality in (b). The last formula follows with Lemma~\ref{lem:legendre}(i). This proof uses only Lemmas~\ref{lem:comparable}--\ref{lem:adapted}, not Theorem~\ref{thm:main}.
\end{proof}

We call a pair $(a,b)$ as in (b) \emph{adapted} to the chain of $\gamma$ (or of $F$). The values of $\phi_{a,b}$ on the cancellation vectors are $\phi_{a,b}(e_{0,a}+e_{1,a})=\phi_{a,b}(e_{1,b}+e_{2,b})=\frac1{a-b}$, by Lemma~\ref{lem:localvalues}.

On incomparable root sequences, $\phi_{a,b}$ measures how much the chain gains over other decompositions. For a root sequence $y$ and $u\in\ZZ$ let $\gamma^y_{\bullet,u}=|(u-y_1)(u-y_2)|$ be the entry of $\gamma^y$ in the column $u$, and for a root sequence $z=(a,b)$ put
\[
E_y(z)=\begin{cases}0 & \text{if $z$ and $y$ are comparable},\\[1ex] \dfrac{\min\big(\gamma^y_{\bullet,a},\gamma^y_{\bullet,b}\big)}{a-b} & \text{otherwise}.\end{cases}
\]

\begin{corollary}[the gain formula]\label{cor:gain}
Let $\gamma=\sum_yq_y\gamma^y$, $q_y\ge0$, be any decomposition of a table $\gamma\in\CC$ on $\PP^2$, and let $c=c_\lambda\in\frac12\ZZ$. Then
\[
\Lambda_\gamma(\lambda)=\sum_y2q_y\big(c-c(y)\big)_++\min_{z:\ c(z)=c}\ \sum_yq_y\,E_y(z),
\]
the minimum over all root sequences $z$ with center $c$. It is attained at the root sequences with center $c$ that are comparable with the chain of $\gamma$.
\end{corollary}

\begin{proof}
By Theorem~\ref{thm:three}, $\Lambda_\gamma(\lambda)=\min_z\sum_yq_y\phi_z(\gamma^y)$, with the minimum attained as stated, and $\phi_z(\gamma^y)=2(c-c(y))_++E_y(z)$ for $z$ with $c(z)=c$, by Theorem~\ref{thm:three}(a) for comparable $y$ and $z$ and by \eqref{eq:incomparable} otherwise.
\end{proof}

So the chain improves on any decomposition by the \emph{gain} $\min_z\sum_yq_yE_y(z)$, which is positive exactly when no root sequence with center $c$ is comparable with all the $y$ with $q_y>0$. This quantifies Theorem~\ref{thm:concave}(a) on $\PP^2$. For the table $\gamma^{(3,0)}+\gamma^{(2,1)}$ of Example~\ref{ex:sum} and $c=\frac32$, the first sum vanishes and the gain is $\frac23$, attained at $z=(3,0)$; indeed $\Lambda_\gamma(-3)=\frac23$ for the chain decomposition found there.

\subsection{The four defects}
On $\PP^2$ the defect of Proposition~\ref{prop:marginPn} is a sum of nonnegative terms.

\begin{proposition}[four defects]\label{prop:fourdefect}
Let $F$ be a vector bundle on $\PP^2$, $F'\subset F$ a nonzero subsheaf of rank $N'$ and slope $\mu'$ with torsion-free quotient $F''=F/F'$, $\lambda\in\frac12\ZZ$, and $(a,b)$ a pair adapted to the chain of $F$ with $a+b=-2\lambda-3$. Write $h^j_t(\cG)=h^j(\cG(t))$, and let $\partial^0_a\colon H^0(F''(a))\to H^1(F'(a))$ and $\partial^1_b\colon H^1(F''(b))\to H^2(F'(b))$ be the connecting homomorphisms. Then
\begin{multline*}
(a-b)\big(\Lambda_F(\lambda)-N'(\mu'-\lambda)\big)\\=\big(h^2_b(F')-h^2_a(F')\big)+\big(h^0_a(F'')-h^0_b(F'')\big)+\big(h^1_a(F')-\rk\partial^0_a\big)+\big(h^1_b(F'')-\rk\partial^1_b\big),
\end{multline*}
and all four terms are nonnegative.
\end{proposition}

\begin{proof}
By Theorem~\ref{thm:three}(b), $(a-b)\Lambda_F(\lambda)=h^0_a(F)-h^0_b(F)+h^1_b(F)$. The long exact cohomology sequences give $h^0_a(F)=h^0_a(F')+h^0_a(F'')-\rk\partial^0_a$ and $h^1_b(F)-h^0_b(F)=h^1_b(F')-h^0_b(F')+h^1_b(F'')-h^0_b(F'')-\rk\partial^1_b$. Writing $h^0_a(F')=\chi(F'(a))+h^1_a(F')-h^2_a(F')$ and $h^1_b(F')-h^0_b(F')=h^2_b(F')-\chi(F'(b))$, we obtain the identity, since $\chi(F'(a))-\chi(F'(b))=\frac{N'}2(a-b)(a+b+2\mu'+3)=(a-b)\,N'(\mu'-\lambda)$ by \eqref{eq:RRP2}. The last two terms are nonnegative trivially. A nonzero linear form is a nonzerodivisor on the torsion-free sheaves $F'$ and $F''$; multiplication by it gives injections $H^0(F''(t))\to H^0(F''(t+1))$, and surjections $H^2(F'(t))\to H^2(F'(t+1))$, since the cokernel of $F'(t)\to F'(t+1)$ is supported on a line. So $h^0_t(F'')$ is nondecreasing and $h^2_t(F')$ is nonincreasing in $t$, and $a>b$.
\end{proof}

For $n=2$, \eqref{eq:marginPn} is this identity before the rearrangement: $B_1=h^1_a(F')-h^2_a(F')+h^0_a(F'')-\rk\partial^0_a$ and $B_2=h^2_b(F')+h^1_b(F'')-h^0_b(F'')-\rk\partial^1_b$; so $B_1$ can be negative, e.g.\ $B_1=-1$, $B_2=3$ for $0\to\OO\to\OO^2\to\OO\to0$ and $(a,b)=(-3,-4)$, $\lambda=2$. Proposition~\ref{prop:fourdefect}, Theorem~\ref{thm:three}(b) and Lemma~\ref{lem:legendre} give an elementary proof of Theorem~\ref{thm:main} on $\PP^2$; it corresponds to the test sheaf $\OO_{(a-b)L}(a)$, with resolution $\OO(b)\to\OO(a)$. They also describe the equality cases: the supporting-line equality $\Lambda_F(\lambda)=\deg F'-\lambda\rk F'$ at a step $F'$ of the Harder--Narasimhan filtration holds iff all four terms vanish for a pair $(a,b)$ adapted at $\lambda$.

The inequality $P_F\le\Psi_F$ cannot be obtained from the generic splitting type on a line. That type dominates $P_F$, but it is in general more spread out than $\Psi_F$: for $T_{\PP^2}(-1)$ it is $(1,0)$, while the Boij--Söderberg slopes are $(\frac12,\frac12)$.

The same computation applies to coherent sheaves on $\PP^2$, with $\Lambda$ and $\Psi$ defined as in Section~\ref{sec:coherent}.

\begin{proposition}[Theorem~D for coherent sheaves]\label{prop:cohD}
Let $\cF$ be a coherent sheaf on $\PP^2$, $\lambda\in\frac12\ZZ$, and $a>b$ integers with $a+b=-2\lambda-3$. Then $\Lambda_\cF(\lambda)\le\phi_{a,b}(\gamma(\cF))$, with equality if $(a,b)$ is comparable with every element of the chain of $\cF$. For every coherent subsheaf $\cF'\subset\cF$, with $\cF''=\cF/\cF'$ and the notation of Proposition~\ref{prop:fourdefect},
\begin{multline*}
(a-b)\big(\phi_{a,b}(\gamma(\cF))-\deg\cF'+\lambda\rk\cF'\big)\\=\big(h^2_b(\cF')-h^2_a(\cF')\big)+\big(h^0_a(\cF'')-h^0_b(\cF'')\big)+\big(h^1_a(\cF')-\rk\partial^0_a\big)+\big(h^1_b(\cF'')-\rk\partial^1_b\big),
\end{multline*}
and all four terms are nonnegative.
\end{proposition}

\begin{proof}
The inequality is the definition of $\Lambda_\cF$, and the equality is Proposition~\ref{prop:cohseries}(a). Apart from its first sentence, the proof of Proposition~\ref{prop:fourdefect} uses only the definition of $\phi_{a,b}$, and it applies to $\cF'\subset\cF$, with $\deg\cF'-\lambda\rk\cF'$ in place of $N'(\mu'-\lambda)$: Riemann--Roch is used in the form of Section~\ref{sec:tables}; $H^2(\cF'(t))\to H^2(\cF'(t+1))$ is surjective for every nonzero linear form, since $H^2$ is right exact on $\PP^2$ and the cokernel of $\cF'(t)\to\cF'(t+1)$ is supported on a line; and a general linear form is a nonzerodivisor on $\cF''$.
\end{proof}

In particular Proposition~\ref{prop:cohD} gives Proposition~\ref{prop:cohA}(a) on $\PP^2$ without \cite{MiSp08}.

\subsection{Tables of extensions}
Proposition~\ref{prop:fourdefect} works at the level of tables. This gives necessary conditions on the cohomology tables of extensions of two given bundles, one of the problems of Erman and Sam \cite[Problem~4.3]{ErS}: if $\gamma_1=\gamma(F')$ and $\gamma_2=\gamma(F'')$ are tables of vector bundles, every extension of $F''$ by $F'$ has a table of the following form, with $d^j_t=\rk\partial^j_t$.

\begin{corollary}[tables of extensions]\label{cor:conicQ}
Let $\gamma_1,\gamma_2\in\CC$ be tables on $\PP^2$, $\gamma_1\ne0$, and $D=\sum_td^0_t(e_{0,t}+e_{1,t})+\sum_td^1_t(e_{1,t}+e_{2,t})$ with $0\le d^0_t\le\gamma_{1;1,t}$ and $0\le d^1_t\le\gamma_{2;1,t}$, where $\gamma_{i;j,t}$ are the entries of $\gamma_i$. If $\gamma=\gamma_1+\gamma_2-D\in\CC$, then $\Psi_\gamma(\rk\gamma_1)\ge\deg\gamma_1$.
\end{corollary}

\begin{proof}
Since $D$ is $\chi$-neutral, $\rk\gamma=\rk\gamma_1+\rk\gamma_2$ and $\deg\gamma=\deg\gamma_1+\deg\gamma_2$. As in Lemma~\ref{lem:legendre}(i),(iii), it suffices to show $\Lambda_\gamma(\lambda)\ge\rk\gamma_1\,(\mu(\gamma_1)-\lambda)$ for $\lambda\in\frac12\ZZ$, where $\mu(\gamma_1)=\deg\gamma_1/\rk\gamma_1$. Let $c=c_\lambda$ and choose $(a,b)$ adapted to $\gamma$. Then
\[
\gamma_{0,a}=\gamma_{1;0,a}+\gamma_{2;0,a}-d^0_a\ \ge\ \gamma_{1;0,a}+\gamma_{2;0,a}-\gamma_{1;1,a},
\]
\[
\gamma_{1,b}-\gamma_{0,b}=(\gamma_{1;1,b}-\gamma_{1;0,b})+(\gamma_{2;1,b}-\gamma_{2;0,b})-d^1_b\ \ge\ \gamma_{1;1,b}-\gamma_{1;0,b}-\gamma_{2;0,b}.
\]
Every table in $\CC$ has a nondecreasing row $0$ and a nonincreasing row $2$, because each $\gamma^z$ has. As in the proof of Proposition~\ref{prop:fourdefect},
\[
(a-b)\Lambda_\gamma(\lambda)\ \ge\ \chi_{\gamma_1}(a)-\chi_{\gamma_1}(b)+(\gamma_{1;2,b}-\gamma_{1;2,a})+(\gamma_{2;0,a}-\gamma_{2;0,b})\ \ge\ \chi_{\gamma_1}(a)-\chi_{\gamma_1}(b),
\]
and by the formula for the Hilbert polynomial in Section~\ref{sec:tables}, $\chi_{\gamma_1}(a)-\chi_{\gamma_1}(b)=(a-b)\rk\gamma_1\,(\mu(\gamma_1)-\lambda)$.
\end{proof}

The tables of actual extensions also satisfy $d^0_t\le\gamma_{2;0,t}$ and $d^1_t\le\gamma_{1;2,t}$. Even with these bounds, the conditions of Corollary~\ref{cor:conicQ} are necessary but not sufficient. Let $F'=\Omega^1(1)\oplus\OO(1)$ and $F''=\OO(2)$. Then $\Ext^1(F'',F')=H^1(\Omega^1(-1))\oplus H^1(\OO(-1))=0$, so the only extension is split. The Euler sequence gives $\gamma(F')=3\gamma(\OO)+e_{0,-1}+e_{1,-1}$, and the table $\gamma(F')+\gamma(F'')-(e_{0,-1}+e_{1,-1})=3\gamma(\OO)+\gamma(\OO(2))$ satisfies all the bounds, with $d^0_{-1}=1\le\min\big(h^1(F'(-1)),h^0(F''(-1))\big)$; but it is not the table of an extension of $F''$ by $F'$.

\begin{remark}\label{rem:linear}
The chain enters through the choice of an adapted pair; a linear lower bound valid on all of $\CC$ cannot replace it. A linear functional $L$ with $L(\gamma^z)\le2(\mu(z)-\lambda)_+$ for all $z$ gives a lower bound $L(\gamma(F))\le\Lambda_F(\lambda)$. But then $L(\gamma(F))\le\sum_j2q'_j(\mu(y_j)-\lambda)_+$ for \emph{every} decomposition $\gamma(F)=\sum_jq'_j\gamma^{y_j}$; by Theorem~\ref{thm:concave}(a) the chain gives the largest of these values, and the smallest can be below the vertex target. Let $F_1=\coker(\OO^2\to\OO(1)^2\oplus\OO(2)^2)$ and $F_2=\coker(\OO(-1)\to\OO\oplus\OO(1)^2)$ be general, of slopes $3$ and $\frac32$, and let $F=F_1\oplus F_2$. By Serre duality $H^1(F_1(t))$ and $H^1(F_2(t))$ are dual to the graded pieces of degree $-t-3$ of the cokernels of the transposed maps $S(-1)^2\oplus S(-2)^2\to S^2$ and $S\oplus S(-1)^2\to S(1)$, $S=K[x_0,x_1,x_2]$. For general maps these cokernels have finite length, with Hilbert functions $2,4,4,2$ in degrees $0,\dots,3$ (from the Buchsbaum--Rim complex, a resolution here) and $1,2,1$ in degrees $-1,0,1$ (a complete intersection of type $(1,2,2)$); the $h^0$- and $h^2$-rows then follow from the presentations, which are injective on global sections, and Riemann--Roch. Then $\gamma(F)=\gamma(F_1)+\gamma(F_2)$, with chain
\[
\tfrac12\gamma^{(-1,-4)}+\tfrac12\gamma^{(-2,-5)}+\tfrac25\gamma^{(-2,-6)}+\tfrac15\gamma^{(-2,-7)}+\tfrac25\gamma^{(-3,-7)},
\]
and $\Lambda_F(2)=2=\rk F_1\cdot(\mu_1-2)$. But also
\[
\gamma(F)=\gamma^{(-3,-4)}+\tfrac12\gamma^{(-2,-7)}+\tfrac14\gamma^{(-2,-5)}+\tfrac14\gamma^{(-1,-6)},
\]
which gives the value $1$ at $\lambda=2$. The functionals $\phi_{a,b}$ work the other way: each is an upper bound for $\Lambda$ on all of $\CC$, and $\Lambda$ is their minimum.
\end{remark}

Theorem~\ref{thm:main} says that $\Lambda_F\ge\Lambda_{P_F}$. On $\PP^2$ the table of every step of the Harder--Narasimhan filtration gives a stronger lower bound.

\begin{lemma}\label{lem:lowpair}
Let $\gamma\in\CC$ on $\PP^2$, $\gamma\ne0$, and let $B=M_2+1$ be the second root of the first element of its chain (Corollary~\textup{\ref{cor:endsP2}}). For every $c\in\frac12\ZZ$ there are integers $a>b$ with $a+b=2c$ and $b\le B$ such that $(a,b)$ is comparable with every element of the chain.
\end{lemma}

\begin{proof}
Let $z^{\max}$ and $z^{\min}$ be the first and the last element of the chain, so $z^{\max}_2=B$. If $c\ge c(z^{\max})$, take $(a,b)=(2c-B,B)$; then $(a,b)\ge z^{\max}$ and $a\ge z^{\max}_1>B$. If $c\le c(z^{\min})$, take $(a,b)=(z^{\min}_1,2c-z^{\min}_1)$; then $(a,b)\le z^{\min}$, and $b\le z^{\min}_2<a$, $b\le z^{\max}_2=B$. Otherwise $c(z^{k+1})\le c\le c(z^k)$ for some $k$; take $a=\min(z^k_1,\,2c-z^{k+1}_2)$ and $b=2c-a$. Then $z^{k+1}\le(a,b)\le z^k$ and $a>b$, so $(a,b)$ is comparable with the whole chain, and $b\le z^k_2\le B$.
\end{proof}

\begin{proposition}\label{prop:Mcrit}
Let $0\to \cG\to F\to Q\to0$ be an exact sequence of coherent sheaves on $\PP^2$, with $F$ a vector bundle, $\cG\ne0$ and $Q\ne0$ torsion-free, and let $M_2(F)=\max\{t:H^2(F(t))\ne0\}$. If $H^0\big(Q(M_2(F)+1)\big)=0$, then $\cG$ is locally free and $\Lambda_\cG\le\Lambda_F$; equivalently, $\Psi_\cG\le\Psi_F$ on $[0,\rk \cG]$. The hypothesis holds in each of the following cases:
\begin{enumerate}[label=\textup{(\alph*)}]
\item $\mu_{\max}(Q)<\mu_{\min}(F)+2$;
\item $Q$ is $\mu$-semistable and $\mu(Q)<\mu_{\min}(\cG)+2$.
\end{enumerate}
\end{proposition}

\begin{proof}
$\cG$ is saturated in $F$, hence reflexive, hence locally free. Let $\lambda\in\frac12\ZZ$, $c=-\lambda-\frac32$, and let $(a,b)$ be a pair for $\gamma(F)$ as in Lemma~\ref{lem:lowpair}. Multiplication by a linear form shows that $h^0(Q(t))$ is nondecreasing in $t$, since $Q$ is torsion-free; so $H^0(Q(b))=0$. The long exact cohomology sequences give $h^0_a(F)=h^0_a(\cG)+h^0_a(Q)-\rk\partial^0_a$, $h^0_b(F)=h^0_b(\cG)$ and $h^1_b(F)=h^1_b(\cG)+h^1_b(Q)-\rk\partial^1_b$, so
\[
(a-b)\big(\phi_{a,b}(\gamma(F))-\phi_{a,b}(\gamma(\cG))\big)=\big(h^0_a(Q)-\rk\partial^0_a\big)+\big(h^1_b(Q)-\rk\partial^1_b\big)\ \ge\ 0 .
\]
By Theorem~\ref{thm:three}, $\Lambda_F(\lambda)=\phi_{a,b}(\gamma(F))\ge\phi_{a,b}(\gamma(\cG))\ge\Lambda_\cG(\lambda)$. Both sides are linear between consecutive points of $\frac12\ZZ$, so $\Lambda_\cG\le\Lambda_F$ everywhere. Since $\Lambda_\gamma(\lambda)=\max_{0\le m\le\rk\gamma}(\Psi_\gamma(m)-m\lambda)$ and $\Psi_\gamma(m)=\min_\lambda(m\lambda+\Lambda_\gamma(\lambda))$ (Lemma~\ref{lem:legendre}(i)), this is equivalent to $\Psi_\cG\le\Psi_F$ on $[0,\rk \cG]$.

If $H^2(F(t))\ne0$, Serre duality gives a nonzero map $F\to\OO(-t-3)$, whose image is a torsion-free quotient of rank one and slope at most $-t-3$; so $M_2(F)+1\le-\mu_{\min}(F)-2$. In case (a), $\mu_{\max}(Q(M_2(F)+1))<0$, so $Q(M_2(F)+1)$ has no sections. Case (b) implies (a): the last Harder--Narasimhan factor $R$ of $F$, of slope $\mu_{\min}(F)$, is semistable, so $\mu(R)\ge\mu_{\min}(\cG)$ if $\cG\to R$ is nonzero, and otherwise $R$ is a quotient of $Q$; hence $\mu_{\min}(F)\ge\min(\mu_{\min}(\cG),\mu(Q))$, and $\mu(Q)<\mu_{\min}(\cG)+2$ gives $\mu_{\max}(Q)=\mu(Q)<\mu_{\min}(F)+2$.
\end{proof}

\begin{theorem}[monotonicity along the Harder--Narasimhan filtration]\label{thm:HNmono}
Let $0=F_0\subset F_1\subset\dots\subset F_s=F$ be the Harder--Narasimhan filtration of a vector bundle $F$ on $\PP^2$. Then
\[
\Lambda_{F_1}\le\Lambda_{F_2}\le\dots\le\Lambda_{F_s}=\Lambda_F,\qquad\text{i.e.}\qquad \Psi_{F_i}\le\Psi_{F_{i'}}\ \text{on }[0,\rk F_i]\ \ (i\le i').
\]
\end{theorem}

\begin{proof}
Each $F_i$ is saturated in $F$, hence locally free. For $i\ge2$ the quotient $\cQ_i=F_i/F_{i-1}$ is $\mu$-semistable with $\mu(\cQ_i)<\mu(\cQ_{i-1})=\mu_{\min}(F_{i-1})$; apply Proposition~\ref{prop:Mcrit}(b) to $0\to F_{i-1}\to F_i\to\cQ_i\to0$.
\end{proof}

\begin{remark}\label{rem:HNmono}
\begin{enumerate}[label=(\arabic*)]
\item For $\lambda\ge\mu_{i+1}$ (with $\mu_{s+1}=-\infty$), $\Lambda_{P_F}=\Lambda_{P_{F_i}}\le\Lambda_{F_i}\le\Lambda_F$ by Theorem~\ref{thm:main} for $F_i$, so in this range Theorem~\ref{thm:HNmono} strengthens Theorem~\ref{thm:main}: it compares the whole table of each step with that of $F$, not only its rank and degree. Its proof uses Theorem~\ref{thm:three}, Corollary~\ref{cor:endsP2} and Serre duality, but neither Theorem~\ref{thm:main} nor Bogomolov's inequality, and it holds over any algebraically closed field. It uses the two extreme rows of the plane, and does not extend as it stands to $\PP^n$, $n\ge3$.
\item The constant $2$ in Proposition~\ref{prop:Mcrit}(a) is optimal: in Remark~\ref{rem:nomono}, $\cG=\OO\oplus\ker(\OO^3\to\OO(2))\subset F=\OO^4$ has $Q=\OO(2)$, so $\mu_{\max}(Q)=\mu_{\min}(F)+2$, while $\Lambda_\cG(0)=\frac13>0=\Lambda_F(0)$. In (b), $\mu_{\min}(\cG)$ cannot be replaced by $\mu(\cG)$: adding $\OO(11)$ to $\cG$ and to $F$ gives $\mu(\cG\oplus\OO(11))+2=\frac{17}4>2=\mu(Q)$, and the chains, which are the unions of those of $\OO(11)$ and of $\cG$, resp.\ $F$, still give $\Lambda_{\cG\oplus\OO(11)}(0)>\Lambda_{F\oplus\OO(11)}(0)$. Proposition~\ref{prop:firstedge}, under the different hypothesis $\mu_{\max}(F/F_1)<\mu(F_1)$, only compares the ends of the chains.
\end{enumerate}
\end{remark}

On $\PP^3$ the analogue of Corollary~\ref{cor:conicQ} is false, even with both bounds $d^j_t\le\min(\gamma_{1;j+1,t},\gamma_{2;j,t})$ that hold for actual extensions. So no argument that uses only the tables of $F'$ and $F''$ and these numerical bounds for the ranks of the connecting homomorphisms can prove Theorem~\ref{thm:main} on $\PP^3$; some geometric input beyond the tables is needed, such as the test sheaves.

\begin{example}[tables do not suffice on $\PP^3$]\label{ex:conicP3}
On $\PP^3$ put
\begin{align*}
\gamma_1&=\gamma^{(5,4,-1)}+4\gamma^{(5,4,-2)}+15\gamma^{(1,0,-3)}+6\gamma^{(1,0,-4)},\\
\gamma_2&=4\gamma^{(4,0,-1)}+11\gamma^{(3,0,-1)}+50\gamma^{(2,0,-1)}+9\gamma^{(1,-4,-5)},\\
\gamma&=100\,\gamma^{(1,0,-1)}=600\,\gamma(\OO(-2)).
\end{align*}
A direct computation from the definition of $\gamma^z$ gives
\[
\gamma_1+\gamma_2-\gamma=\sum_{t=-3}^{3}d_t\,\big(e_{1,t}+e_{2,t}\big),\qquad (d_{-3},\dots,d_3)=(72,162,216,180,168,114,48),
\]
with $d_t=\gamma_{1;2,t}=\gamma_{2;1,t}$; moreover $\gamma_1$ has no $h^1$ row and $\gamma_2$ no $h^2$ row. So $\gamma$ is what the long exact sequences give for an extension of a bundle with table $\gamma_2$ by a bundle with table $\gamma_1$ in which $\partial^0=\partial^2=0$ and every $\partial^1_t\colon H^1(F''(t))\to H^2(F'(t))$ is an isomorphism; and up to a common positive integer multiple, $\gamma_1$ and $\gamma_2$ are tables of direct sums of supernatural bundles. All hypotheses of the analogue of Corollary~\ref{cor:conicQ} hold. But $\rk\gamma_1=156$ and $\deg\gamma_1=-288$, while $\Psi_\gamma$ is linear of slope $-2$, so $\Psi_\gamma(156)=-312<\deg\gamma_1$.

Similarly, with
\begin{align*}
\gamma'_1&=4\gamma^{(6,5,-1)}+20\gamma^{(2,0,-2)}+5\gamma^{(2,0,-3)}+2\gamma^{(2,0,-4)}+\gamma^{(2,0,-5)},\\
\gamma'_2&=2\gamma^{(5,0,-1)}+4\gamma^{(4,0,-1)}+10\gamma^{(3,0,-1)}+2\gamma^{(2,-5,-6)},\\
\gamma'&=50\,\gamma^{(2,0,-1)},
\end{align*}
one finds
\[
\gamma'_1+\gamma'_2-\gamma'=\sum_{t=-4}^{4}d'_t\,\big(e_{1,t}+e_{2,t}\big),\qquad (d'_{-4},\dots,d'_4)=(24,60,96,120,120,160,144,96,40),
\]
with $d'_t=\min(\gamma'_{1;2,t},\gamma'_{2;1,t})$, while $\rk\gamma'_1=192$ and $\deg\gamma'_1=-440>-448=\Psi_{\gamma'}(192)$.
\end{example}

The stated properties follow directly from the definition of $\gamma^z$. By Theorem~\ref{thm:main} neither pair comes from an actual subsheaf, and the test sheaves of Lemma~\ref{lem:testsheaf} show this directly. A vector bundle $F$ with table $\gamma$ has no intermediate cohomology, so it is $\OO(-2)^{\oplus600}$ by Horrocks' criterion \cite{Ho64}. For $z=(1,0,-1)$ the test sheaf is $\Tt=\OO_L(1)$, with resolution $\OO(-1)\to\OO^2\to\OO(1)$, and $h^0(F\otimes \Tt)=0$; but a subbundle $F'\subset F$ with table $\gamma_1$ and locally free quotient, as in an extension, would give, for $L$ general, $0<24=\deg F'+2\rk F'=\chi(F'\otimes\Tt)\le h^0(F'\otimes\Tt)\le h^0(F\otimes\Tt)$. For $\gamma'$, which is a multiple of the supernatural table $\gamma^{(2,0,-1)}$, every bundle with table $\gamma'$ is semistable by Corollary~\ref{cor:supernatural}; the test sheaf is $\Tt=(\OO/\cI_L^2)(2)$, with resolution $\OO(-1)^2\to\OO^3\to\OO(2)$, and $\phi_{(2,0,-1)}(\gamma')=0$ gives $h^0(F\otimes \Tt)=0$, while such a subbundle with table $\gamma'_1$ would give $\chi(F'\otimes\Tt)=3(\deg F'+\frac73\rk F')=24>0$. The cancellation runs through $\partial^1$, between two intermediate rows. On $\PP^2$ every connecting homomorphism has its source in row $0$ or its target in row $2$, and these are the rows whose monotonicity Proposition~\ref{prop:fourdefect} uses; on $\PP^3$ the injection $F'\otimes \Tt\to F\otimes \Tt$ carries information that the tables do not record.

\subsection{The equality cases of Theorem~A on \texorpdfstring{$\PP^2$}{P2}}\label{sec:equality}
Let $0=F_0\subset\dots\subset F_s=F$ be the Harder--Narasimhan filtration of a vector bundle $F$ on $\PP^2$, with factors $\cQ_i=F_i/F_{i-1}$ of slopes $\mu_1>\dots>\mu_s$. We write $N_i=\rk F_i$, and we call a vertex $(N_i,\deg F_i)$ of $P_F$ \emph{tight} if $\Psi_F(N_i)=\deg F_i$.

\begin{lemma}[contact]\label{lem:contact}
If $\Psi_F(m)=P_F(m)$ for some $m$ in the interior of an edge of $P_F$, then $\Psi_F=P_F$ on that edge. So the contact set $\{m:\Psi_F(m)=P_F(m)\}$ is the union of the tight vertices and of the edges of $P_F$ whose endpoints are both tight and on which $\Psi_F$ is linear.
\end{lemma}

\begin{proof}
On the edge, $\Psi_F-P_F$ is concave, nonnegative by Theorem~\ref{thm:main}, and vanishes at an interior point; a nonnegative concave function on an interval that vanishes at an interior point vanishes identically. For the last assertion, if $\Psi_F=P_F$ on an edge then $\Psi_F$ is linear there, and conversely.
\end{proof}

\begin{lemma}[Chern characters at tight vertices]\label{lem:Ktight}
Let $(N_i,\deg F_i)$ be a tight vertex, $0<i<s$, and let $\lambda\in\frac12\ZZ$ be a slope of $\Psi_F$ at $N_i$, i.e.\ $\Psi_F(N_i)=\min_{\lambda'}\big(N_i\lambda'+\Lambda_F(\lambda')\big)$ is attained at $\lambda$. Let $(a,b)$ be a pair adapted to the chain of $F$ with $a+b=-2\lambda-3$, as in Theorem~\textup{\ref{thm:three}(b)}, and let
\[
x=\sum_{\mu(z^k)>\lambda}q_k\gamma^{z^k}+f\,\gamma^{(a,b)},\qquad f\ge0 \text{ such that } \rk x=N_i .
\]
Then $x$ is well defined, $\chi_x(t)=\chi(F_i(t))$ for all $t$, i.e.\ $x$ and $F_i$ have the same Chern character, and
\[
h^0\big((F/F_i)(a)\big)=0,\qquad h^2\big(F_i(b)\big)=0 .
\]
\end{lemma}

\begin{proof}
Put $c=-\lambda-\frac32$, so $c(a,b)=c$. Since $\lambda$ minimizes $N_i\lambda'+\Lambda_F(\lambda')$, whose one-sided derivatives at $\lambda$ are $N_i-\sum_{\mu(z^k)\ge\lambda}2q_k\le0\le N_i-\sum_{\mu(z^k)>\lambda}2q_k$, the pieces of slope $>\lambda$ have total rank $\le N_i$, and those of slope $\ge\lambda$ total rank $\ge N_i$. The pieces of slope $>\lambda$ are those with center $<c$; being comparable with $(a,b)$, they are the $z^k<(a,b)$. Let $x'$ be their sum. If $\rk x'<N_i$, the chain contains a piece of slope exactly $\lambda$, of rank at least $N_i-\rk x'$; it is comparable with $(a,b)$ and has the same center, so it is $(a,b)$. So $x$ is well defined, $f$ is at most the coefficient of $\gamma^{(a,b)}$ in the chain, $\rk x=N_i$, and $\deg x=\deg x'+\lambda(N_i-\rk x')=N_i\lambda+\Lambda_F(\lambda)=\Psi_F(N_i)=\deg F_i$.

For a root sequence $z$ comparable with $(a,b)$ put $P_z(t)=(t-z_1)(t-z_2)$. If $z<(a,b)$ then $z_1\le a$, so the column $a$ of $\gamma^z$ has its entry in row $0$, equal to $P_z(a)$; and $z_2\le b$, so the column $b$ has its entry in row $0$ or $1$, and $\gamma^z_{0,b}-\gamma^z_{1,b}=P_z(b)$ in both cases. If $z\ge(a,b)$ then $\gamma^z_{0,a}=0$ and $\gamma^z_{0,b}=\gamma^z_{1,b}=0$. Since $P_{(a,b)}(a)=P_{(a,b)}(b)=0$, summing over the chain gives
\begin{equation}\label{eq:h0chi}
h^0(F(a))=\chi_x(a),\qquad h^0(F(b))-h^1(F(b))=\chi_x(b).
\end{equation}
The polynomials $\chi_x$ and $\chi(F_i(\,\cdot\,))$ have the same coefficients of $t^2$ and $t$, because $x$ and $F_i$ have the same rank and degree; so $D=\chi_x-\chi(F_i(\,\cdot\,))$ is a constant, equal to $\operatorname{ch}_2(x)-\operatorname{ch}_2(F_i)$.

Put $F'=F_i$, $F''=F/F_i$, which is torsion-free, and write $h^j_t(\cG)=h^j(\cG(t))$. By Lemma~\ref{lem:legendre}(i), $\Lambda_F(\lambda)=\Psi_F(N_i)-N_i\lambda=N_i(\mu(F_i)-\lambda)$, so the four terms of Proposition~\ref{prop:fourdefect} vanish; in particular $\rk\partial^0_a=h^1_a(F')$ and $\rk\partial^1_b=h^1_b(F'')$. From the long exact sequences and \eqref{eq:h0chi},
\begin{align*}
D&=h^0_a(F)-\chi(F'(a))=h^0_a(F')+h^0_a(F'')-\rk\partial^0_a-\chi(F'(a))=h^0_a(F'')-h^2_a(F'),\\
D&=h^0_b(F)-h^1_b(F)-\chi(F'(b))=h^0_b(F'')-h^1_b(F'')+\rk\partial^1_b-h^2_b(F')=h^0_b(F'')-h^2_b(F').
\end{align*}
By Theorem~\ref{thm:main}, $\Psi_F\ge P_F$ with equality at $N_i$, so the slopes of $\Psi_F$ at $N_i$ lie between the slopes $\mu_{i+1}$ and $\mu_i$ of $P_F$ there; in particular $\mu_{i+1}\le\lambda\le\mu_i$. Now $h^0_b(F'')=0$: a nonzero section of $F''(b)$ gives a nonzero map $\OO(-b)\to F''$, so $-b\le\mu_{\max}(F'')=\mu_{i+1}\le\lambda$, while $b<c=-\lambda-\frac32$ gives $-b>\lambda$. And $h^2_a(F')=h^0\big(F'^\vee(-a-3)\big)=0$: a nonzero section gives $\OO(a+3)\to F'^\vee$, so $a+3\le\mu_{\max}(F'^\vee)=-\mu_{\min}(F')=-\mu_i$, while $a>c$ gives $-a-3<-c-3=\lambda-\frac32<\mu_i$. Hence $D=h^0_a(F'')\ge0$ and $D=-h^2_b(F')\le0$, so $D=0$, $h^0_a(F'')=0$ and $h^2_b(F')=0$.
\end{proof}

\begin{lemma}[gap and discriminant]\label{lem:gapdisc}
Let $z=(a,b)$ be a root sequence with gap $g=a-b$, and let $\cQ'$ be a $\mu$-semistable torsion-free sheaf on $\PP^2$ of slope $\mu(z)$ with $H^0(\cQ'(a))=0$. Then
\[
\delta(\cQ')\ \ge\ \frac{g^2-1}4=\delta(E_z),
\]
with equality if and only if $H^1(\cQ'(a))=0$; here $E_z$ is a supernatural bundle with table $\gamma^z$.
\end{lemma}

\begin{proof}
By Serre duality $H^2(\cQ'(a))$ is dual to $H^0(\cQ'^\vee(-a-3))$, and $\cQ'^\vee$ is $\mu$-semistable of slope $-\mu(z)$. A nonzero section would give $a+3\le-\mu(z)=\frac{a+b}2+\frac32$, i.e.\ $a+3\le b$, which is false. So $\chi(\cQ'(a))=-h^1(\cQ'(a))\le0$. By \eqref{eq:RRP2}, with $a+\mu(z)+\frac32=a-c(z)=\frac g2$,
\[
\frac{\chi(\cQ'(a))}{\rk\cQ'}=\frac{g^2}8-\frac18-\frac{\delta(\cQ')}2\ \le\ 0 .
\]
For $E_z$, Lemma~\ref{lem:basic}(ii) gives $4\delta(E_z)=G(E_z)-1=g^2-1$.
\end{proof}

\begin{proposition}[contact edges]\label{prop:contactedge}
Let $F$ be a vector bundle on $\PP^2$, and suppose that $\Psi_F=P_F$ on the edge $[N_{i-1},N_i]$ of $P_F$. Then $\cQ_i=F_i/F_{i-1}$ is a vector bundle whose table is a multiple of $\gamma^z$, where $z$ is the unique element of the chain of $F$ of slope $\mu_i$.
\end{proposition}

\begin{proof}
Two comparable root sequences with the same center are equal, so $z=(a,b)$ is unique, and it is the only pair adapted to the chain with $a+b=-2\mu_i-3$. Put $g=a-b$ and $r=\rk\cQ_i$. Since $\Psi_F$ is linear of slope $\mu_i$ on the edge, $\mu_i$ is a slope of $\Psi_F$ at $N_{i-1}$ and at $N_i$. At an interior endpoint, Lemma~\ref{lem:Ktight} with $\lambda=\mu_i$ and the pair $(a,b)$ identifies the Chern character of $F_{i-1}$, resp.\ $F_i$, with that of the prefix of the chain consisting of the pieces of slope $>\mu_i$ and a fraction of $\gamma^z$; at the endpoints $0$ and $\rk F$ the prefixes are $0$ and $\gamma(F)$, and the pieces of slope $\ne\mu_i$ do not occur between them because $\Psi_F$ has slope $\mu_i$ on the whole edge. Subtracting, $\chi(\cQ_i(t))=\frac r2(t-a)(t-b)$, and $\delta(\cQ_i)=\frac{g^2-1}4$.

If $i\ge2$, Lemma~\ref{lem:Ktight} at $F_{i-1}$ gives $H^0((F/F_{i-1})(a))=0$, hence $H^0(\cQ_i(a))=0$; if $i=1$, $z$ is the last element of the chain, every element $z^k$ has $z^k_1\ge a$, so $H^0(F(a))=0$ and again $H^0(\cQ_1(a))=0$. Similarly, if $i<s$, Lemma~\ref{lem:Ktight} at $F_i$ gives $H^2(F_i(b))=0$, and $H^2(F_i(b))\to H^2(\cQ_i(b))$ is surjective; if $i=s$, $z$ is the first element of the chain, every $z^k$ has $z^k_2\le b$, so $H^2(F(b))=0$, which surjects onto $H^2(\cQ_s(b))$. Since $b<a$, and $h^0(\cQ_i(t))$ is nondecreasing and $h^2(\cQ_i(t))$ nonincreasing in $t$ for the torsion-free sheaf $\cQ_i$ (as in the proof of Proposition~\ref{prop:fourdefect}), also $H^0(\cQ_i(b))=0$ and $H^2(\cQ_i(a))=0$. As $\chi(\cQ_i(a))=\chi(\cQ_i(b))=0$, the sheaf $\cQ_i$ has no cohomology at the twists $a$ and $b$.

Let $E=\cQ_i^{\vee\vee}$ and $\ell=\operatorname{length}(E/\cQ_i)$, so that $\delta(E)=\delta(\cQ_i)-2\ell/r$. The dual $E^\vee=\cQ_i^\vee$ is semistable of slope $-\mu_i$, which is the slope of the root sequence $(-b-3,-a-3)$, of gap $g$, and $H^0(E^\vee(-b-3))$ is dual to $H^2(E(b))\cong H^2(\cQ_i(b))=0$, since $E$ is locally free and $E/\cQ_i$ has finite support. Lemma~\ref{lem:gapdisc} gives $\delta(E)=\delta(E^\vee)\ge\frac{g^2-1}4=\delta(\cQ_i)$, so $\ell=0$ and $\cQ_i$ is locally free. The columns $a$ and $b$ of its table vanish, so every element of its chain is $(a,b)$.
\end{proof}

\begin{corollary}\label{cor:equalityP2}
For a vector bundle $F$ on $\PP^2$, $P_F=\Psi_F$ if and only if the Harder--Narasimhan factors of $F$ are vector bundles whose tables are multiples of supernatural tables $\gamma^{z(1)},\dots,\gamma^{z(s)}$ with pairwise comparable root sequences. In that case $\gamma(F)=\sum_i\gamma(\cQ_i)$ is the chain decomposition.
\end{corollary}

\begin{proof}
If $P_F=\Psi_F$, every edge is a contact edge, and Proposition~\ref{prop:contactedge} applies; the root sequences belong to the chain of $F$, so they are comparable. The converse holds on every $\PP^n$ (Remark~\ref{rem:main}(1)).
\end{proof}

Lemma~\ref{lem:Ktight} says that the equality cases of Theorem~\ref{thm:main} on $\PP^2$ are visible in $K$-theory: at a tight vertex the pieces above the vertex have the Chern character of the Harder--Narasimhan step, and in particular the same discriminant. Lemma~\ref{lem:gapdisc} is a discriminant bound under a vanishing condition, valid in any characteristic: a $\mu$-semistable sheaf $\cQ'$ of slope $\mu(z)$ with $H^0(\cQ'(a))=0$ has discriminant at least that of the supernatural bundles with roots $z$.

\section{Bridgeland stability and the large volume limit}\label{sec:lvl}

On $\PP^2$ the Harder--Narasimhan polygon is the large volume limit of the Harder--Narasimhan polygons of Bridgeland stability. In this section we show that Theorem~\ref{thm:main} survives the limit: for the stability conditions $\sigma_{\alpha,\beta}$ with $\alpha\gg0$, the Bridgeland Harder--Narasimhan polygon of a vector bundle lies inside a polygon built from its Boij--Söderberg decomposition. The leading order is Theorem~\ref{thm:main}; the second order rests on the results of Section~\ref{sec:equality}: the equality cases are determined at the level of Chern characters (Lemma~\ref{lem:Ktight}), and a Harder--Narasimhan factor along a contact edge is a supernatural bundle (Proposition~\ref{prop:contactedge}, which rests on a Bogomolov-type bound for the gap, Lemma~\ref{lem:gapdisc}). For Theorem~\ref{thm:lvl} and in Remark~\ref{rem:lvl} the ground field is $\mathbb C$, as in the references.

\subsection{Stability conditions and polygons}
We use the stability conditions $\sigma_{\alpha,\beta}$ of \cite{Br08,AB13,ABCH13}, with heart $\operatorname{Coh}^\beta(\PP^2)$ and central charge $Z_{\alpha,\beta}$; see \cite{MS17} for an exposition. For $\beta\in\RR$, the heart $\operatorname{Coh}^\beta(\PP^2)$ is the tilt of $\operatorname{Coh}(\PP^2)$ at the torsion pair $(\mathcal T^\beta,\mathcal F^\beta)$: $\mathcal T^\beta$ consists of the sheaves whose torsion-free part has $\mu_{\min}>\beta$, and $\mathcal F^\beta$ of the torsion-free sheaves with $\mu_{\max}\le\beta$. For $\alpha>0$ the central charge is $Z_{\alpha,\beta}(E)=-\int_{\PP^2}e^{-(\beta+i\alpha)H}\operatorname{ch}(E)$, i.e.
\[
Z_{\alpha,\beta}(E)=-\Big(\operatorname{ch}_2(E)-\beta c_1(E)+\frac{\beta^2-\alpha^2}2\,r(E)\Big)+i\alpha\big(c_1(E)-\beta\,r(E)\big),
\]
which depends only on the class $v(E)=(r,c_1,\operatorname{ch}_2)$ of $E$. Every nonzero object $E$ of the heart has $Z_{\alpha,\beta}(E)\in\RR_{>0}\,e^{i\pi\theta}$ with a phase $\theta=\theta(E)\in(0,1]$, and a Harder--Narasimhan filtration with $\sigma_{\alpha,\beta}$-semistable factors $E_1,\dots,E_m$ of strictly decreasing phases \cite{Br07}. The \emph{Harder--Narasimhan polygon} $P^\sigma_E$ is the polygonal path with vertices $0,\ Z(E_1),\ Z(E_1)+Z(E_2),\ \dots,\ Z(E)$; together with the segment $[0,Z(E)]$ it bounds a convex region, which we denote $|P^\sigma_E|$. A vector bundle $F$ with $\mu_{\min}(F)>\beta$ lies in $\operatorname{Coh}^\beta(\PP^2)$.

The Boij--Söderberg decomposition gives a second polygon in the same plane. A supernatural bundle with table $\gamma^z$ has Hilbert polynomial $(t-z_1)(t-z_2)$, so by Riemann--Roch its class is
\[
v(z)=\big(2,\ 2\mu(z),\ z_1z_2-3\mu(z)-2\big),\qquad \mu(z)=-\tfrac{z_1+z_2}2-\tfrac32 .
\]

\begin{definition}\label{def:PsiSigma}
Let $F$ be a vector bundle on $\PP^2$ with chain decomposition $\gamma(F)=\sum_kq_k\gamma^{z^k}$, let $\sigma=\sigma_{\alpha,\beta}$, and suppose that the vectors $q_kZ_{\alpha,\beta}(v(z^k))$ lie in an open half-plane $\mathcal H$ bounded by a line through $0$. The \emph{Boij--Söderberg polygon} $\Psi^\sigma_F$ is the polygonal path from $0$ to $Z_{\alpha,\beta}(F)$ whose edges are these vectors, arranged in order of decreasing argument, measured in an interval of length $\pi$ corresponding to $\mathcal H$; $|\Psi^\sigma_F|$ is the convex region it bounds with the segment $[0,Z_{\alpha,\beta}(F)]$.
\end{definition}

Without the half-plane condition the path need not bound a region: for $F=\OO(2)\oplus\OO\oplus\OO(-2)$ and $(\alpha,\beta)=(1,0)$ the edges are $-\frac32+2i$, $\frac12$ and $-\frac32-2i$, and the first and the last edge meet at $-\frac54+\frac53i$. The condition holds in the two cases we use. For $\alpha\gg0$ all edges have positive real part, by \eqref{eq:Zxy} below; and if $\beta<\mu^\BS_{\min}(F)=\min_k\mu(z^k)$, all edges have positive imaginary part. In the second case the supernatural bundles $E_{z^k}$ lie in $\operatorname{Coh}^\beta(\PP^2)$, and if they are $\sigma$-semistable with distinct phases, then $\Psi^\sigma_F$ is the Harder--Narasimhan polygon of the supernatural model $\bigoplus_kE_{z^k}^{\oplus m q_k}$, up to the factor $m$. Note that $\mu^\BS_{\min}(F)\le\mu_{\min}(F)$ by Corollary~\ref{cor:supernatural}(b), so the hypothesis $\beta<\mu_{\min}(F)$ of Theorem~\ref{thm:lvl} below does not imply $\beta<\mu^\BS_{\min}(F)$.

\begin{theorem}[large volume limit]\label{thm:lvl}
Let $F$ be a vector bundle on $\PP^2$ over $\mathbb C$ and $\beta<\mu_{\min}(F)$ a real number. There is $\alpha_0$ such that for all $\alpha\ge\alpha_0$, with $\sigma=\sigma_{\alpha,\beta}$,
\[
|P^\sigma_F|\ \subseteq\ |\Psi^\sigma_F| .
\]
\end{theorem}

For $r>0$ and $v=(r,c_1,\operatorname{ch}_2)$ write $\mu=c_1/r$ and $\operatorname{ch}_2^\beta=\operatorname{ch}_2-\beta c_1+\beta^2r/2$. Then
\begin{equation}\label{eq:Zxy}
\operatorname{Re}Z_{\alpha,\beta}(v)=\frac{\alpha^2r}2-\operatorname{ch}_2^\beta,\qquad
\operatorname{Im}Z_{\alpha,\beta}(v)=\alpha r(\mu-\beta),\qquad
\frac{\operatorname{Im}Z}{\operatorname{Re}Z}=\frac{2(\mu-\beta)}{\alpha-2\operatorname{ch}_2^\beta/(\alpha r)} .
\end{equation}
So for a fixed finite set of classes of positive rank and $\alpha\gg0$, the real parts are positive, the order of the arguments is the order of the slopes $\mu$, and among classes with the same $\mu>\beta$ the argument increases with $\operatorname{ch}_2/r$: on classes with $\mu>\beta$ the order of phases is the Gieseker order. (For $\mu<\beta$ the second comparison is reversed, but it is not needed: no two elements of a chain have the same slope.) The large volume limit enters through the following lemma, the form in which we use the analogue for surfaces of \cite[Prop.~14.2]{Br08} (see \cite[\S6, Exercise~6.27]{MS17}, and \cite[\S6]{ABCH13} for $\PP^2$).

\begin{lemma}[large volume limit]\label{lem:LV}
Let $\beta\in\RR$, and let $\cG$ be a Gieseker semistable torsion-free sheaf on $\PP^2$ with $\mu(\cG)>\beta$. Then $\cG$ lies in $\operatorname{Coh}^\beta(\PP^2)$, and there is a number $M$, depending only on $\beta$ and on the class $v(\cG)$, such that $\cG$ is $\sigma_{\alpha,\beta}$-semistable for all $\alpha\ge M$.
\end{lemma}

\begin{proof}
We use \cite{LQ14} for $X=\PP^2$, $\omega=H$ and the class $\beta H$. In the notation there, the heart $\mathcal A^\sharp_{(\omega,\beta\omega)}$ is $\operatorname{Coh}^\beta(\PP^2)$, and the stability condition $(Z_m,\mathcal P_m)$ with $m=\alpha$ is $\sigma_{\alpha,\beta}$ \cite[(2.5)]{LQ14}. An object $E$ of the heart is semistable for the polynomial stability condition $(Z_\Omega,\mathcal P_\Omega)$ of Bayer \cite{Bay09} if for every subobject $A\subset E$ in the heart the phase of $Z_{\alpha,\beta}(A)$ is at most the phase of $Z_{\alpha,\beta}(E)$ for $\alpha\gg0$, where the bound on $\alpha$ may depend on $A$; this is how the condition is used in \cite{LQ14}, see Remark~2.8 and the proofs of Lemmas~2.11 and~5.1 there. By \cite[Lemma~4.3]{LQ14}, which is stated for $\omega$ rational and $\beta$ real, a $(Z_\Omega,\mathcal P_\Omega)$-semistable object $E$ is $\sigma_{\alpha,\beta}$-semistable for all $\alpha\ge M$, where $M$ depends only on $\beta$ and on the class of $E$. So it suffices to show that $\cG$ is $(Z_\Omega,\mathcal P_\Omega)$-semistable.

Since $\cG$ is $\mu$-semistable with $\mu(\cG)>\beta$, it lies in $\mathcal T^\beta$. Let $A\subset\cG$ be a nonzero proper subobject in $\operatorname{Coh}^\beta(\PP^2)$, with quotient $B$, and let $\mathcal H^{-1}(B)$, $\mathcal H^0(B)$ be the cohomology sheaves of $B$. The long exact sequence of cohomology sheaves shows that $A$ is a sheaf in $\mathcal T^\beta$ and gives an exact sequence of sheaves
\[
0\to\mathcal K\to A\to\cG\to\mathcal H^0(B)\to0,\qquad \mathcal K=\mathcal H^{-1}(B)\in\mathcal F^\beta .
\]
Let $I\subset\cG$ be the image of $A$. The torsion of $A$ maps to zero in $\cG$, so it lies in the torsion-free sheaf $\mathcal K$; hence $A$ is torsion-free. Moreover $I\ne0$, because $\mathcal T^\beta\cap\mathcal F^\beta=0$; so $A$ has positive rank, and $\mu(A)>\beta$ because $A\in\mathcal T^\beta$. By \eqref{eq:Zxy} and the discussion after it, the phase of $A$ is at most that of $\cG$ for $\alpha\gg0$ if $\mu(A)<\mu(\cG)$, and also if $\mu(A)=\mu(\cG)$ and $\operatorname{ch}_2(A)/r(A)\le\operatorname{ch}_2(\cG)/r(\cG)$. If $\mathcal K\ne0$, then $\mu(\mathcal K)\le\beta<\mu(\cG)$, and $\mu(I)\le\mu(\cG)$ because $\cG$ is $\mu$-semistable; so $\mu(A)<\mu(\cG)$. If $\mathcal K=0$, then $A$ is a subsheaf of $\cG$, and the Gieseker semistability of $\cG$ says that $\mu(A)<\mu(\cG)$, or $\mu(A)=\mu(\cG)$ and $\operatorname{ch}_2(A)/r(A)\le\operatorname{ch}_2(\cG)/r(\cG)$, since by Riemann--Roch the reduced Hilbert polynomial of a torsion-free sheaf $\cE$ of rank $r$ on $\PP^2$ is $\frac{t^2}2+\big(\mu(\cE)+\frac32\big)t+\frac{\operatorname{ch}_2(\cE)}r+\frac32\mu(\cE)+1$. In both cases the phase of $A$ is at most that of $\cG$ for $\alpha\gg0$, so $\cG$ is $(Z_\Omega,\mathcal P_\Omega)$-semistable. For coprime rank and degree this is \cite[Lemma~5.1]{LQ14}.
\end{proof}

In \cite{Br08,MS17} the large volume limit is formulated with twisted Gieseker stability. On $\PP^2$ the twist plays no role: for every $b\in\RR$, a torsion-free sheaf is $bH$-twisted Gieseker semistable if and only if it is Gieseker semistable. Indeed, the $bH$-twisted reduced Hilbert polynomial of a torsion-free sheaf $\cE$ of rank $r>0$ is $\frac1r\int\operatorname{ch}(\cE)e^{-bH}e^{tH}\operatorname{td}(\PP^2)$. Its coefficient of $t$ is $\mu(\cE)-b+\frac32$, and for fixed $\mu(\cE)$ its constant term is $\operatorname{ch}_2(\cE)/r$ plus a number depending only on $\mu(\cE)$ and $b$. So twisted and untwisted reduced Hilbert polynomials compare $\mu$ first and then, for equal $\mu$, $\operatorname{ch}_2/r$, in the same way. The argument is unchanged if $\operatorname{td}(\PP^2)$ is replaced by $\sqrt{\operatorname{td}(\PP^2)}$, as in some definitions of twisted stability.

By Lemma~\ref{lem:LV}, if $F$ is a vector bundle with $\mu_{\min}(F)>\beta$ and $0=\cG_0\subset \cG_1\subset\dots\subset \cG_p=F$ is its Gieseker Harder--Narasimhan filtration, then for $\alpha\gg0$, depending on the finitely many classes of the factors, the factors $\cG_j/\cG_{j-1}$ lie in the heart, are $\sigma_{\alpha,\beta}$-semistable, and have strictly decreasing phases; so this filtration is the $\sigma_{\alpha,\beta}$-Harder--Narasimhan filtration of $F$, and $P^\sigma_F$ has vertices $Z(\cG_j)$. The Gieseker filtration refines the $\mu$-Harder--Narasimhan filtration $0=F_0\subset\dots\subset F_s=F$: on each $\mu$-semistable factor $\cQ_i=F_i/F_{i-1}$ it induces the Gieseker filtration of $\cQ_i$, whose factors are $\mu$-semistable of slope $\mu_i$ with decreasing $\operatorname{ch}_2/r$.

\subsection{Proof of Theorem~\ref{thm:lvl}}
We use the notation $N_i=\rk F_i$ and the tight vertices of Section~\ref{sec:equality}. Fix $\beta<\mu_{\min}(F)$ and write $Z=Z_{\alpha,\beta}$, $x(v)=\operatorname{Re}Z(v)$, $y(v)=\operatorname{Im}Z(v)$. Let $\alpha$ be large enough that $P^\sigma_F$ is the Gieseker polygon, with vertices $Z(\cG_j)$, and that the arguments of the edges of both polygons are ordered as described after \eqref{eq:Zxy}. Then the edges of $\Psi^\sigma_F$ are the pieces in decreasing order of slope. All edges of both polygons have $x>0$ and decreasing $y/x$, so both paths are graphs of concave piecewise linear functions $p$ and $\psi$ on $[0,x(F)]$, and $|P^\sigma_F|\subseteq|\Psi^\sigma_F|$ if and only if $p\le\psi$.

Since $\psi$ is concave and $p$ is linear between consecutive vertices $Z(\cG_j)$, it suffices to show that each of the finitely many vertices $Z(\cG_j)$ lies on or below the path $\Psi^\sigma_F$ for $\alpha$ large; then one $\alpha_0$ works for all of them. Put $r_j=\rk \cG_j$; since the Gieseker filtration refines the Harder--Narasimhan filtration, $\deg \cG_j=P_F(r_j)$.

The affine map $\Phi(x,y)=\big(2x/\alpha^2,\ y/\alpha+2\beta x/\alpha^2\big)$ sends $Z(v)$ to $(r,c_1)+O(\alpha^{-2})$, uniformly for $v$ in the finite set of classes of the Gieseker factors and of the pieces, and it preserves the vertical order of points with the same first coordinate. So $\Phi(Z(\cG_j))=(r_j,P_F(r_j))+O(\alpha^{-2})$, and $\Phi(\Psi^\sigma_F)$, as a graph over the first coordinate, lies within $O(\alpha^{-2})$ of $\Psi_F$. Hence if $\Psi_F(r_j)>P_F(r_j)$, the vertex $Z(\cG_j)$ lies strictly below $\Psi^\sigma_F$ for $\alpha$ large.

Suppose $\Psi_F(r_j)=P_F(r_j)$. If $r_j$ lies in the interior of an edge $[N_{i-1},N_i]$ of $P_F$, then $\Psi_F=P_F$ on that edge by Lemma~\ref{lem:contact}, so $\cQ_i$ is a supernatural bundle by Proposition~\ref{prop:contactedge}, hence Gieseker semistable by Proposition~\ref{prop:sameslope}; then the Gieseker filtration has no step strictly between $F_{i-1}$ and $F_i$, a contradiction. So $\cG_j=F_i$ and $(N_i,\deg F_i)$ is a tight vertex. If $0<i<s$, Lemma~\ref{lem:Ktight} gives $v(F_i)=v(x)$ for a prefix $x$ of the chain, consisting of the pieces of slope $>\lambda$ and a fraction of the piece of slope $\lambda$, for a slope $\lambda$ of $\Psi_F$ at $N_i$; since the edges of $\Psi^\sigma_F$ are the pieces in order of slope, $Z(F_i)=Z(x)$ lies on $\Psi^\sigma_F$. The endpoints $0$ and $Z(F)$ lie on both paths. This proves $p\le\psi$ on $[0,x(F)]$ for $\alpha\ge\alpha_0$. \qed

\begin{lemma}[support functions]\label{lem:support}
Let $\sigma=\sigma_{\alpha,\beta}$, identify $\mathbb C$ with $\RR^2$, with $\langle u,w\rangle=\operatorname{Re}(\bar uw)$, and for $X$ in $\operatorname{Coh}^\beta(\PP^2)$ let $h_X$ be the support function of $|P^\sigma_X|$. For $\theta\in(0,1)$ put $\ell_\theta=\operatorname{Im}(e^{-i\pi\theta}\,\cdot\,)=\langle e^{i\pi(\theta+1/2)},\cdot\,\rangle$, and let $X_{>\theta}\subset X$ be the step of the Harder--Narasimhan filtration with factors of phase $>\theta$.
\begin{enumerate}[label=\textup{(\alph*)}]
\item $h_X(e^{i\pi(\theta+1/2)})=\ell_\theta(Z(X_{>\theta}))$, and $h_X(u)=\max(0,\langle u,Z(X)\rangle)$ if $\operatorname{Re}u\ge0$. In particular $h_{X\oplus Y}(u)=h_X(u)+h_Y(u)$ for $\operatorname{Re}u\le0$.
\item For every subobject $S\subset X$, $\ell_\theta(Z(S))\le\ell_\theta(Z(X_{>\theta}))$.
\item If $0\to A\to B\to C\to0$ is exact in $\operatorname{Coh}^\beta(\PP^2)$, then $|P^\sigma_{B\oplus D}|\subseteq|P^\sigma_{A\oplus C\oplus D}|$ for every $D$ in $\operatorname{Coh}^\beta(\PP^2)$.
\end{enumerate}
\end{lemma}

\begin{proof}
(a) Along the path $P^\sigma_X$ the increments of $\ell_\theta$ are first positive and then $\le0$, so the maximum of $\ell_\theta$ on $|P^\sigma_X|$ is $\ell_\theta(Z(X_{>\theta}))$; since $(X\oplus Y)_{>\theta}=X_{>\theta}\oplus Y_{>\theta}$, $h_{X\oplus Y}=h_X+h_Y$ in these directions. In the remaining directions $u$, those with $\operatorname{Re}u\ge0$, the increments of $\langle u,\cdot\,\rangle$ along the path are first $\le0$ and then $\ge0$, so $h_X(u)=\max(0,\langle u,Z(X)\rangle)$ depends only on $Z(X)$. Additivity on the boundary $\operatorname{Re}u=0$ follows by continuity from the left half-plane; no additivity is asserted in the right half-plane. (b) We have $S_{>\theta}\subset X_{>\theta}$, since $\operatorname{Hom}(S_{>\theta},X/X_{>\theta})=0$, and $X_{>\theta}/S_{>\theta}$ has factors of phase $>\theta$, while $S/S_{>\theta}$ has factors of phase $\le\theta$; so $\ell_\theta(Z(S))\le\ell_\theta(Z(S_{>\theta}))\le\ell_\theta(Z(X_{>\theta}))$. (c) Applied to $B_{>\theta}\cap A\subset A$ and to its image in $C$, (b) gives $h_B\le h_A+h_C=h_{A\oplus C}$ in the directions $e^{i\pi(\theta+1/2)}$, and $h_B=h_{A\oplus C}$ in the other directions, since $Z(B)=Z(A\oplus C)$. For $\operatorname{Re}u<0$, adding $h_D(u)$ and using (a) gives $h_{B\oplus D}(u)\le h_{A\oplus C\oplus D}(u)$. For $\operatorname{Re}u\ge0$, the two support functions are equal by (a), since $Z(B\oplus D)=Z(A\oplus C\oplus D)$. The support-function criterion for containment of compact convex sets proves the assertion.
\end{proof}

\begin{remark}\label{rem:lvl}
\begin{enumerate}[label=(\arabic*)]
\item The statement does not extend to all $(\alpha,\beta)$ at the level of the Gieseker polygon. Let $S=\coker(\OO(-2)^2\to\OO(-1)^2\oplus\OO^2)$ for a general map, a stable bundle of slope $1$ with $h^0(S(t))=2\big(h^0(\OO(t-1))+h^0(\OO(t))-h^0(\OO(t-2))\big)$, $h^2(S(t))=h^0(S(-5-t))$ and $\chi(S(t))=t^2+5t+2$, and let $F=\OO(2)\oplus S\oplus T(-2)$, where here and below $T=T_{\PP^2}$ is the tangent bundle. Its Harder--Narasimhan filtration is $\OO(2)\subset\OO(2)\oplus S\subset F$, and its chain is
\[
\gamma^{(0,-2)}+\tfrac35\gamma^{(0,-4)}+\tfrac1{20}\gamma^{(0,-5)}+\tfrac{13}{30}\gamma^{(-1,-5)}+\tfrac1{12}\gamma^{(-2,-5)}+\tfrac13\gamma^{(-3,-5)} .
\]
The vertex at rank $3$ is tight, and $v(\OO(2)\oplus S)=(3,4,-1)$ is the class of the five pieces of slope $\ge\frac12$, which are the pieces above this vertex in the sense of Lemma~\ref{lem:Ktight} (with $\lambda=\frac12$ and $f=\frac35$, the whole piece $\gamma^{(0,-4)}$), as that lemma predicts. The vertex at rank $1$ is not: $\Psi_F(1)=\frac94>2$. For $(\alpha,\beta)=(\frac15,-\frac{11}{20})$, where all edges of $\Psi^\sigma_F$ have positive imaginary part, neither the Gieseker polygon of $F$ nor its $\sigma$-Harder--Narasimhan polygon is contained in $|\Psi^\sigma_F|$. The sections of $S$ give an exact sequence $0\to\OO^2\to S\to Q\to0$, with $Q=\coker(\OO(-2)^2\to\OO(-1)^2)$ supported on a conic; as $Q$ is a torsion sheaf and $\mu(T(-2))=-\frac12>\beta$, both lie in $\mathcal T^\beta$, and $\cG=\OO(2)\oplus\OO^2$ is a subobject of $F$ in $\operatorname{Coh}^\beta(\PP^2)$. Its central charge $Z(\cG)=-\frac{559}{160}+\frac{73}{100}\,i$ lies strictly beyond the edge $\frac{13}{30}Z(v(-1,-5))$ of $\Psi^\sigma_F$, of phase $\theta\approx0.648$, and hence so does the vertex $Z(F_{>\theta})$ of $P^\sigma_F$, since $\ell_\theta(Z(\cG))\le\ell_\theta(Z(F_{>\theta}))$ by Lemma~\ref{lem:support}(b). (So $S$ is $\sigma$-unstable: by the sequence above, $\OO\subset S$ is a subobject in $\operatorname{Coh}^\beta(\PP^2)$, and $\OO$ and $S$ have phases $0.778$ and $0.115$.) In the same example, $v(\OO(2))=(1,2,2)$ is not of the form $\sum_kt_kq_kv(z^k)$ with $0\le t_k\le1$: the combinations with rank $1$ and degree $2$ have $\operatorname{ch}_2\le\frac{31}{16}$. So the containment for all $(\alpha,\beta)$, which would follow if every Harder--Narasimhan step were such a fractional sum of pieces, cannot be proved by an argument on classes alone.
\item The containment cannot hold for all $\sigma$ in the geometric chamber $\{\sigma_{\alpha,\beta}:\alpha>0,\ \beta\in\RR\}$, because the pieces themselves can cross walls. For $F=T(-1)$, with $\gamma(F)=\gamma^{(-1,-3)}$, the polygon $\Psi^\sigma_F$ is the segment $[0,Z(T(-1))]$ for every $\sigma$. But for $-1\le\beta<0$ the Euler sequence becomes the exact sequence $0\to\OO^3\to T(-1)\to\OO(-1)[1]\to0$ in $\operatorname{Coh}^\beta(\PP^2)$, and inside the semicircle $(\beta+\frac12)^2+\alpha^2<\frac14$ the subobject $\OO^3$ has larger phase than $T(-1)$ (at $(\alpha,\beta)=(\frac3{10},-\frac12)$ the phases are $0.656$ and $0.583$), so $T(-1)$ is unstable and $|P^\sigma_{T(-1)}|$ has nonempty interior. Below the walls of the pieces, the polygon to compare with is the Harder--Narasimhan polygon of a supernatural model, which is an object of $\operatorname{Coh}^\beta(\PP^2)$ when $\beta<\mu^\BS_{\min}(F)$. We ask: for a vector bundle $F$ on $\PP^2$ with chain $\sum_kq_k\gamma^{z^k}$, are there supernatural bundles $E_{z^k}$ with tables $\gamma^{z^k}$, and $m$ with $mq_k\in\ZZ$, such that
\[
|P^\sigma_F|\ \subseteq\ \tfrac1m\,\big|P^\sigma_{\bigoplus_kE_{z^k}^{\oplus mq_k}}\big|
\]
for every $\sigma=\sigma_{\alpha,\beta}$ with $\alpha>0$ and $\beta<\mu^\BS_{\min}(F)$? For $\alpha\gg0$ this holds for every choice, by Theorem~\ref{thm:lvl}.
\item The answer is yes for $F=\OO\oplus E$, where $E=\coker(\OO(-2)^2\to\OO(-1)^4)$ is a stable bundle with $c_1=0$, $c_2=2$ and table $\gamma^{(0,-3)}$. The chain of $F$ is $\frac13\gamma^{(0,-2)}+\frac56\gamma^{(0,-3)}+\frac13\gamma^{(-1,-3)}$, so with $m=6$ a supernatural model is $T(-2)^{\oplus2}\oplus E^{\oplus5}\oplus T(-1)^{\oplus2}$, and $\mu^\BS_{\min}(F)=-\frac12$. The four rows of the presentation of $E$ are linearly independent in $\operatorname{Hom}(\OO(-2)^2,\OO(-1))\cong H^0(\OO(1))^{\oplus2}$, because $\operatorname{Hom}(E,\OO(-1))=H^0(E(-1))=0$; here and below we use $E^\vee\cong E$, as $E$ has rank $2$ and $c_1=0$. Completing them to a basis gives a presentation $\OO(-2)^2\to\OO(-1)^6$ of $T(-2)^{\oplus2}$, and projecting to the first four coordinates gives an exact sequence $0\to\OO(-1)^2\to T(-2)^{\oplus2}\to E\to0$. Since $\operatorname{Ext}^1(E,\OO)=H^1(E)=0$, its class is the image of a map $E\to T(-1)^{\oplus2}$ under the connecting map of the Euler sequence $0\to\OO(-1)^2\to\OO^6\to T(-1)^{\oplus2}\to0$; so $T(-2)^{\oplus2}$ is the fiber product of $\OO^6$ and $E$ over $T(-1)^{\oplus2}$, and
\[
0\to T(-2)^{\oplus2}\to\OO^{6}\oplus E\to T(-1)^{\oplus2}\to0
\]
is exact. For $\beta<-\frac12$ all its terms lie in $\mathcal T^\beta$, so it is exact in $\operatorname{Coh}^\beta(\PP^2)$. By Lemma~\ref{lem:support}(c), with $B=\OO^6\oplus E$, $A=T(-2)^{\oplus2}$, $C=T(-1)^{\oplus2}$ and $D=E^{\oplus5}$, we get $6\,|P^\sigma_F|=|P^\sigma_{F^{\oplus6}}|\subseteq|P^\sigma_{T(-2)^{\oplus2}\oplus E^{\oplus5}\oplus T(-1)^{\oplus2}}|$ for all $\alpha>0$ and $\beta<-\frac12$.
\end{enumerate}
\end{remark}

\section{Sheaves of lower dimension}\label{sec:lowerdim}

Let $1\le n\le N$. For a coherent sheaf $\cF$ of dimension at most $n$ on $\PP^N$ write
\[
\chi(\cF(t))=\sum_{i=0}^n(-1)^i\,e_i(\cF)\binom{t+n-i}{n-i},
\]
and put $\rk_n\cF=e_0(\cF)$ and $\deg_n\cF=-e_1(\cF)$. Both are additive on exact sequences. The number $\rk_n\cF\ge0$ is the multiplicity of the $n$-dimensional part of $\cF$; if it vanishes, $\deg_n\cF\ge0$ is the multiplicity of the $(n-1)$-dimensional part. Since $\binom{t+n}n=\frac{t^n}{n!}+\frac{n+1}2\,\frac{t^{n-1}}{(n-1)!}+\dots$, for $N=n$ Riemann--Roch (Section~\ref{sec:tables}) gives $\rk_n\cF=\rk\cF$ and $\deg_n\cF=\deg\cF$. In general we use the normalized slope $\mu_n(\cF)=\deg_n\cF/\rk_n\cF=\hat\mu(\cF)-\frac{n+1}2$, where $\hat\mu$ is the slope of \cite[\S1.6]{HL10}. Let $P_\cF$ be the upper concave envelope of the points $(\rk_n\cF',\deg_n\cF')$ for coherent subsheaves $\cF'\subset\cF$; for $\cF$ pure of dimension $n$ it is the Harder--Narasimhan polygon of $\cF$ with respect to $\mu_n$. Equivalently, if $P^{\hat\mu}_{\cF}$ denotes the polygon with slopes $\hat\mu$ and horizontal lengths $e_0$, then $P_\cF(m)=P^{\hat\mu}_\cF(m)-\frac{n+1}2m$ \cite[\S1.6]{HL10}.

The table $\gamma(\cF)$ has rows $0,\dots,n$. Root sequences $z=(z_1>\dots>z_n)$, their centers $c(z)$ and slopes $\mu(z)$, and the functionals $\phi_z$ of \eqref{eq:phi} are those of Sections~\ref{sec:tables} and~\ref{sec:main} for $\PP^n$, and, as in Section~\ref{sec:coherent}, we put
\begin{align*}
\Lambda_\cF(\lambda)&=\inf\big\{\phi_z(\gamma(\cF)):\ \mu(z)=\lambda\big\},\qquad \lambda\in\tfrac1n\ZZ,\\
\Psi_\cF(m)&=\inf_{\lambda\in\frac1n\ZZ}\big(m\lambda+\Lambda_\cF(\lambda)\big),\qquad 0\le m\le\rk_n\cF .
\end{align*}
For $N=n$ these are the functions of Section~\ref{sec:coherent}.

\begin{lemma}[projection]\label{lem:projection}
Let $\cF$ be a coherent sheaf of dimension at most $n$ on $\PP^N$, let $\Pi\subset\PP^N$ be a linear subspace of dimension $N-n-1$ disjoint from $\operatorname{Supp}\cF$ (empty if $N=n$), and let $p\colon\PP^N\setminus\Pi\to\PP^n$ be the projection from $\Pi$. On coherent sheaves $\mathcal G$ with $\operatorname{Supp}\mathcal G\subset\operatorname{Supp}\cF$ the functor $p_*$ is exact, and:
\begin{enumerate}[label=\textup{(\roman*)}]
\item $H^j(\PP^n,(p_*\mathcal G)(t))\cong H^j(\PP^N,\mathcal G(t))$ for all $j$ and $t$, functorially in $\mathcal G$; in particular $\gamma(p_*\mathcal G)=\gamma(\mathcal G)$, $\rk p_*\mathcal G=\rk_n\mathcal G$ and $\deg p_*\mathcal G=\deg_n\mathcal G$;
\item if $\mathcal G$ is pure of dimension $n$, then $p_*\mathcal G$ is torsion-free;
\item if $\mathcal G_x$ is a Cohen--Macaulay module of dimension $n$ for every closed point $x\in\operatorname{Supp}\mathcal G$, then $p_*\mathcal G$ is locally free.
\end{enumerate}
\end{lemma}

\begin{proof}
Such a $\mathcal G$ is a module over a closed subscheme $Z\subset\PP^N$ with support $\operatorname{Supp}\cF$, for instance the one defined by a power of the ideal of $\operatorname{Supp}\cF$, and $Z\subset\PP^N\setminus\Pi$. For $y\in\PP^n$ the linear space $L_y$ spanned by $\Pi$ and $y$ contains $\Pi$ as a hyperplane, and $p^{-1}(y)\cap Z=Z\cap L_y$ is projective and contained in the affine space $L_y\setminus\Pi$, hence finite. So $p|_Z$ is proper and quasi-finite, hence finite; $p_*\mathcal G=(p|_Z)_*\mathcal G$ is coherent, and $p_*$ is exact, without higher direct images. Since $p^*\OO(1)=\OO(1)$, the projection formula gives $p_*(\mathcal G(t))=(p_*\mathcal G)(t)$, and (i) follows; the rank and the degree are read off from $\chi((p_*\mathcal G)(t))=\chi(\mathcal G(t))$ by Riemann--Roch on $\PP^n$. (ii) Since $p|_Z$ is finite, the associated points of $p_*\mathcal G$ are the images of those of $\mathcal G$, which are generic points of $n$-dimensional subvarieties; so the generic point of $\PP^n$ is the only associated point of $p_*\mathcal G$. (iii) For a closed point $y\in\PP^n$, the stalk $(p_*\mathcal G)_y$ is a finite module over the regular local ring $\OO_{\PP^n,y}$ of dimension $n$, and its depth is the minimum of the depths of the $\mathcal G_x$, $x\in p^{-1}(y)\cap Z$, that is $n$ (or the stalk is zero). By the Auslander--Buchsbaum formula it is free.
\end{proof}

\Needspace{6\baselineskip}
\begin{theorem}\label{thm:lowerdim}
Let $\cF$ be a coherent sheaf of dimension at most $n$ on $\PP^N$.
\begin{enumerate}[label=\textup{(\alph*)}]
\item For every coherent subsheaf $\cF'\subset\cF$ and every root sequence $z$ of length $n$,
\[
\deg_n\cF'-\mu(z)\rk_n\cF'\le\phi_z(\gamma(\cF)).
\]
Hence $\deg_n\cF'\le\Psi_\cF(\rk_n\cF')$, i.e.\ $P_\cF\le\Psi_\cF$. For $\cF$ pure of dimension $n$, the Harder--Narasimhan polygon of $\cF$ with respect to $\mu_n=-e_1/e_0$ lies below $\Psi_\cF$.
\item If $\cF_x$ is a Cohen--Macaulay module of dimension $n$ at every closed point $x$ of its support, for instance if $\cF$ is a vector bundle on a Cohen--Macaulay subvariety of dimension $n$, then $\gamma(\cF)$ is the table of a vector bundle on $\PP^n$, and $\Psi_\cF$ is the Boij--Söderberg polygon of its chain decomposition, with finitely many edges and endpoints $(0,0)$ and $(\rk_n\cF,\deg_n\cF)$.
\item $\Psi_\cF\le\Psi_{\cF'\oplus\cF''}$ for every exact sequence $0\to\cF'\to\cF\to\cF''\to0$. For $n=2$, Theorem~D holds for $\cF$ in the form of Proposition~\ref{prop:cohD}, with $\rk_2$ and $\deg_2$ in place of the rank and the degree.
\end{enumerate}
\end{theorem}

\begin{proof}
Choose $\Pi$ and $p$ as in Lemma~\ref{lem:projection}; a general $\Pi$ of dimension $N-n-1$ is disjoint from $\operatorname{Supp}\cF$. Every subsheaf and every quotient of $\cF$ is supported in $\operatorname{Supp}\cF$, so $p_*$ carries $0\to\cF'\to\cF\to\cF''\to0$ to an exact sequence on $\PP^n$, and by Lemma~\ref{lem:projection}(i) it preserves the tables, the long exact cohomology sequences, and the numbers $\rk_n$ and $\deg_n$, which become the rank and the degree on $\PP^n$. In particular $\Lambda_\cF$ and $\Psi_\cF$ are the functions of Section~\ref{sec:coherent} for the coherent sheaf $p_*\cF$ on $\PP^n$. (a) is Proposition~\ref{prop:cohA}(a) for $p_*\cF'\subset p_*\cF$. (b) By Lemma~\ref{lem:projection}(iii), $p_*\cF$ is a vector bundle on $\PP^n$, and by Proposition~\ref{prop:minformula} and Lemma~\ref{lem:legendre}(i) its function $\Psi$ is its Boij--Söderberg polygon; a vector bundle on a Cohen--Macaulay subvariety $X$ of dimension $n$ has stalks free over the Cohen--Macaulay local rings $\OO_{X,x}$ of dimension $n$. (c) is Proposition~\ref{prop:cohC}(a) and Proposition~\ref{prop:cohD} for $p_*$ of the sequence.
\end{proof}

\begin{remark}\label{rem:lowerdim}
\begin{enumerate}[label=(\arabic*)]
\item \emph{Vector bundles on subvarieties.} Let $X\subset\PP^N$ be a smooth subvariety of dimension $n$, $H=\OO_X(1)$, and $\cE$ a vector bundle on $X$, viewed as a sheaf on $\PP^N$. The subsheaves of $\cE$ on $\PP^N$ are the subsheaves on $X$, and by Riemann--Roch on $X$ a torsion-free subsheaf $\cE'$ of rank $r'$ has $\rk_n\cE'=r'H^n$ and $\deg_n\cE'=c_1(\cE')\cdot H^{n-1}-\frac{r'}2(K_X+(n+1)H)\cdot H^{n-1}$. So $P_\cE$ is the $\mu_H$-Harder--Narasimhan polygon of $\cE$ after a linear change of coordinates, and Theorem~\ref{thm:lowerdim} bounds it by the Boij--Söderberg polygon of the vector bundle $p_*\cE$ on $\PP^n$, whose table is $(h^j(\cE(tH)))_{j,t}$. If $H$ is ample and base-point-free, the proof applies to a finite morphism $\pi\colon X\to\PP^n$ with $\pi^*\OO(1)=\OO_X(H)$ in place of $p$; this is how the comparison is used in the proof of Corollary~\ref{cor:cases}(a).
\item \emph{Hilbert coefficients.} Let $M$ be a finitely generated graded module of Krull dimension $n+1\ge2$ over $S=K[x_0,\dots,x_N]$, with Hilbert polynomial $\sum_{i=0}^n(-1)^ie_i(M)\binom{t+n-i}{n-i}$, and let $M'\subset M$ be a graded submodule of the same dimension. Then $\widetilde{M'}\subset\widetilde M$, the Hilbert polynomials of $M'$ and $\widetilde{M'}$ agree, and Theorem~\ref{thm:lowerdim}(a) gives
\[
e_1(M')\ \ge\ -\Psi_{\widetilde M}\big(e_0(M')\big).
\]
The table of $\widetilde M$, hence the right-hand side, is determined by the Hilbert functions of $M$ and of its local cohomology modules $H^j_{\mathfrak m}(M)$, since $h^0(\widetilde M(t))=\dim M_t-\dim H^0_{\mathfrak m}(M)_t+\dim H^1_{\mathfrak m}(M)_t$ and $h^j(\widetilde M(t))=\dim H^{j+1}_{\mathfrak m}(M)_t$ for $j\ge1$. If $M$ is Cohen--Macaulay, then $M$ is free over a polynomial ring $K[y_0,\dots,y_n]$ generated by linear forms, with $h_k$ basis elements in degree $k$, where $\sum_kh_kt^k/(1-t)^{n+1}$ is the Hilbert series of $M$; so $p_*\widetilde M\cong\bigoplus_k\OO(-k)^{h_k}$ for the projection $p$ from $V(y_0,\dots,y_n)$, and the bound says that $e_1(M')$ is at least the sum of the $e_0(M')$ smallest elements of the multiset in which $k$ occurs $h_k$ times. In this case it reduces to the elementary bound for subsheaves of a split bundle.
\item \emph{The Betti side.} The multiplicity conjecture of Herzog, Huneke and Srinivasan \cite{HS98}, proved by Eisenbud and Schreyer \cite{ES09} and, for the upper bound in the non-Cohen--Macaulay case, by Boij and Söderberg \cite{BS12}, bounds $e_0(M)$ through the pure decomposition of the Betti table of $M$. Theorem~\ref{thm:lowerdim} is of a different kind: it bounds $e_0$ and $e_1$ of all submodules through the cohomology table. The naive analogue on the Betti side, with the number of generators in place of the rank, fails. For $N\ge2$ and $M=S/(x_0)\oplus S/(x_1)$, the Betti table is twice the pure diagram of type $(0,1)$, whose pieces have multiplicity $1$ per generator, while the cyclic submodule $M'=S\cdot(1,1)\cong S/(x_0x_1)$ has $e_0(M')=2$. Theorem~\ref{thm:lowerdim} holds in this example, with $n=N-1$: $e_0(M')=e_0(M)=2$ and $e_1(M')=1>0=e_1(M)$.
\end{enumerate}
\end{remark}

\section{Polarized varieties}\label{sec:polarized}

In this section $n\ge2$ is the dimension of a polarized variety, which we project to $\PP^n$.

\subsection{The defect identity}
Let $X$ be a smooth projective variety of dimension $n\ge2$ over $K$ and $H$ an ample, base-point-free divisor. Then $n+1$ general sections of $\OO_X(H)$ define a finite morphism $\pi\colon X\to\PP^n$ with $\pi^*\OO(1)=\OO_X(H)$, flat because $X$ is Cohen--Macaulay and $\PP^n$ is smooth. Put
\[
d=H^n,\qquad \tau=\frac{(K_X^2-2c_2(X))\cdot H^{n-2}}{n+1}.
\]
For $n=2$, Noether's formula gives $\tau=K_X^2-8\chi(\OO_X)$, the topological signature of $X$ over $\mathbb C$. For a vector bundle $\cE$ of rank $r$, $\pi_*\cE$ is a vector bundle of rank $N=rd$, and $H^j(X,\cE(tH))=H^j(\PP^n,\pi_*\cE(t))$. We set $G(\cE)=G(\pi_*\cE)$ and $V_c(\cE)=V_c(\pi_*\cE)$, with the gap moment and the variance of Section~\ref{sec:gap}. The bundle $\cE$ is aCM, i.e.\ $H^j(X,\cE(tH))=0$ for $0<j<n$ and all $t$, iff $\pi_*\cE$ splits, iff $G(\cE)=1$ (Lemma~\ref{lem:gapPn}(iii)).

For a torsion-free sheaf $\cQ$ of rank $m$ on $X$ put $v=c_1(\cQ)/m\in\operatorname{NS}(X)_\QQ$, $\mu_H(\cQ)=v\cdot H^{n-1}$ and
\[
\delta(\cQ)=\Big(v^2-\frac{2\operatorname{ch}_2(\cQ)}m\Big)\cdot H^{n-2}=\frac{\Delta(\cQ)\cdot H^{n-2}}{m^2},\qquad
\eta(\cQ)=\frac{\mu_H(\cQ)-\frac12K_X\cdot H^{n-1}}d,
\]
\[
\xi(\cQ)=v-\frac{K_X}2-\eta(\cQ)H ,
\]
so that $\xi(\cQ)\cdot H^{n-1}=0$. The number $\tau$, the class $\xi(\cQ)$ and the number $\eta(\cQ)$ are denoted $\sigma$, $\alpha(\cQ)$ and $\beta(\cQ)$ in \cite{An1,An2}. In the computations below every product of two classes is understood to be intersected with $H^{n-2}$, and a product of a class with $H$ to be intersected with $H^{n-2}$ as well; in particular $\xi(\cQ)^2$ denotes the number $\xi(\cQ)^2\cdot H^{n-2}$ and $v\cdot H$ the number $v\cdot H^{n-1}$. If $\cQ$ is $\mu_H$-semistable, then $\delta(\cQ)\ge0$ by Bogomolov's inequality \cite{Bo78}, \cite[Thm.~7.3.1]{HL10}. Always $\xi(\cQ)^2\le0$: choose $k_0\ge1$ with $k_0H$ very ample; by Bertini, $n-2$ general members of $|k_0H|$ cut out a smooth connected surface $S$, and $k_0^{n-2}\xi(\cQ)^2=(\xi(\cQ)|_S)^2\le0$ by the Hodge index theorem on $S$, since $\xi(\cQ)|_S\cdot H|_S=0$.

Let $0=\cE_0\subset\dots\subset\cE_s=\cE$ be the $\mu_H$-Harder--Narasimhan filtration of $\cE$, with factors $\cQ_i$ of ranks $r_i$ and slopes $\mu_i$. Put $w_i=r_i/r$, $\eta_i=\eta(\cQ_i)$, $\bar\eta=\sum_iw_i\eta_i$ and $V_\eta(\cE)=\sum_iw_i(\eta_i-\bar\eta)^2$.

\begin{lemma}\label{lem:push}
Let $\cQ$ be a torsion-free sheaf of rank $m$ on $X$.
\begin{enumerate}[label=\textup{(\alph*)}]
\item $\pi_*\cQ$ is a torsion-free sheaf of rank $md$ and slope $\eta(\cQ)-\frac{n+1}2$ on $\PP^n$.
\item $\displaystyle\frac{\chi(\cQ(tH))}m=d\,\frac{t^n}{n!}+d\,\eta(\cQ)\,\frac{t^{n-1}}{(n-1)!}+\Big(\frac d2\eta(\cQ)^2-\frac{(n+1)\,\tau}{24}+\frac{\xi(\cQ)^2-\delta(\cQ)}2\Big)\frac{t^{n-2}}{(n-2)!}+O(t^{n-3})$.
\end{enumerate}
For $n=2$ this is \cite[Lemma~2.1 and (2)]{An2}; in every dimension it is equivalent to \cite[Prop.~10.6]{An3}.
\end{lemma}

\begin{proof}
(b) By Riemann--Roch, $\chi(\cQ(tH))=\int_X\operatorname{ch}(\cQ)e^{tH}\operatorname{td}(X)$ with $\operatorname{td}(X)=1-\frac{K_X}2+\frac{K_X^2+c_2(X)}{12}+\dots$. The coefficient of $t^{n-1}/(n-1)!$ is $(c_1(\cQ)-\frac m2K_X)\cdot H^{n-1}=m\,d\,\eta(\cQ)$. The coefficient of $t^{n-2}/(n-2)!$ is $\big(\operatorname{ch}_2(\cQ)-\frac12c_1(\cQ)K_X+\frac m{12}(K_X^2+c_2(X))\big)\cdot H^{n-2}$; dividing by $m$ and substituting $\operatorname{ch}_2(\cQ)\cdot H^{n-2}=\frac m2(v^2-\delta(\cQ))$ and $v\cdot K_X=v^2+\frac{K_X^2}4-2\eta v\cdot H+\eta\,K_X\cdot H+\eta^2d-\xi(\cQ)^2$, with $v\cdot H=\eta d+\frac12K_X\cdot H$, gives the stated coefficient.

(a) Since $\pi$ is finite and $X$ is integral, $\pi_*\cQ$ is torsion-free, of rank $md$. As $\chi(\pi_*\cQ(t))=\chi(\cQ(tH))$, the coefficient of $t^{n-1}/(n-1)!$ computed on $\PP^n$, $\deg\pi_*\cQ+md\,\frac{n+1}2$, equals $md\,\eta(\cQ)$.
\end{proof}

\begin{proposition}[defect identity]\label{prop:defect}
For every vector bundle $\cE$ on $X$,
\[
\frac{n+1}{12}\big(d\,G(\cE)-\tau\big)=d\,\big(V_c(\cE)-V_\eta(\cE)\big)+\sum_iw_i\big(\delta(\cQ_i)-\xi(\cQ_i)^2\big).
\]
\end{proposition}

\begin{proof}
Summing Lemma~\ref{lem:push}(b) over the $\cQ_i$ with weights $w_i$ gives $\chi(\cE(tH))/r$. On the other hand $\chi(\cE(tH))/r=d\,\chi(\pi_*\cE(t))/N$, which is given by Lemma~\ref{lem:gapPn}(i). Comparing the coefficients of $t^{n-1}$ gives $\bar c=-\bar\eta$; comparing the coefficients of $t^{n-2}/(n-2)!$ and subtracting $\frac d2\bar c^{\,2}=\frac d2\bar\eta^2$ gives
\[
\frac d2V_c-\frac{(n+1)\,dG}{24}=\frac d2V_\eta-\frac{(n+1)\,\tau}{24}+\frac12\sum_iw_i(\xi_i^2-\delta_i). \qedhere
\]
\end{proof}

In terms of $\delta(\pi_*\cE)$, whose negative is the variance of the push-forward in \cite{An3}, the part of Proposition~\ref{prop:defect} that lives on $X$ is the identity $\frac{n+1}{12}(d-\tau)=-d\big(\delta(\pi_*\cE)+V_\eta\big)+\sum_iw_i(\delta_i-\xi_i^2)$ of \cite[(10.3)]{An3} (for $n=2$, \cite[Prop.~9.1]{An3}; for aCM bundles, \cite[Prop.~2.6]{An2}); Lemma~\ref{lem:gapPn}(ii), $\delta(\pi_*\cE)=\frac{n+1}{12}(G-1)-V_c$, converts it into Proposition~\ref{prop:defect}. What the table adds is the term $V_c$, which Theorem~\ref{thm:main} compares with $V_\eta$ for every bundle.

\begin{corollary}\label{cor:cases}
Let $\cE$ be a vector bundle of rank $r$ on $X$.
\begin{enumerate}[label=\textup{(\alph*)}]
\item $G(\cE)\ge\tau/d$, with equality iff $V_c(\cE)=V_\eta(\cE)$ and $\delta(\cQ_i)=\xi(\cQ_i)^2=0$ for all $i$. Moreover
\[
\sum_{0<j<n}\sum_th^j(\cE(tH))\ \ge\ C_n\,r\,(\tau-d)_+ ,
\]
with $C_n$ as in Lemma~\ref{lem:gapPn}(iv); for $n=2$ also $\sum_th^1(\cE(tH))\ge rd\,L(\kappa)$, where $\kappa=\max(1,\tau/d)$ and $L\colon[1,\infty)\to\RR$ is the function with $L(g^2)=\frac{g^3-g}{12}$ for integers $g\ge1$ that is linear on each interval $[g^2,(g+1)^2]$.
\item If $\cE$ is $\mu_H$-semistable, then $G(\cE)\ge\tau/d+\frac{12}{n+1}V_c(\cE)$.
\item \textup{(aCM, \cite[Thm.~1.1, Prop.~5.1]{An2}, \cite[\S10.7]{An3})} If $\cE$ is aCM, then $d\ge\tau$, i.e.\ $(n+1)H^n\ge(K_X^2-2c_2(X))\cdot H^{n-2}$.
\end{enumerate}
\end{corollary}

\begin{proof}
(a) The subsheaves $\pi_*\cE_i\subset\pi_*\cE$ form a filtration with torsion-free quotients $\pi_*\cQ_i$ of slopes $\eta_i-\frac{n+1}2$ (Lemma~\ref{lem:push}(a)), which strictly decrease. So its polygon $P_\eta$, with these slopes and lengths $dr_i$, is concave, has the same endpoints as $\Psi_{\pi_*\cE}$, and its vertices are the points $(\rk\pi_*\cE_i,\deg\pi_*\cE_i)$. By Theorem~\ref{thm:main} they lie below $\Psi_{\pi_*\cE}$, so $P_\eta\le\Psi_{\pi_*\cE}$ by Lemma~\ref{lem:legendre}(ii) (compare Remark~\ref{rem:lowerdim}(1)). As in the proof of Corollary~\ref{cor:supernatural}(c), Karamata's inequality gives $V_\eta(\cE)\le V_c(\cE)$. Both terms on the right of Proposition~\ref{prop:defect} are now nonnegative, and they vanish iff $V_c=V_\eta$ and $\delta(\cQ_i)=\xi(\cQ_i)^2=0$ for all $i$. The bound on the intermediate cohomology follows from Lemma~\ref{lem:gapPn}(iv) with $N=rd$, since $G(\cE)-1\ge(\tau-d)/d$; for $n=2$, Lemma~\ref{lem:basic}(iv) gives $\sum_th^1(\cE(tH))=\frac{rd}{12}\sum_kp_k(g_k^3-g_k)=rd\sum_kp_kL(g_k^2)$; the slopes of $L$ on the intervals $[g^2,(g+1)^2]$ are $\frac{g(g+1)}{4(2g+1)}$, which increase with $g$, so $L$ is convex and increasing, and Jensen's inequality gives $\sum_th^1(\cE(tH))\ge rd\,L(G(\cE))\ge rd\,L(\kappa)$.

(b) Here $V_\eta=0$, and the last sum in Proposition~\ref{prop:defect} is $\ge0$.

(c) If $\cE$ is aCM, then $G(\cE)=1$, and (a) gives $d\ge\tau$.
\end{proof}

\begin{remark}
\begin{enumerate}[label=(\arabic*)]
\item Corollary~\ref{cor:cases}(b) is the Boij--Söderberg form of Bogomolov's inequality. The point of (a) is uniformity: all cohomology tables of bundles on $(X,H)$ lie in the region $G\ge\tau/d$, a uniform lower bound for the gap moment, while the aCM bundles are those with $G=1$. This is not a bound for a distance between tables: on $\PP^2$ the normalized tables $(1-p)\gamma^{(-1,-2)}+p\gamma^{(m,0)}$ in $\CC$, of total rank $2$, with $a\in\QQ_{>0}$ fixed, $p=a/(m^2-1)$ and $m^2>a+1$, have $G=1+a$ for every $m$, while they converge entrywise to $\gamma^{(-1,-2)}=2\gamma(\OO)$ as $m\to\infty$. So for $\kappa>1$ the region $\{G\ge\kappa\}$ of the cone $\CC$ is not closed for entrywise convergence; these tables are not claimed to come from bundles on a fixed $X$. For the Hesse surfaces $Y_\nu$ with $H_\nu=4A_\nu-E_\nu$ of \cite[Cor.~1.4]{An2} (denoted $Y_n$ there), $\tau/d=(3\nu^2-11)/7$, which is $\frac{16}7$ for $\nu=3$ and $\frac{37}7$ for $\nu=4$. Since $L$ is the convex piecewise linear interpolation of $s\mapsto\frac{\sqrt s(s-1)}{12}$ at the squares, it dominates this function and $s\mapsto\frac{s-1}6$; for $\nu=4$ it gives $\sum_th^1(\cE(tH))\ge\frac{31}{35}rd$, against $\frac57rd$ from the linear bound.
\item If $S$ is a smooth surface cut out by $n-2$ general members of $|H|$, adjunction gives $K_S^2-2c_2(S)=(n+1)\tau-(n-2)d$, i.e.\ $3(\tau_S-d)=(n+1)(\tau-d)$; so $\tau>d$ holds for $(X,H)$ iff it holds for its surface sections; this is the case $a=1$ of the formula for the excess $(n+1)(\tau-d)$ under sections in $|aH|$ \cite[Rem.~5.7]{An2}, \cite[Lemma~3.3]{An3}, compare the proof of \cite[Lem.~3.1]{Lo26}. For $X=Y\times\PP^{n-2}$ with $H=H_Y\boxtimes\OO(1)$ one finds $(n+1)(\tau-d)=3\big(\tau_Y-(n-1)^2H_Y^2\big)$ \cite[Cor.~5.3]{An2}, \cite[Cor.~10.9]{An3}; for $Y=Y_\nu$ this is positive iff $(3\nu^2-11)/7>(n-1)^2$, e.g.\ for $\nu\ge4$ when $n=3$. Lopez \cite[Thm.~1]{Lo26} constructs from the Hesse surfaces varieties $Y_\nu\times Z$, $\nu\gg0$, of every dimension $n\ge3$, and their hyperplane sections, with $\tau>d$ for very ample $H$. On all of them (a) bounds the intermediate cohomology of every vector bundle from below.
\item The rank of $\pi_*\cE$ is $rd\ge d$, and $\pi_*\cQ_i$ need not be semistable, so $P_\eta$ is in general not the Harder--Narasimhan polygon of $\pi_*\cE$; the proof uses only that its vertices come from subsheaves. When the $\cQ_i$ are locally free, Corollary~\ref{cor:monoext}(a) gives the opposite bound $V_c(\cE)\le V_c(\gr\cE)$: cancellation can only decrease the spread of the centers, and Theorem~\ref{thm:main} says that it cannot decrease it below the spread $V_\eta$ of the slopes.
\item For $n=2$, Lemma~\ref{lem:twotwist} transports verbatim: $G(\cE)-1\ge\frac4{rd}\big(h^1(\cE(tH))+h^1(\cE((t+1)H))\big)$ for every $t$.
\item The proof of (c) in \cite[Prop.~5.1]{An2} restricts to a surface section; \cite[\S10.7]{An3} uses the projection to $\PP^n$ directly, with Horrocks' splitting and the majorization for split bundles \cite[Lemma~10.8]{An3}. Here Theorem~\ref{thm:main} replaces that majorization, so the identity of Proposition~\ref{prop:defect} measures the defect of every bundle, not only of the aCM ones.
\item Characteristic zero enters (a)--(c) only through Bogomolov's inequality, and it cannot be dropped. Lemmas~\ref{lem:gapPn} and~\ref{lem:push}, Proposition~\ref{prop:defect} and Theorem~\ref{thm:main} hold over any algebraically closed field. Over a field of characteristic $p>0$, every smooth surface with a very ample polarization carries Ulrich bundles \cite[Thm.~1.1]{Pi26a} (for $H$ ample see \cite[Lemma~2.11]{Pi26a}), and on a smooth variety with a very ample polarization a Frobenius direct image $F^e_*(\cE(aH))$ of any vector bundle $\cE$ is aCM for suitable $e$ and $a$ \cite[Prop.~4.11]{Pi26b}. The last statement holds for every ample $H$, by the same argument: by Serre vanishing and Serre duality there is $M$ with $H^j(\cE(tH))=0$ for $0<j<n$ and $|t|\ge M$; take $a=M$ and $p^e\ge2M$. Then $|a+p^et|\ge M$ for all $t\in\ZZ$, and by the projection formula $H^j\big(F^e_*(\cE(aH))(tH)\big)=H^j\big(\cE((a+p^et)H)\big)=0$ for $0<j<n$; and $F^e_*(\cE(aH))$ is a vector bundle because $X$ is smooth. These bundles have $G=1$ whatever $\tau/d$ is. For them Proposition~\ref{prop:defect} gives $\sum_iw_i\delta(\cQ_i)\le\frac{n+1}{12}(d-\tau)$, so if $\tau>d$ some Harder--Narasimhan factor violates Bogomolov's inequality and is therefore not strongly semistable \cite{La04}. For the Ulrich bundles of \cite{Pi26a} on a surface this is the computation of \cite[Rem.~2.9]{Pi26a}. Conversely, in characteristic $p$ the bound $G(\cE)\ge\tau/d$ of (a) holds for every $\cE$ whose Harder--Narasimhan factors are strongly semistable. The contrast also shows up in the cone question of the next subsection: if $\cU$ is Ulrich on $(X,H)$, then $\pi_*(\cU\otimes\pi^*S)\cong S^{\oplus\rk\pi_*\cU}$ for every vector bundle $S$ on $\PP^n$, so the tables $\gamma(\pi_*\cE)$ of vector bundles $\cE$ on $(X,H)$ generate the cone $\CC$ (see \cite[Thm.~5.3]{ESi}); for surfaces in characteristic $p$ this is \cite[Cor.~1.2]{Pi26a}. In every dimension, over an infinite perfect field of characteristic $p$ and for $H$ ample and base-point-free, the cones of cohomology tables of coherent sheaves on $X$ and on $\PP^n$ have the same closure under entrywise limits, by a theorem of Iyengar, Ma and Walker \cite[Thm.~3.4]{IMW23}. In characteristic zero with $\tau>d$, by (a) the tables $\gamma(\pi_*\cE)$ lie in the proper subcone $\{G\ge\tau/d\}$; by part~(1) of this remark, this region of $\CC$ is not closed for entrywise convergence.
\end{enumerate}
\end{remark}

\subsection{Vanishing at \texorpdfstring{$n$}{n} twists}\label{sec:twists}
After \cite{An1}, D.~Eisenbud and F.-O.~Schreyer asked the first author (personal communication) whether a surface $Y$ carries vector bundles $\cE$ with $h^j(\cE(a_1H))=h^j(\cE(a_2H))=0$ for $j=0,1,2$ at two non-consecutive integers $a_1<a_2$, the first case being $a_1=-1$, $a_2=1$, and what the Boij--Söderberg cone of cohomology tables of aCM bundles on $Y$ is. For consecutive integers such bundles are, up to twist, the Ulrich bundles; the aCM part is answered in \cite{An2}. For $H^2<\tau$ the question becomes one about all vector bundles, and \cite[\S12, Question~4]{An3} asks whether their cone is cut out from that of $\PP^2$ by conditions like $-\delta(\pi_*\cE)\le\frac14-\frac\tau{4d}$, which holds for $\mu_H$-semistable $\cE$, once unstable bundles are taken into account. Corollary~\ref{cor:cases}(a) gives such a condition for every bundle, $G\ge\tau/d$. Supernatural bundles with bounded roots on polarized varieties are studied in \cite{GL26}; what Corollary~\ref{cor:ES} adds is the explicit dependence on the normalized spread $\rho$ of the roots, through the defect identity. The structure of $\Psi$ gives a necessary condition for each configuration of $n$ twists.

\begin{corollary}\label{cor:ES}
Let $\cE\ne0$ be a vector bundle on $X$ and $a_1<\dots<a_n$ integers with $h^j(\cE(a_iH))=0$ for all $i$ and $j$. Then $\gamma(\pi_*\cE)$ is a multiple of $\gamma^{z}$, $z=(a_n,\dots,a_1)$, $\cE$ is $\mu_H$-semistable, and
\[
\frac{n+1}{12}\big(d\,\rho(z)-\tau\big)\ =\ \delta(\cE)-\xi(\cE)^2\ \ge\ 0 .
\]
In other words, a supernatural ray $\RR_{>0}\gamma^{z}$ contains the table of $\pi_*\cE$ for some vector bundle $\cE$ on $X$ only if $\rho(z)\ge\tau/d$. If the $a_i$ form an arithmetic progression of difference $e$, then $\rho(z)=e^2$ and the condition reads $e^2d\ge\tau$. In particular, if $\tau>4d$ there is no such bundle on a surface with $a_2-a_1=2$, e.g.\ with $a_1=-1$ and $a_2=1$. On the Hesse surfaces $Y_\nu$ of \cite[Cor.~1.4]{An2} there is none with $(a_2-a_1)^2<(3\nu^2-11)/7$; for instance none with $a_2-a_1=2$ once $\nu\ge4$.
\end{corollary}

\begin{proof}
The columns $a_1,\dots,a_n$ of $\gamma(\pi_*\cE)=\sum_kq_k\gamma^{z^k}$ vanish. Its entries are sums of nonnegative terms, so these columns vanish in every $\gamma^{z^k}$: the $a_i$ are roots of every $z^k$, i.e.\ $z^k=z$. By Corollary~\ref{cor:supernatural}(a) the bundle $\pi_*\cE$ is semistable. If $\cE$ were not $\mu_H$-semistable, its maximal destabilizing subsheaf $\cE_1$ would give $\pi_*\cE_1\subset\pi_*\cE$ of slope $\eta_1-\frac{n+1}2>\bar\eta-\frac{n+1}2=\mu(\pi_*\cE)$ (Lemma~\ref{lem:push}(a)), a contradiction. Now $G(\cE)=\rho(z)$, $V_c(\cE)=0$ and $V_\eta(\cE)=0$, and Proposition~\ref{prop:defect} gives the identity; $\delta(\cE)\ge0$ by Bogomolov's inequality and $\xi(\cE)^2\le0$ by the Hodge index theorem. A table of $\pi_*\cE$ on the ray of $\gamma^{z}$ is exactly one whose columns $a_1,\dots,a_n$ vanish. For an arithmetic progression, $\sum_i(a_i-c)^2=e^2\sum_i(i-\frac{n+1}2)^2=e^2\,\frac{n(n^2-1)}{12}$.
\end{proof}

Let $a_i=a_1+(i-1)e$. The bundle $\cU=\cE((a_n+e)H)$ has no cohomology at the twists $-eH,\dots,-neH$, so it satisfies the cohomological condition for Ulrich bundles on the polarized variety $(X,eH)$ (compare Remark~\ref{rem:main}(3)). When $eH$ is very ample, the semistability of $\cE$ also follows from \cite[Thm.~2.9]{CH12}, and the inequality $e^2d\ge\tau$, i.e.\ $(n+1)(eH)^n\ge(K_X^2-2c_2(X))\cdot(eH)^{n-2}$, from the inequality of Lopez for Ulrich bundles \cite[Lem.~3.1]{Lo26} applied to the polarization $eH$, whose proof is the computation of Proposition~\ref{prop:defect} in this case; see also \cite{An1} for surfaces. What Corollary~\ref{cor:ES} adds is the case of root sequences not in arithmetic progression, which occurs only for $n\ge3$ (two roots always form one) and does not reduce by a twist to the Ulrich condition for a multiple of $H$, and the interpretation of the inequality as a necessary numerical condition for the realization of a supernatural ray. For $e=1$ it gives $d\ge\tau$, which is the Ulrich case: then $\pi_*\cE\cong\OO(-a_n-1)^{\oplus rd}$, and $\cE((a_n+1)H)$ is an Ulrich bundle.

\appendix
\section{Supernatural bundles on \texorpdfstring{$\PP^3$}{P3} with small root windows}\label{app:P3}

By Corollary~\ref{cor:supernatural}, every vector bundle on $\PP^3$ whose table is a multiple of a supernatural table $\gamma^z$, $z=(z_1>z_2>z_3)$, is semistable. Two classes of such bundles are familiar. If $z_1-z_3\le3$, they are direct sums of copies of a twist of $\OO$, $\Omega^1$ or $\Omega^2$ (see below). If the roots form an arithmetic progression, a twist is an Ulrich bundle for a multiple of $\OO(1)$ (Remark~\ref{rem:main}(3)). The smallest root sequences outside both classes have $z_1-z_3=4$ and gaps $(z_1-z_2,z_2-z_3)=(3,1)$ or $(1,3)$. In this appendix we describe the bundles with $z_1-z_3=4$, up to a twist. For the gaps $(3,1)$ they are the duals of certain Steiner bundles (Lemma~\ref{lem:beil31}). For the gaps $(2,2)$ they are instanton bundles (Lemma~\ref{lem:beil22} and Remark~\ref{rem:beil22}), and a comparison with the semistability theorem of Jardim \cite{Ja06} shows that the vanishing at all three roots is needed. The consequences for stability are collected in Corollary~\ref{cor:appP3}.

The tool is Beilinson's theorem \cite{Be78}: a coherent sheaf $\cF$ on $\PP^n$ is the abutment of a spectral sequence with $E_1^{p,j}=H^j(\cF(p))\otimes\Omega^{-p}(-p)$, $-n\le p\le0$. If $z_1-z_3\le3$, the three roots lie in a window $\{k,k-1,k-2,k-3\}$, and after twisting by $\OO(k)$, which moves the roots into $\{0,-1,-2,-3\}$, only one column of the spectral sequence survives; so a vector bundle with table $q\gamma^z$ is a direct sum of copies of one of $\OO$, $\Omega^1(1)$, $\Omega^2(2)\cong T(-2)$, $\OO(-1)$, suitably twisted. If $z_1-z_3=4$, the gaps are $(2,2)$, $(3,1)$ or $(1,3)$, and the last two cases are dual to each other. We write $U=H^0(\OO_{\PP^3}(1))$, so that $\PP^3=\PP(U^\vee)$, $H^0(T(-1))=U^\vee$ and $\operatorname{Hom}(T(-2),\OO)=H^0(\Omega^1(2))=\Lambda^2U$, and we use that $T(-3)$ and $T(-5)$ are acyclic; here $T=T_{\PP^3}$.

\begin{lemma}[gaps $(2,2)$]\label{lem:beil22}
Let $z=(1,-1,-3)$ and $q\in\frac13\ZZ_{>0}$. A vector bundle $F$ on $\PP^3$ has cohomology table $q\gamma^z$ if and only if $F\cong\ker\varphi$ for a surjection $\varphi\colon V\otimes T(-2)\to W\otimes\OO$ with $\dim V=\dim W=3q$ such that
\[
B=H^0(\varphi(1))\colon V\otimes U^\vee\longrightarrow W\otimes U
\]
is an isomorphism. The surjectivity of $\varphi$ follows from that of $B$. For $q=\frac13$ these are the twists $N(-1)$ of the null correlation bundles.
\end{lemma}

\begin{proof}
If $\gamma(F)=q\gamma^z$, then $F(-1)$ and $F(-3)$ are acyclic, so the columns $p=-1,-3$ of the Beilinson spectral sequence vanish; the column $p=0$ is $H^1(F)\otimes\OO$ in the row $1$, and the column $p=-2$ is $H^2(F(-2))\otimes\Omega^2(2)$ in the row $2$, with $\Omega^2(2)\cong T(-2)$ on $\PP^3$. The only differential is $d_2\colon H^2(F(-2))\otimes T(-2)\to H^1(F)\otimes\OO$, of total degrees $0\to1$, and convergence to $F$ forces it to be surjective with kernel $F$. Here $h^2(F(-2))=h^1(F)=3q$. By the twisted sequence $0\to F(1)\to V\otimes T(-1)\to W\otimes\OO(1)\to0$ and the vanishing $H^i(T(-1))=H^i(\OO(1))=0$ for $i>0$, the acyclicity of $F(1)$ says exactly that $B$ is bijective. If $B$ is surjective, the image $I$ of $\varphi(1)$ has $H^0(I)=H^0(W\otimes\OO(1))$, and $W\otimes\OO(1)$ is globally generated, so $\varphi$ is surjective. Conversely, let $\cG=\ker\varphi$. It is locally free, and $\cG(t)$ is acyclic for $t=1$ by the hypothesis, and for $t=-1,-3$ because $T(-3)$, $T(-5)$, $\OO(-1)$, $\OO(-3)$ are. Every piece of the chain of $\cG$ has zero columns at $1,-1,-3$, hence root sequence $z$, so $\gamma(\cG)=q'\gamma^z$ with $6q'=\rk \cG=6q$. For $\dim V=1$, $\varphi$ is a nonzero $2$-form $\omega\in\Lambda^2U$, and $B$ is $x\mapsto\omega(\cdot,x)$, bijective iff $\omega$ is nondegenerate; then $\ker\varphi=N(-1)$ \cite[I, \S4.2]{OSS}.
\end{proof}

\begin{remark}\label{rem:beil22}
The bundle $E=F(1)$ is the cohomology of the linear monad $V\otimes\OO(-1)\to V\otimes U^\vee\otimes\OO\to W\otimes\OO(1)$, so it is an instanton bundle of rank $6q$ and charge $3q$ in the sense of \cite{Ja06}, with $E$, $E(-2)$ and $E(-4)$ acyclic; and $F(3)$ is an Ulrich bundle on $(\PP^3,\OO(2))$. For the correspondence between Ulrich bundles on $(\PP^3,\OO(d))$ and instanton bundles see \cite[Thm.~3.6]{CMR21}. Jardim proves that locally free instanton sheaves of rank at most $2n-1$ on $\PP^n$ are semistable, and that this bound is sharp \cite{Ja06}. A simple example on $\PP^3$ is $T(-1)\oplus\Omega^1(1)$, the cohomology of $\OO(-1)\to U^\vee\otimes\OO\oplus U\otimes\OO\to\OO(1)$: an instanton bundle of rank $6$ and charge $1$, destabilized by $T(-1)$, with $h^0=4$ in the column where our tables have a root. So the supernatural hypothesis in Corollary~\ref{cor:supernatural}(a) cannot be weakened to the vanishing conditions that define instanton bundles.
\end{remark}

\begin{lemma}[gaps $(3,1)$]\label{lem:beil31}
Let $q\in\frac12\ZZ_{>0}$. A vector bundle $F$ on $\PP^3$ has cohomology table $q\gamma^{(1,-2,-3)}$ if and only if $F\cong\ker\theta$ for a surjection $\theta\colon V\otimes\Omega^1(1)\to W\otimes\OO$ with $\dim V=4q$, $\dim W=6q$, such that $H^0(\theta(1))\colon V\otimes\Lambda^2U\to W\otimes U$ is bijective. Every map $\Omega^1(1)\to\OO$ is the restriction of a linear form on $U$, so $\theta$ is the restriction of a linear map $\tilde\theta\colon V\otimes U\to W$; put $S=\ker\tilde\theta$. Then $\dim S=10q$,
\[
0\longrightarrow F\longrightarrow S\otimes\OO\longrightarrow V\otimes\OO(1)\longrightarrow0
\]
is exact, the second map being the evaluation, and the condition on $\theta$ says that $S\otimes U\to V\otimes S^2U$ is bijective. So $F^\vee$ is the Steiner bundle $\coker\big(V^\vee\otimes\OO(-1)\to S^\vee\otimes\OO\big)$, of rank $6q$; Steiner bundles of rank $n$ on $\PP^n$ were introduced in \cite{DK93}. For $q=\frac12$ an example is the Hankel space $S=\langle e_1\otimes u_{j}+e_2\otimes u_{j-1}:0\le j\le4\rangle$, for a basis $e_1,e_2$ of $V$ and a basis $u_0,\dots,u_3$ of $U$, with $u_{-1}=u_4=0$.
\end{lemma}

\begin{proof}
If $\gamma(F)=q\gamma^z$, $z=(1,-2,-3)$, then $F(-2)$ and $F(-3)$ are acyclic, and the cohomology of $F$ and $F(-1)$ is concentrated in degree one, with $h^1(F)=6q$ and $h^1(F(-1))=4q$. So the only terms of the Beilinson spectral sequence are $E_1^{-1,1}=H^1(F(-1))\otimes\Omega^1(1)$ and $E_1^{0,1}=H^1(F)\otimes\OO$, and as in Lemma~\ref{lem:beil22} the differential is a surjection $\theta$ with kernel $F$. The sequence $0\to F(1)\to V\otimes\Omega^1(2)\to W\otimes\OO(1)\to0$ and the vanishing of $H^i(\Omega^1(2))$ and $H^i(\OO(1))$ for $i>0$ show that $F(1)$ is acyclic if and only if $H^0(\theta(1))$ is bijective. The assertion on maps $\Omega^1(1)\to\OO$ follows from the Euler sequence $0\to\Omega^1(1)\to U\otimes\OO\to\OO(1)\to0$, since $\operatorname{Hom}(\OO(1),\OO)=\Ext^1(\OO(1),\OO)=0$. Since $\theta$ is surjective, so is $\tilde\theta$, hence $\dim S=10q$, and the map $V\otimes U\otimes\OO\to W\otimes\OO\oplus V\otimes\OO(1)$ given by $\tilde\theta$ and the evaluation is surjective, by the snake lemma applied to the Euler sequence tensored with $V$; restricting it to $S\otimes\OO$ gives the displayed sequence. Its cohomology at the twist $1$ is $0\to H^0(F(1))\to S\otimes U\to V\otimes S^2U\to H^1(F(1))\to0$, which gives the second form of the condition. Conversely, if $\theta$ is a surjection with $H^0(\theta(1))$ bijective, then $F(1)$ is acyclic, and $F(-2)$, $F(-3)$ are acyclic because $\OO(-1),\OO(-2),\OO(-3)$ are; as in Lemma~\ref{lem:beil22}, every piece of the chain of $F$ has the roots $1,-2,-3$, so $\gamma(F)=q'\gamma^z$ with $6q'=\rk F=6q$. Finally, the condition on $S$ alone suffices: if $S\otimes U\to V\otimes S^2U$ is surjective, then, as in Lemma~\ref{lem:beil22}, the image of $S\otimes\OO(1)\to V\otimes\OO(2)$ contains all sections of the globally generated sheaf $V\otimes\OO(2)$, so the evaluation $S\otimes\OO\to V\otimes\OO(1)$ is surjective; with $W=(V\otimes U)/S$, the Euler sequence identifies the cokernel and the kernel of the induced map $\theta\colon V\otimes\Omega^1(1)\to W\otimes\OO$ with those of the evaluation on $S\otimes\OO$, so $\theta$ is surjective with kernel $F$. For the Hankel space, the $2\times5$ matrix of linear forms of the evaluation map has rows $(u_0,u_1,u_2,u_3,0)$ and $(0,u_0,u_1,u_2,u_3)$, which are proportional at no point of $\PP^3$, so the evaluation is surjective; and $S\otimes U\to V\otimes S^2U$, a map between spaces of dimension $20$, is injective over every field. Indeed, with $s_j=e_1\otimes u_j+e_2\otimes u_{j-1}$, an element $\sum_{j=0}^4s_j\otimes\ell_j$, $\ell_j=\sum_ka_{jk}u_k$, maps to $e_1\otimes\sum_{j=0}^3u_j\ell_j+e_2\otimes\sum_{j=0}^3u_j\ell_{j+1}$. The kernel of the multiplication $U\otimes U\to S^2U$ consists of the alternating tensors, so if the image is zero, the matrices $(a_{jk})_{0\le j,k\le3}$ and $(a_{j+1,k})_{0\le j,k\le3}$ are alternating. Then $a_{j+1,k}=-a_{k+1,j}=a_{j,k+1}$ for $0\le j\le3$, $0\le k\le2$, so $a_{jk}=h(j+k)$ depends only on $j+k\in[0,7]$, and the vanishing of the diagonals gives $h(2j)=a_{jj}=0$ and $h(2j+1)=a_{j+1,j}=0$ for $0\le j\le3$. So all $\ell_j$ vanish.
\end{proof}

\begin{corollary}\label{cor:appP3}
\begin{enumerate}[label=\textup{(\alph*)}]
\item For every $q\in\frac13\ZZ_{>0}$, the bundles $F$ of Lemma~\ref{lem:beil22} are semistable. So the instanton bundles $E=F(1)$ of Remark~\ref{rem:beil22}, of rank $6q$ and charge $3q$, are semistable also when their rank exceeds the bound $2n-1=5$ of \cite{Ja06}.
\item For every $q\in\frac12\ZZ_{>0}$, the bundles $F$ of Lemma~\ref{lem:beil31} and their duals, the Steiner bundles $\coker\big(V^\vee\otimes\OO(-1)\to S^\vee\otimes\OO\big)$ of rank $6q$ for which $S\otimes U\to V\otimes S^2U$ is bijective, are semistable. For $q=\frac12$ they are stable, and in characteristic zero the general one is stable for every $q$.
\end{enumerate}
\end{corollary}

\begin{proof}
The semistability is Corollary~\ref{cor:supernatural}(a), in any characteristic; in (a) it also follows from \cite[Thm.~2.9]{CH12}, since $F(3)$ is an Ulrich bundle on $(\PP^3,\OO(2))$. In (b), for $q=\frac12$, $\rk F=3$ and $c_1(F)=-2$ are coprime, so every such $F$ is stable. The subspaces $S\subset V\otimes U$ of dimension $10q$ for which $S\otimes U\to V\otimes S^2U$ is bijective form a nonempty open subset of the Grassmannian; it contains the direct sums of $2q$ Hankel spaces. A choice of basis of $S$ identifies the frame bundle over this open set with a nonempty open subset of the space of $10q\times4q$ matrices of linear forms, so a general point of it is a general matrix of linear forms in the sense of Coskun, Huizenga and Smith, whose theorem \cite[Thm.~5.4]{CHS} (in characteristic zero) applies to the Steiner bundle $F^\vee$ of rank $r=6q\ge3$ with $t=4q$ relations, since $r<(1+\sqrt3)\,t$. So the general vector bundle with table $q\gamma^{(1,-2,-3)}$ is stable, and dually for the gaps $(1,3)$. For $q\ge1$, direct sums give strictly semistable ones.
\end{proof}

\section*{Acknowledgments}
This paper grew out of an attempt to answer the questions on the Boij--Söderberg cone of cohomology tables of vector bundles on the surfaces of \cite{An1} that David Eisenbud and Frank-Olaf Schreyer asked the first author after that paper appeared (see Section~\ref{sec:twists}). We are grateful for their interest and encouragement.

\section*{AI Declaration}
ChatGPT and Claude were used for exploratory computations, literature searches, and interactive feedback on earlier drafts of the manuscript. The authors take full responsibility for all content.

\end{document}